\documentclass[nonblindrev]{informs3}

\OneAndAHalfSpacedXI
\usepackage{mathtools}
\usepackage{natbib}
 \bibpunct[, ]{(}{)}{,}{a}{}{,}%
 \def\bibfont{\small}%

\TheoremsNumberedThrough     %
\EquationsNumberedThrough    %

\usepackage{caption,subcaption}
\usepackage{graphicx}
\usepackage{enumitem}
\usepackage{bbm}
\usepackage{palatino}
\usepackage[para,online,flushleft]{threeparttable}
\usepackage{booktabs,multirow,makecell,tabularx}
\usepackage{longtable}
\usepackage[para, referable, flushleft]{threeparttablex}
\usepackage{adjustbox}

\usepackage{tikz}
\usepackage{pgfplots}
\pgfplotsset{compat=1.10}
\usetikzlibrary{positioning}

\definecolor{darkgreen}{rgb}{0.31, 0.47, 0.26}

\newcolumntype{L}[1]{>{\raggedright\let\newline\\arraybackslash\hspace{0pt}}m{#1}}
\newcolumntype{C}[1]{>{\centering\let\newline\\arraybackslash\hspace{0pt}}m{#1}}
\newcolumntype{R}[1]{>{\raggedleft\let\newline\\arraybackslash\hspace{0pt}}m{#1}}

\newlist{inlinelist}{enumerate*}{1}
\setlist[inlinelist]{label=(\roman*)} %

\makeatletter
\newcommand\footnoteref[1]%
\makeatother

\def\Indc{{\mathbbm{1}}}
\def\Expect{{\mathbb E}}
\def\Prob{{\mathbb P}}

\def\subto{{\rm s.\rm t.}}
\def\min{\mathop{\rm min}}
\def\max{\mathop{\rm max}}

\def\argmin{\mathop{\rm argmin}}
\def\sup{\mathop{\rm sup}}
\def\inf{\mathop{\rm inf}}

\def\Indc{{\mathbbm{1}}}
\renewcommand{\Re}{\mathbb{R}}
\newcommand{\X}{\mathcal{X}}
\newcommand{\Y}{\mathcal{Y}}
\newcommand{\Z}{\mathcal{Z}}

\newcommand{\Dist}{\mathcal{D}}
\newcommand{\Dirac}{\delta}

\newcommand{\mymbox}[1]{\mbox{\scriptsize #1}}
\newcommand{\quoteIt}[1]{``#1''}

\newcommand{\z}{\boldsymbol{z}}
\newcommand{\x}{\boldsymbol{x}}
\newcommand{\y}{\boldsymbol{y}}
\newcommand{\bv}{\boldsymbol{v}}

\newcommand{\V}[1]{\boldsymbol{#1}}
\renewcommand{\Xi}{\Y}

\newcommand{\Vxi}{\boldsymbol{y}}

\newcommand{\btheta}{{\boldsymbol{\theta}}}

\newcommand{\Pemp}{\hat{\mathbb{P}}_N}

\usepackage{pifont}
\newcommand{\cross}{\ding{55}}
\newcommand{\tick}{\ding{52}}

\renewcommand{\bibfont}{\small}
\usepackage[colorlinks=true,linkcolor=blue!80!black,citecolor=green!50!black,urlcolor=magenta]{hyperref}

\newcommand{\model}[1]{\textbf{\texttt{#1}}}

\let\oldcross\cross
\renewcommand{\cross}{\textcolor{gray}{\oldcross}}

\usepackage{etoolbox}
\makeatletter
\pretocmd{\NAT@citex}{%
  \let\NAT@hyper@\NAT@hyper@citex
  \def\NAT@postnote{#2}%
  \setcounter{NAT@total@cites}{0}%
  \setcounter{NAT@count@cites}{0}%
    \forcsvlist{\stepcounter{NAT@total@cites}\@gobble}{#3}}{}{}
\newcounter{NAT@total@cites}
\newcounter{NAT@count@cites}
\def\NAT@postnote{}
\def\NAT@hyper@citex#1{%
  \stepcounter{NAT@count@cites}%
  \hyper@natlinkstart{\@citeb\@extra@b@citeb}#1%
  \ifnumequal{\value{NAT@count@cites}}{\value{NAT@total@cites}}
    {\ifNAT@swa\else\if*\NAT@postnote*\else%
     \NAT@cmt\NAT@postnote\global\def\NAT@postnote{}\fi\fi}{}%
  \ifNAT@swa\else\if\relax\NAT@date\relax
  \else\NAT@@close\global\let\NAT@nm\@empty\fi\fi%
  \hyper@natlinkend}
\renewcommand\hyper@natlinkbreak[2]{#1}
\patchcmd{\NAT@citex}
  {\ifNAT@swa\else\if*#2*\else\NAT@cmt#2\fi
   \if\relax\NAT@date\relax\else\NAT@@close\fi\fi}{}{}{}
\patchcmd{\NAT@citex}
  {\if\relax\NAT@date\relax\NAT@def@citea\else\NAT@def@citea@close\fi}
  {\if\relax\NAT@date\relax\NAT@def@citea\else\NAT@def@citea@space\fi}{}{}
\makeatother

\newcommand{\revised}[1]{#1}

\begin{document}

\TITLE{A Survey of Contextual Optimization Methods for Decision-Making under Uncertainty}

\RUNTITLE{A Survey of Contextual Optimization Methods for Decision-Making under Uncertainty}

\RUNAUTHOR{Sadana et al.}
\ARTICLEAUTHORS{\AUTHOR{Utsav Sadana}
                \AFF{Department of Computer Science and Operations Research\\ Université de Montréal,  Québec, Canada\\ \EMAIL{utsav.sadana@umontreal.ca}}
                \AUTHOR{Abhilash Chenreddy}
                \AFF{GERAD \& Department of Decision Sciences\\ HEC Montréal, Québec, Canada\\ \EMAIL{abhilash.chenreddy@hec.ca}}
                \AUTHOR{Erick Delage}
                \AFF{GERAD \& Department of Decision Sciences\\ HEC Montréal, Québec, Canada\\ \EMAIL{erick.delage@hec.ca}}
                \AUTHOR{Alexandre Forel}
                \AFF{CIRRELT \& SCALE-AI Chair in Data-Driven Supply Chains, \\Department of Mathematical and Industrial Engineering, Polytechnique Montréal, Québec, Canada\\ \EMAIL{alexandre.forel@polymtl.ca}}
                \AUTHOR{Emma Frejinger}
                \AFF{CIRRELT \& Department of Computer Science and Operations Research
\\ Université de Montréal, Québec, Canada \\\EMAIL{emma.frejinger@umontreal.ca}}
                \AUTHOR{Thibaut Vidal}
                \AFF{CIRRELT \& SCALE-AI Chair in Data-Driven Supply Chains, \\Department of Mathematical and Industrial Engineering, Polytechnique Montréal, Québec, Canada\\ \EMAIL{thibaut.vidal@polymtl.ca}}

}
\ABSTRACT{%
Recently there has been a surge of interest in operations research~(OR) and the machine learning~(ML) community in combining prediction algorithms and optimization techniques to solve decision-making problems in the face of uncertainty. This gave rise to the field of contextual optimization, under which data-driven procedures are developed to prescribe actions to the decision-maker that make the best use of the most recently updated information. A large variety of models and methods have been presented in both OR and ML literature under a variety of names, including data-driven optimization, prescriptive optimization, predictive stochastic programming, policy optimization, (smart) predict/estimate-then-optimize, decision-focused learning, (task-based) end-to-end learning/forecasting/optimization, etc. \revised{This survey} article \revised{unifies these models under the lens of contextual stochastic optimization, thus providing a general presentation of a large variety of problems.} We identify three main frameworks for learning policies from data and present the existing models and methods under a uniform notation and terminology. Our objective with this survey is to both strengthen the general understanding of this active field of research and stimulate further theoretical and algorithmic advancements in integrating ML and stochastic programming.
}%

\maketitle
\tableofcontents
\section{Introduction}
This article surveys the literature on single and two-stage contextual optimization. In contextual optimization, a decision-maker faces a decision-making problem with uncertainty where the distribution of uncertain parameters that affect the objective and the constraints is unknown, although correlated side information (covariates or features) can be exploited. 
The usefulness of side information in inferring relevant statistics of uncertain parameters and, thereby, in decision-making is evident in many different fields. For example, weather and time of day can help resolve uncertainty about congestion on a road network and aid in finding the shortest path traversing a city. In portfolio optimization, stock returns may depend on historical prices and sentiments posted on Twitter \citep{xu-cohen-2018-stock}. Harnessing this information can allow decision-makers to build a less risky portfolio. Similarly, a retailer facing uncertain demand for summer products can infer whether the demand will be low or high depending on the forecasted weather conditions \citep{martinez2021comes}.

\revised{In these applications, the decision-maker has access to historical data, that is, past values of the covariates (e.g., weather) and the corresponding uncertain parameter (e.g., congestion). Data-driven contextual optimization methods use this data to estimate the conditional distribution of the uncertain parameter (or a sufficient statistic) based on the covariate. Conversely, t}raditional stochastic optimization models ignore contextual information and use unconditional distributions of the uncertain parameters to make a decision \citep{birge2011introduction}. Such a decision may be suboptimal \citep{ban_big_2019} and, in some cases, even at the risk of being infeasible \citep{Rahimian2022}. The availability of data and huge computational power combined with advancements in machine learning (ML) and optimization techniques have resulted in a shift of paradigm to contextual optimization \citep{mivsic2020data}.

Making prescriptions using the side information requires a decision rule that maps the observed \revised{covariate} to an action. We identify three different paradigms for learning this mapping.
\begin{itemize}
    \item \textbf{Decision rule optimization}: This approach was introduced to the operation research community in \cite{Liyanage_Shanthikumar_2005} for data-driven optimization and popularized in 
    \cite{ban_big_2019} for big data environments, although a similar idea was already common practice in reinforcement learning under the name of policy gradient methods \citep[see][and literature that followed]{Sutton1999}. It consists in employing a parameterized mapping as the decision rule and in identifying the parameter that achieves the best empirical performance based on the available data. The decision rule can be formed as a linear combination of functions of the covariates or even using a deep neural network~(DNN). When the data available is limited, some form of regularization might also be needed.
    \item \textbf{Sequential learning and optimization~(SLO)}: \cite{Bertsimas_Kallus_2020} appears to be the first to have formalized this two-stage procedure (also referred to as predict/estimate-then-optimize or prescriptive optimization/stochastic programming) that first uses a trained model to predict a conditional distribution for the uncertain parameters given the covariates, and then solves an associated \revised{contextual} stochastic optimization~(CSO) problem to obtain the optimal action. This procedure can be robustified to reduce post-decision disappointment \citep{smith2006optimizers} caused by model overfitting or misspecification by using proper regularization at training time or by adapting the CSO problem formulation.
    \item \textbf{Integrated learning and optimization~(ILO)}: In the spirit of identifying the best decision rule, one might question in SLO the need for high precision predictors when one is instead mostly interested in the quality of the resulting prescribed action. This idea motivates an integrated version of learning and optimization that searches for the predictive model that guides the CSO problem toward the best-performing actions. The ILO paradigm appears as early as in \cite{Bengio_1997} and has seen a resurgence recently in active streams of literature \revised{under various names such as smart predict-then-optimize, decision-focused learning, and (task-based) end-to-end learning/forecasting/optimization.}
\end{itemize}

The outline of the survey goes as follows. Section \ref{sec:notations} rigorously defines the three frameworks for \revised{identifying the best mapping from covariate to action based on data}: decision rule optimization, SLO, and ILO. Section~\ref{sec:decision_rules} reviews the literature on decision rule optimization with linear and non-linear decision rules. Section~\ref{sec:sequential} focuses on SLO, including the models that lead to robust decisions, and Section~\ref{sec:integrated_learning} describes the models based on the ILO framework and the algorithms used to train them. \revised{Because ILO is the more recent and less explored framework of the three identified, we provide a separate subsection of applications of ILO to diverse problems such as logistics and energy management.} Section~\ref{sec:directions} provides an overview of active research directions being pursued both from a theoretical and applications perspective. \revised{Section~\ref{sec:conclusion} concludes our survey with a summary of our contributions.}

We note that there are other surveys and tutorials in the literature that are complementary to ours. \citet{mivsic2020data} survey the applications of the SLO framework to problems in supply chain management, revenue management, and healthcare operations. \cite{qi2022integrating} is a tutorial that mainly focuses on the application of ILO to expected value-based models with limited discussions on more general approaches. It summarizes the most popular methods and some of their theoretical guarantees. \cite{Kotary_Fioretto_Van_Hentenryck_Wilder_2021} provide a comprehensive survey of the literature proposing ML methods to accelerate the resolution of constrained optimization models \citep[see also][]{Bengio_2021}. It also reviews some of the earlier literature on ILO applied to \revised{what we define as ``expected value-based models'', a subset of CSO problems (see Definition~\ref{def:exp_model}) where uncertainty can be completely described by a sufficient statistic before the optimization problem is solved. \citet{mandi2023decision}\footnote{\revised{\cite{mandi2023decision} appeared online soon after the submission of our paper.}} include more recent expected value-based models and provide a comprehensive evaluation of their methods, complementing the toolbox of \citet{Tang_Khalil_2022} that provided an interface for solving expected value-based models.}

Our survey of ILO \revised{literature} goes beyond the expected value-based models and reflects better the more modern literature by casting the contextual decision problem as a CSO problem and presenting a comprehensive overview of the current state of this rapidly progressing field of research. We establish links between approaches that minimize regret \citep{elmachtoub_smart_2022}, (task-based) end-to-end learning \citep{donti_task-based_2019} and imitation-based models \citep{Kong2022energy}. Further, we create a taxonomy based on the training procedure for a general ILO framework encompassing recent theoretical and algorithmic progresses in designing differentiable surrogates and optimizers and improving training procedures based on unrolling and implicit differentiation.

\section{Contextual optimization: An overview} \label{sec:notations}
The contextual optimization paradigm seeks a decision (i.e., an action) $\z$ in a feasible set $\Z \subseteq \Re^{\revised{d_{\z}}}$ that minimizes a cost function $c(\z, \y )$ with uncertain parameters $\y \in \Y \subseteq \Re^{\revised{d_{\y}}}$. The uncertain parameters are unknown when making the decision. However, a vector of relevant \revised{covariates}~\mbox{$\x \in \mathcal{X} \subseteq \Re^{\revised{d_{\x}}}$}, which is correlated with the uncertain parameters $\y$, is revealed before having to choose~$\z$. The joint distribution of the \revised{covariates} in $\X$ and uncertain parameters in $\Y$ is denoted by $\Prob$.

\subsection{Contextual problem and policy}
In general, uncertainty can appear in the objectives and constraints of the problem. In the main sections of this paper, we focus on problems with uncertain objectives and consider that the decision-maker is risk-neutral. We broaden the discussion to risk-averse settings and uncertain constraints in Section~\ref{sec:directions}.

\paragraph{The CSO problem.} Given a \revised{covariate} described by a vector of \revised{covariates} $\x$ and the joint distribution $\Prob$ of the \revised{covariates} $\x$ and uncertain parameter $\y$, a risk-neutral decision-maker is interested in finding an optimal action $z^*(\x) \in \Z$ that minimizes the expected costs conditioned on the covariate $\x$. Formally, the optimal action is a solution to the \revised{CSO} problem given by:
\begin{align}
    \label{eq:CSO}
    \text{(CSO)} \quad &   z^*(\x) \in \argmin_{\z \in \Z}\Expect_{\Prob(\y\rvert \x)}\left[c \left(\z, \y \right) \right],
\end{align}
where $\Prob(\y \rvert \x)$ denotes the conditional distribution of $\y$ given the covariate $\x$ and it is assumed that a minimizer exists. For instance, a minimizer exists when $\Z$ is compact, $\Prob(\y\rvert \x)$ has bounded support and $c(\z,\y)$ is continuous in $\z$ almost surely \citep[see][for more details]{vanparys2021data}.

Problem~\eqref{eq:CSO} can equivalently be written in a compact form using the expected cost operator $h(\cdot, \cdot)$ that receives a\revised{n action} as a first argument and a distribution as a second argument:
\begin{align}
    \label{eq:CSO_h}
    z^*(\x) \in \argmin_{\z \in \Z} h(\z, \Prob(\y \rvert \x)) :=\Expect_{\Prob(\y \rvert \x)}\left[c \left(\z, \y \right) \right].
\end{align}

\paragraph{Optimal policy.} In general, the decision-maker repeatedly solves CSO problems in many different contexts. Hence, the decision-maker is interested in finding the policy that provides the lowest long-term expected cost, that is:
\begin{equation}
    \label{eq:CSO:policy}
   \pi^* \in \argmin_{\pi \in \Pi}\Expect_{\Prob} \big[c(\pi(\x), \y) \big] = \argmin_{\pi \in \Pi}\Expect_{\Prob} \big[h(\pi(\x), \Prob(\y \rvert \x))\big],
\end{equation}
\revised{where $\Pi:=\{\pi:\X\rightarrow\Z\}$ denotes the class of all feasible policies.}

Note that the optimal policy does not need to be obtained explicitly in a closed form. Indeed, based on the interchangeability property \citep[see Theorem 14.60 of][]{rockafellar1998variational}, solving the CSO problem~\eqref{eq:CSO} \revised{for} any \revised{covariate} $\x$ naturally identifies an optimal policy:
\begin{equation*}
    \bar{\pi}(\x)\in\argmin_{\z \in \Z} h(\z, \Prob(\y \rvert \x)) \; \text{a.s.} \Leftrightarrow \Expect_{\Prob} \big[h(\bar{\pi}(\x), \Prob(\y \rvert \x))\big] = \min_{\pi \in \Pi}\Expect_{\Prob} \big[h(\pi(\x), \Prob(\y \rvert \x))\big],
\end{equation*}
assuming that a minimizer of $h(\z, \Prob(\y \rvert \x))$ exists almost surely. Thus, the two optimal policies $\pi^*$ and $z^*(\cdot)$ coincide.

\subsection{Mapping \revised{covariate} to \revised{actions} in a data-driven environment}
Unfortunately, the joint distribution $\Prob$ is generally unknown. Instead, the decision-maker possesses historical data $\Dist_N = \left\lbrace (\x_i, \y_i) \right\rbrace_{i=1}^N$ that is assumed to be made of independent and identically distributed realizations of $(\x,\y)\in \mathcal{X} \times \mathcal{Y}$. Using this data, the decision-maker aims to find a policy that approximates well the optimal policy given by \eqref{eq:CSO:policy}. Many approaches have been proposed to find effective approximate policies. Most of them can be classified into the three main frameworks that we introduce below:
\begin{inlinelist}
    \item decision rule optimization,
    \item sequential learning and optimization, and
    \item integrated learning and optimization.
\end{inlinelist}

\subsubsection{Decision rule optimization.} In this framework, the policy is assumed to belong to a hypothesis class $\Pi^{\btheta}:=\{\pi_{\btheta}\}_{\btheta\in\Theta}\subseteq \Pi$ that contains a family of parametric policies $\pi_\btheta:\mathcal{X} \rightarrow\mathcal{Z}$ (e.g., linear functions or decision trees). The parameterized policy $\pi_\btheta$ maps directly any \revised{covariate} $\x$ to an \revised{action} $\pi_\btheta(\x)$ and will be referred to as a decision rule.

Denote by $\Pemp$ the empirical distribution of $(\x,\y)$ given historical data $\Dist_N$. One can identify the \quoteIt{best} parameterization of the policy in $\Pi^\btheta$ by solving the following empirical risk minimization~(ERM) problem:
\begin{align}
    \label{eq:policy_opt}
   \text{(ERM)} \quad & \btheta^*\in\argmin_{\btheta} H(\pi_{\btheta},\Pemp) := \Expect_{\Pemp} \big[c(\pi_\btheta(\x),\y)\big].
\end{align} 

\begin{figure}[htb]
    \centering
    \begin{tikzpicture}[scale=0.9, every node/.style={scale=0.9}]
    \node (A) at (0, 0) {};
    \node [ draw,
            minimum width=2cm,
            minimum height=1.2cm,
            right=2cm of A,
            align=center
           ]  (ml) {Decision\\ rule};
    
    \draw[-stealth] (A) -- node[label= above:Context] {} node[label= below:$\x$] {} (ml.west);
    \draw[-stealth] (ml.east) -- node[label= above:Decision] {} node[label= below:$\pi_\btheta(\x)$] {} ++ (2,0) node (B) {};
    
    \node (A) [right=1cm of B] {};
    \node [ draw,
            minimum width=2cm,
            minimum height=1.2cm,
            right=2cm of A,
            align=center
           ]  (ml) {Decision\\ rule};
    
    \node [draw,%
            minimum width=2cm, 
            minimum height=1.2cm,
            right=2cm of ml,
            align=center
          ] (loss) {Task loss};
          
    \draw[-stealth] (A) -- node[label= above:Context] {} node[label= below:$\x_i$] {} (ml.west);
    \draw[-stealth] (ml.east) -- node[label= above:Decision] {} node[label= below:$\pi_\btheta(\x_i)$] {} (loss.west);
    
    \draw[-stealth, thick, dashed] (loss.south) -- ++(0,-1cm)  -| node [below,pos=0.25] {$\nabla_\theta c(\pi_\btheta(\x_i), \V{y}_i)$}(ml.south);
    
    \draw[-stealth] (A) ++(0,1.7cm) -| (loss.north) node[label= above:Parameter, pos=0.1] {} node[label= below:$\V{y}_i$, pos=0.1] {};
    
\end{tikzpicture}
    \caption{Decision and training pipelines based on the decision rule paradigm: (left) the decision pipeline and (right) the training pipeline for a given training example $(\x_i, \y_i)$.}
    \label{fig:decision_rule}
\end{figure}
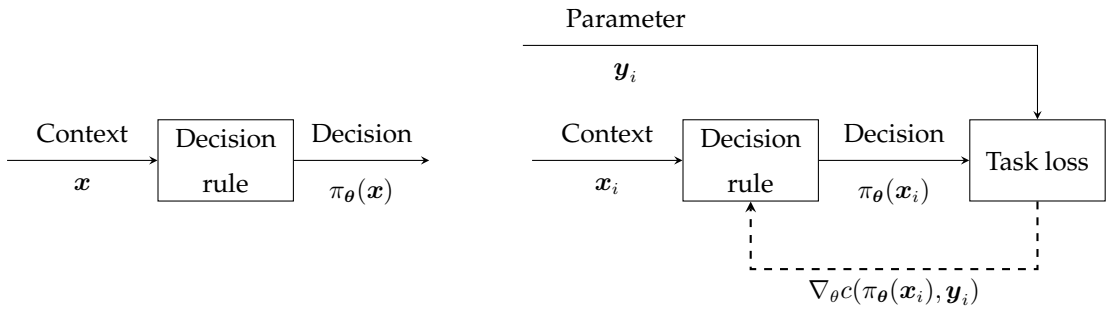

In simple terms, Problem~\eqref{eq:policy_opt} finds the policy $\pi_{\btheta^*} \in \Pi^{\btheta}$ that minimizes the expected costs over the training data. This decision pipeline is shown in Figure~\ref{fig:decision_rule}. Notice that there are two approximations of Problem~\eqref{eq:CSO:policy} made by Problem~\eqref{eq:policy_opt}. First, the policy is restricted to a hypothesis class that may not contain the true optimal policy. Second, the long-term expected costs are calculated over the empirical distribution $\Pemp$ rather than the true distribution~$\Prob$. Furthermore, Problem~\eqref{eq:policy_opt} focuses its policy optimization efforts on the overall performance (averaged over different \revised{covariates}) and disregards the question of making the policy achieve a good performance uniformly from one \revised{covariate} to another.

\subsubsection{Learning and optimization.} \label{sec:learnOptim}
The second and third frameworks combine a predictive component and an optimization component. The predictive component is a general model~$f_{\btheta}$, parameterized by~$\btheta$, whose role is to provide the input of the optimization component. \revised{For} any covariate $\x$, the intermediate input $f_{\btheta}(\x)$ can be interpreted as a predicted distribution that approximates the true conditional distribution $\Prob(\y \rvert \x)$ \revised{(or a sufficient statistic in the case of expected value-based models)}. The predictive component is typically learned from historical data.

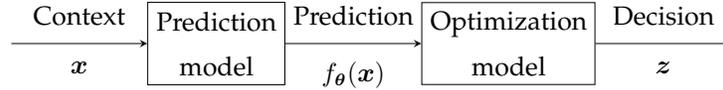
\begin{figure}[htb]
    \centering
    \resizebox{0.6\linewidth}{!}{\begin{tikzpicture}
 
 \node (A) at (0, 0) {};
 
\node [ draw,
        minimum width=2cm,
        minimum height=1.2cm,
        right=2cm of A,
        align=center
       ]  (ml) {Prediction\\model};
 
\node [draw,
        minimum width=2cm, 
        minimum height=1.2cm,
        right=2cm of ml,
        align=center
      ] (opt) {Optimization \\ model};

\draw[-stealth] (A) -- node[label= above:Context] {} node[label= below:$\x$] {} (ml.west);

\draw[-stealth] (ml.east) -- node[label= above:Prediction] {} node[label= below:$f_\btheta(\x)$] {} (opt.west);

\draw[-stealth] (opt.east) -- node[label= above:Decision] {} node[label= below:$\z$] {} ++ (2,0);

\end{tikzpicture}}
    \caption{Decision pipeline for learning and optimization.}
    \label{fig:learn_opt_pipeline}
\end{figure}

At decision time, a learning and optimization decision pipeline (see Figure~\ref{fig:learn_opt_pipeline}) solves the CSO problem under $f_{\btheta}(\x)$, namely:
\begin{equation}
    \label{eq:optim}
   z^*(\x, f_{\btheta}) \in \argmin_{\z \in \Z}h (\z, f_{\btheta}(\x)).
\end{equation}
The solution of Problem~\eqref{eq:optim} minimizes the expected cost with respect to the predicted distribution $f_{\btheta}(\x)$. Notice that the only approximation between Problem~\eqref{eq:optim} and the true CSO problem in \eqref{eq:CSO_h} lies in $f_{\btheta}$ being an approximation of $\Prob(\y \rvert \x)$. Since the predicted distribution changes with the \revised{covariate} $\x$, this pipeline also provides a policy. In fact, if the predictive component were able to perfectly predict the true conditional distribution $\Prob(\y\rvert \x)$ for any $\x$, the pipeline would recover the optimal policy $\pi^*$ given in \eqref{eq:CSO:policy}.

We now detail the second and third frameworks to address contextual optimization: SLO and ILO. They differ in the way the predictor $f_{\btheta}(\x)$ is trained using the historical data.

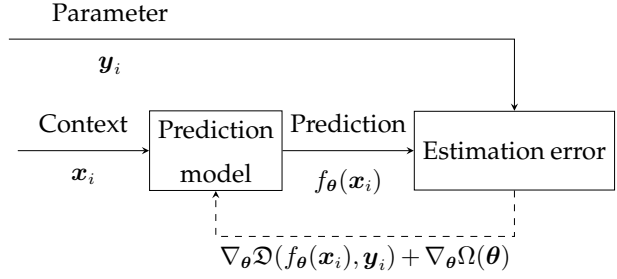
\begin{figure}[htb]
    \centering
\resizebox{0.5\linewidth}{!}{\begin{tikzpicture}
 
 \node (A) at (0, 0) {};
 
\node [ draw,
        minimum width=2cm,
        minimum height=1.2cm,
        right=2cm of A,
        align=center
       ]  (ml) {Prediction\\model};
 
\node [draw, %
        minimum width=2cm, 
        minimum height=1.2cm,
        right=2cm of ml,
        align=center
      ] (loss) {Estimation error};

\draw[-stealth] (A) -- node[label= above:Context] {} node[label= below:$\x_i$] {} (ml.west);
\draw[-stealth] (ml.east) -- node[label= above:Prediction] {} node[label= below:$f_\btheta(\x_i)$] {} (loss.west);
\draw[-stealth] (A) ++(0,1.7cm) -| (loss.north) node[label= above:Parameter, pos=0.1] {} node[label= below:$\V{y}_i$, pos=0.1] {};
\draw[-stealth, dashed] (loss.south) -- ++(0,-0.7cm) -| (ml.south) node [below,pos=0.25] {$\nabla_\btheta \mathfrak{D}(f_{\btheta}(\x_i), \V{y}_i)+ \nabla_\btheta\Omega(\btheta)$};

\end{tikzpicture}}
    \caption{SLO training pipeline for a given training example.}
    \label{fig:sequential_pipeline}
\end{figure}

\paragraph{Sequential learning and optimization.}
In this framework, the contextual predictor is obtained by minimizing an estimation error, $\rho$, between the conditional distribution given by $f_{\btheta}(\x)$ and the true conditional distribution of $\y$ given $\x$. Training the contextual \revised{parametric} predictor \revised{usually} implies solving the following estimation problem:
\begin{equation}
    \label{eq:seq_opt}
    \min_{\btheta} \rho(f_{\btheta}, \Pemp)+\Omega(\btheta)\mbox{ with }\rho(f_{\btheta}, \Pemp):= \Expect_{\Pemp}[\mathfrak{D}(f_{\btheta}(\x),\y)],
\end{equation}
where $\mathfrak{D}$ is a divergence function, e.g., negative log-likelihood and the regularization term $\Omega(\btheta)$ controls the complexity of $f_\btheta$. \revised{The conditional distribution can also take a non-parametric form, e.g. $k$-nearest neighbor, where $\btheta$ then captures hyper-parameters of the non-parametric model (such as the number of neighbors) selected through a form of cross-validation scheme.} The SLO training pipeline is shown in Figure~\ref{fig:sequential_pipeline}.

\begin{definition}[Expected value-based models]\label{def:exp_model}
   When the cost function $c(\x,\y)$ of the decision model is linear in $\y$, the problem of estimating a conditional distribution reduces to finding the expected value of the uncertain parameter given the covariates since $h(\z,\Prob(\y|\x))=\Expect_{\Prob(\y|\x)}[c(\z,\y)]=c(\z,\Expect_{\Prob(\y|\x)}[\y])$. Training the contextual predictor, therefore, reduces to a mean regression problem over a parameterized function $g_{\btheta}(\x)$. Specifically,
    \begin{align}
        \label{eq:regression}
        \min_{\btheta} \rho(g_{\btheta}, \Pemp)+ \Omega(\btheta)\mbox{ with }\rho(g_{\btheta}, \Pemp):= \Expect_{\Pemp} \left[d(g_{\btheta}(\x), \y)  \right],
    \end{align}
    for some distance metric $d$, usually the \revised{mean} squared errors. While the \revised{mean} squared error is known to asymptotically retrieve $g_{\hat{\btheta}}(\x)=\Expect_{\Prob(\y|\x)}[\y]$ \revised{as $N\rightarrow \infty$} under standard conditions, other distance metric or more general loss functions can also be used \citep{Hastie2009}. For any new \revised{covariate} $\x$, \revised{an action} is obtained using:
    \begin{equation}
        \label{eq:dp}
        z^*(\x,g_{\btheta})\in\argmin_{\z\in\Z}h(\z,\Dirac_{g_{\btheta}(\x)})=\argmin_{\z\in\Z}c(\z,g_{\btheta}(\x)),
    \end{equation}
    \revised{where, with a slight abuse of notation, $z^*$ now takes an estimator of the mean of the conditional distribution as the second argument, and} with $\Dirac_{\y}$ being the Dirac distribution putting all its mass at $\y$. In the remainder of this survey, we refer to these approaches as \textbf{expected value-based models}, while the more general models that prescribe using a conditional distribution estimator (i.e. $z^*(\x,f_\btheta)$) will be referred as a \textbf{conditional distribution-based models} when it is not clear from the context.
\end{definition}

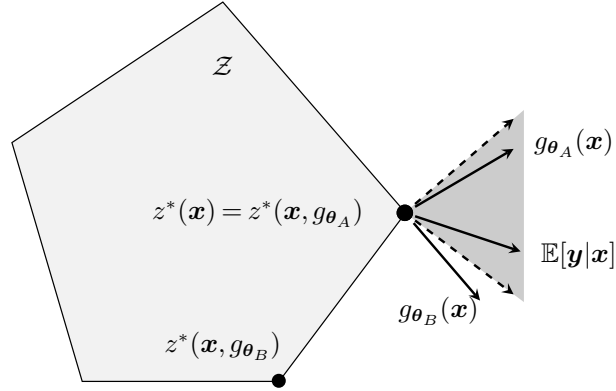
\begin{figure}[ht]
    \centering
   \resizebox{0.5\linewidth}{!}{\begin{tikzpicture}[>=stealth, line width=1pt]
    \node (zStar) at (4.6,2.4) {};
    \node (zB) at (2.8,0) {};
    \node (gA) at (6.3,3.4) {};
    \node (expY) at (6.4,1.8) {};
    \node (gB) at (5.8,1.0) {};

    \draw (0,0) -- (2.8,0) -- (4.6,2.4) -- (2,5.4) -- (-1,3.4) -- cycle;
    \fill[gray!10] (0,0) -- (2.8,0) -- (4.6,2.4) -- (2,5.4) -- (-1,3.4) -- cycle;

    \node (fanSouth) at (6.3, 1.125) {};
    \node (fanNorth) at (6.3, 3.87) {};
    \fill[black!20] (4.6,2.4) -- (6.3, 3.87) -- (6.3, 1.125) -- cycle;
    \draw[dashed, ->] (zStar) -- (fanNorth);
    \draw[dashed, ->] (zStar) -- (fanSouth);

    \node at (2,4.5) {$\mathcal{Z}$};
    \node at (2.5, 2.4) {$z^*(\x)=z^*(\x,g_{\btheta_A})$};
    \node at (2, 0.5) {$z^*(\x,g_{\btheta_B})$};
    
    \draw[->] (zStar) -- (gA) node[right] {$g_{\btheta_A}(\x)$};
    \draw[->] (zStar) -- (expY) node[right] {$\Expect[\y \rvert \x]$};
    \draw[->] (zStar) -- (gB) node[left] {$g_{\btheta_B}(\x)$};
    
    \fill[black] (zStar) circle (0.1);
    \draw (zStar) circle (0.1);
    \fill[black] (zB) circle (0.1);
    \draw (zStar) circle (0.1);
\end{tikzpicture}}
    \caption{Predicting $g_{\btheta_A}(\x)$ results in the optimal \revised{action} $z^*(\x, g_{\btheta_A}) = z^*(\x)$ whereas a small error resulting from predicting $g_{\btheta_B}(\x)$ leads to a suboptimal \revised{action} $z^*(\x, g_{\btheta_B})$ under $c(\x,\y):= -\y^\top \x$, i.e., $h(\z, \Prob(\y \rvert \x)) = -\Expect[\y\rvert \x]^\top \z$ (adapted from \citealt{elmachtoub_smart_2022}).}
    \label{fig:spo_example}
\end{figure}

\paragraph{Integrated learning and optimization.}
Sequential approaches ignore the mismatch between the prediction divergence $\mathfrak{D}$ and the cost function $c(\x,\y)$. Depending on the \revised{covariate}, a small prediction error about $\Prob(\y\rvert \x)$ 
may have a large impact on the prescription. In integrated learning, the goal is to maximize the prescriptive performance of the induced policy. That is, we want to train the predictive component to minimize the task loss (i.e., the downstream costs incurred by the decision) as stated in \eqref{eq:CSO:policy}. The prescriptive performance may guide the estimation procedure toward a \revised{solution with higher MSE (or any distance metric) that nevertheless} produces a nearly-optimal decision. This is illustrated in Figure~\ref{fig:spo_example}.

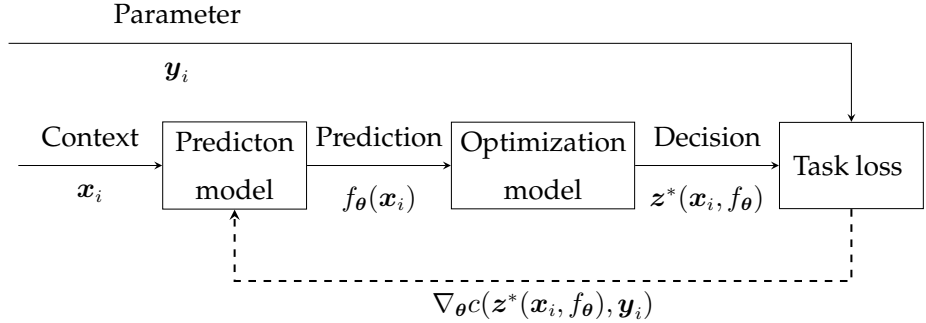
\begin{figure}[htb]
    \centering
    \resizebox{0.75\linewidth}{!}{\begin{tikzpicture}[->, >=stealth]
 \node (A) at (0, 0) {};
 
\node [ draw,
        minimum width=2cm,
        minimum height=1.2cm,
        right=2cm of A,
        align=center
       ]  (ml) {Predicton\\model};
 
\node [draw,
        minimum width=2cm, 
        minimum height=1.2cm,
        right=2cm of ml,
        align=center
      ] (opt) {Optimization\\ model};

\node [draw,
        minimum width=2cm, 
        minimum height=1.2cm,
        right=2cm of opt,
        align=center
      ] (loss) {Task loss };
\draw (A) -- node[label= above:Context] {} node[label= below:$\x_i$] {} (ml.west);
\draw (ml.east) -- node[label= above:Prediction] {} node[label= below:$f_\btheta(\x_i)$] {} (opt.west);
\draw (opt.east) -- node[label= above:Decision] {} node[label= below:{$\z^*(\x_i, f_{\btheta})$}] {} ++ (2,0);

\draw[thick, dashed] (loss.south) -- ++(0,-1cm)  -| node [below,pos=0.25] {$\nabla_\btheta c(\z^*(\x_i, f_{\btheta}), \V{y}_i)$}(ml.south);

\draw (A) ++(0,1.7cm) -| (loss.north) node[label= above:Parameter, pos=0.1] {} node[label= below:$\V{y}_i$, pos=0.1] {};

\end{tikzpicture}}
    \caption{ILO training pipeline for a given training example.}
    \label{fig:integrated_pipeline}
\end{figure}

Finding the best parameterization of a contextual predictor that minimizes the downstream expected costs of the CSO solution can be formulated as the following problem:
\begin{equation}
      \min_{\btheta} H(z^*(\cdot, f_{\btheta}),\Pemp) = \min_{\btheta} \Expect_{\Pemp} \big[c(z^*(\x, f_{\btheta}), \y )\big]. \label{eq:integ_opt:objt}
\end{equation}
The objective function in \eqref{eq:integ_opt:objt} minimizes the \revised{expected} cost of the policy over the empirical distribution. The policy induced by this training problem is thus optimal with respect to the predicted distribution and minimizes the average historical costs over the whole training data. Figure~\ref{fig:integrated_pipeline} describes how the downstream cost is propagated by the predictive model during training. This training procedure necessarily comes at the price of heavier computations because an optimization model needs to be solved for each data point, and differentiation needs to be applied through an $\argmin$ operation.

\subsection{Summary}
This section \revised{presented} the main pipelines proposed in recent years to address contextual optimization. Although these pipelines all include a learning component, they differ significantly in their specific structures and training procedures. Overall, there are several design choices that the decision-maker should make when tackling contextual optimization:
\begin{inlinelist}
    \item \revised{the type of loss function used during training, which defines whether an approach belongs to the decision rule (using ERM), sequential (minimizing the estimation error), or integrated paradigm (minimizing the downstream cost of the policy),}
    \item the class of the predictive model (e.g., linear, neural network, or random forest) and its hyperparameters.
\end{inlinelist}
Each design choice has its own inductive bias and may imply specific methodological challenges, especially for ILO. In general, it is a priori unclear what combination of choices will lead to the best performance with limited data; therefore, pipelines have to be evaluated experimentally.

In the following sections, we survey the recent literature corresponding to the three main frameworks for contextual optimization using the notation introduced so far, which is summarized in Table~\ref{tab:notation_table}. A list of abbreviations used in this survey is given in Appendix \ref{appen:abbrev}.
\begin{table}[htb]
\centering
\caption{Notation: distributions, variables, and operators.}
\label{tab:notation_table}
\begin{adjustbox}{max width=\textwidth}
    \begin{tabular}{*{2}{c}l} 
    \toprule
     & Domain & Description \\
    \midrule
    $\Prob$            & $\mathcal{M}(\X \times \Y)$ & True (unknown) joint distribution of $(\x,\Vxi)$\\
    $\hat{\Prob}_N$    & $\mathcal{M}(\X \times \Y)$ & Joint empirical  distribution of $(\x,\Vxi)$\\
    $\Dirac_y$         & $\mathcal{M}(\Y)$           & Dirac distribution that puts all of its weight on $\Vxi$   \\
    $\x$               & $\X \subseteq \Re^{\revised{d_{\x}}}$    & Contextual information    \\
    $\V{y}$            & $Y \subseteq \Re^{\revised{d_{\y}}}$     & Uncertain parameters      \\
    $\z$               & $\Z \subseteq \Re^{\revised{d_{\z}}}$    & A feasible \revised{action}       \\
    $\btheta$          & $\Theta$                    & Parameters of a prediction model \\
    $\hat{\btheta}$    & $\Theta$                    & Optimal parameter value that minimizes the estimation error \\
    $c(\z, \Vxi )$     & $\Re$                       & Cost of an \revised{action} $\z$ under $\Vxi$ \\
    $h(\z, \mathbb{Q}_{\xi} )$  & $\Re$              & Expected cost of an \revised{action} $\z$ under $\mathbb{Q}_{\xi}$ (a distribution over $\Vxi$) \\
    $H(\pi, \mathbb{Q})$        & $\Re$              & Expected cost of a policy $\pi$ under $\mathbb{Q}$ (a distribution over $(\x,\Vxi)$) \\
    $f_\btheta(\x)$    & \revised{$\mathcal{M}(\Y)$} & Estimate of the conditional distribution of $\V{y}$ given $\x$\\
    $g_\btheta(\x)$    & \revised{$\Re^{d_{\y}}$}       & Estimate of the conditional expectation of $\V{y}$ given $\x$ \\
    $\pi^*(\x)$        & $\Z$                        & Optimal solution of CSO under true conditional distribution $\Prob(\Vxi|\x)$ \\
    $\pi_\btheta(\x)$  & $\Z$                        & \revised{Action} prescribed by a policy parameterized by $\btheta$ for context $\x$\\
    $z^*(\x)$          & $\Z$                        & Optimal solution to the CSO problem under the true conditional distribution $\Prob(\Vxi|\x)$ \\   
    $z^*(\x, f_{\btheta})$ & $\Z$                    & Optimal solution to the CSO problem under the conditional distribution $f_{\btheta}(\x)$\\
    $z^*(\x, g_{\btheta})$ & $\Z$                    & Optimal solution to the CSO problem under the Dirac distribution $\delta_{g_{\btheta}(\x)}$\\
    $\rho(f_{\btheta}, \Pemp)$ & $\Re$               & Expected prediction error for distribution model $f_{\btheta}$ based on empirical distribution $\Pemp$\\
    $\rho(g_{\btheta}, \Pemp)$ & $\Re$               & Expected prediction error for point prediction model $g_{\btheta}$ based on empirical distribution $\Pemp$\\
    \bottomrule
    \end{tabular}
\end{adjustbox}
\end{table}

\section{Decision rule optimization}
\label{sec:decision_rules}
Decision rules obtained by solving the ERM in Problem~\eqref{eq:policy_opt} minimize the cost of a policy on the task, that is, the downstream optimization problem. Policy-based approaches are especially efficient computationally at decision time since it suffices to evaluate the estimated policy. No optimization problem needs to be solved once the policy is trained. 
We defined the decision rule approach as employing a parameterized mapping  $\pi_{\btheta}(\x)$, e.g., linear policies \citep{ban_big_2019} or a neural network \citep{Oroojlooyjadid2020}. Since policies obtained using neural networks lack interpretability, linear decision rules are widely used. 

\subsection{Linear decision rules}
\cite{ban_big_2019} show that an approach based on the sample-average approximation (SAA) that disregards side information can lead to inconsistent decisions (i.e., asymptotically suboptimal) for a newsvendor problem. Using linear decision rules (LDRs), they study two variants of the newsvendor problem with and without regularization:
\begin{align*}
     \min_{\pi:\revised{\exists \btheta\in\Re^{\revised{d_{\x}}}, }\pi(\x) = \revised{\btheta}^\top \x,\revised{\forall \x}} H(\pi, \Pemp) + \Omega(\pi)=\min_{\revised{\btheta}}\frac{1}{N} \sum_{i=1}^N u(y_i -\revised{\btheta}^\top \x_i)^+ + o(\revised{\btheta}^\top \x_i - y_i)^+ + \lambda \lVert \revised{\btheta}\rVert_k^2,
\end{align*}
where $u$ and $o$ denote the per unit backordering (underage) and holding (overage) costs. \cite{ban_big_2019} show that for a linear demand model, the generalization error for the ERM model scales as $\text{O}(\revised{d_{\x}}/\sqrt{N})$ when there is no regularization and as $\text{O}(\revised{d_{\x}}/ (\sqrt{N} \lambda))$ with regularization. However, one needs to balance the trade-off between generalization error and bias due to regularization to get the optimal performance from using LDRs.
\cite{ban_big_2019} consider unconstrained problems
 because it is difficult to ensure the feasibility of policies and maintain computational tractability using the ERM approach.

\revised{\cite{Bertsimas_Kallus_2020} present a general theory for generalization bounds of decision rules based on Rademacher complexity that goes beyond LDR, although their main examples in this context pertain to LDR}. Unfortunately, LDRs may not be asymptotically optimal in general. To generalize LDRs, one can consider decision rules that are linear in the transformations of the covariate vector \citep{ban_big_2019}. It is also possible to lift the covariate vector to a reproducing kernel Hilbert space \citep[RKHS,][]{aronszajn1950theory}, as seen in the next section.

\subsection{RKHS-based decision rules}
To obtain decision rules that are more flexible than linear ones with respect to $\x$, it is possible to lift the covariate vector to an RKHS in which LDRs might achieve better performance. 
Let $K:\X \times \X \rightarrow \Re$ be the symmetric positive definite kernel associated with the chosen RKHS, e.g., the Gaussian kernel $K(\x_1,\x_2):=\exp(-\|\x_1-\x_2\|^2/(2\sigma^2))$. 
Given $K$, the RKHS $\mathcal{H}_K$ is defined as the closure of a set of functions given below:
\begin{equation*}
    \Big\{\varphi: \X \rightarrow \Re \, \rvert \, \exists \, L\in \mathbb{N}, \, \bv_1,\bv_2, \cdots, \bv_L \in \X, \, \varphi(\x) = \sum_{l=1}^L a_l K(\bv_l,\x), \forall \x\in\X \Big\},
\end{equation*}
with the inner product of $\varphi_1(\x)=\sum_{i=1}^{L_1} a_1^iK(\bv_1^i,\x)$ and $\varphi_2(\x)=\sum_{j=1}^{L_2} a_2^jK(\bv_2^j,\x)$ given by:
\[ \left\langle \varphi_1, \varphi_2 \right\rangle = \sum_{i=1}^{L_1}\sum_{j=1}^{L_2} a_1^ia_2^j K(\bv_1^i,\bv_2^j).\]
\cite{Bertsimas_Koduri_2022} approximate the optimal policy with a linear policy in the RKHS, i.e. $\pi_\varphi(\x):= \langle \varphi,K(\x,\cdot)\rangle$ when $\revised{d_{\z}}=1$, and show using the representer theorem \citep[see Theorem 9 in][]{hofmann2008kernel} that the solution of the following regularized problem:
\begin{align*}
    \min_{\varphi\in\mathcal{H}_K} H(\pi_\varphi,\Pemp)   +\lambda \lVert \varphi\rVert_2^2,
\end{align*}
takes the form $\pi^*(\x) = \sum_{i=1}^N K(\x_i, \x) a_i^*$. Hence, this reduces the decision rule problem to:
\begin{align*}
    \min_{\V{a}\in\Re^N} H\left(\sum_{i=1}^N K(\x_i, \cdot) a_i,\Pemp\right)   +\lambda \sum_{i=1}^N\sum_{j=1}^N K(\x_i,\x_j)a_i a_j.
\end{align*}
This can be extended to $\revised{d_{\z}}>1$ by employing one RKHS for each $z_i$.

This RKHS approach appeared earlier in \cite{ban_big_2019} and \cite{Bazier-Matte_Delage_2020} who respectively study the data-driven single item newsvendor and single risky asset portfolio problems and establish bounds on the out-of-sample performance. \cite{Bertsimas_Koduri_2022} show the asymptotic optimality of RKHS-based policies.
\cite{Notz2022} study a two-stage capacity planning problem with multivariate demand and vector-valued capacity decisions for which the underlying demand distribution is difficult to estimate in practice. Similar to \cite{Bazier-Matte_Delage_2020}, the authors 
optimize over policies that are linear in the RKHS associated with the Gaussian kernel and identify generalization error bounds. For large dimensional problems, this kernel is shown to have a slow convergence rate, and as a result, the authors propose instead using a data-dependent random forest kernel.

\subsection{Non-linear decision rules}

Many non-linear decision rule approaches have been experimented with. \cite{ZhangNNNewsvendor}, \citet{Huber2019}, and \citet{Oroojlooyjadid2020} study the value of training a DNN to learn the ordering policy of a newsvendor problem. It is well-known that neural networks enjoy the universal approximation property; that is, they can approximate any continuous function arbitrarily well \citep{cybenko1989approximation, lu2017expressive}. For constrained problems, one can use softmax as the final layer to ensure that decisions lie in a simplex, e.g., in a portfolio optimization problem \citep{zhang2021universal}. Yet, in general, the output of a neural network might not naturally land in the feasible space $\Z$. To circumvent this issue, \citet{Chen_Tanneau_VanHentenryck2023} introduce an application-specific differentiable repair layer that projects the solution back to feasibility. \cite{Rychener_Sutter_2023} show that the decision rule obtained by using the stochastic gradient descent~(SGD) method to train DNN-based policies approximately minimizes the Bayesian posterior loss. 

Exploiting the fact that the optimal solution of a newsvendor problem is a quantile of the demand distribution, \citet{Huber2019} further train an additive ensemble of decision trees using quantile regression to produce the ordering decision. They test these algorithms on a real-world dataset from a large German bakery chain. \cite{doi:10.1287/ijoo.2018.0005}, \cite{doi:10.1287/mnsc.2020.3592}, and \cite{parisa} optimize decision tree-based decision rules to address the multi-item newsvendor, treatment planning, and optimal stopping problems, respectively. A tutorial on DNN-based decision rule optimization is given in \cite{Shlezinger_Eldar_Boyd_2022}.

\citet{zhang2023_pwaffine} introduce piecewise-affine decision rules and provide non-asymptotic and asymptotic consistency results for unconstrained and constrained problems, respectively. The policy is learned through a stochastic majorization-minimization algorithm, and experiments on a constrained newsvendor problem show that piecewise-affine decision rules can outperform the RKHS-based policies.

\subsection{Distributionally robust decision rules}
Most of the literature on policy learning assumes a parametric form $\Pi^\btheta$ for the policy. A notable exception is \cite{Zhang_Yang_Gao}, which studies a distributionally robust contextual newsvendor problem under the type-1 Wasserstein ambiguity set without assuming an explicit structure on the policy class. The type-\revised{$p$} Wasserstein distance (earth mover's distance) between distributions $\Prob_1$ and $\Prob_2$ is given by:
\[W_p(\Prob_1,\Prob_2) = \inf_{\gamma \in \mathcal{M}(\Y^2)} \left( \int_{\Y \times \Y} \lVert y_1-y_2\rVert^p \gamma(dy_1,dy_2) \right)^{\frac{1}{p}},\]
where $\gamma$ is a joint distribution of $y_1$ and $y_2$ with marginals $\Prob_1$ and $\Prob_2$. The distributionally robust model in \cite{Zhang_Yang_Gao} avoids the degeneracies of ERM for generic $\Pi$ by defining an optimal \quoteIt{Shapley} extension to the scenario-based optimal policy. Mathematically,
\begin{align*}
    \min_{\pi\in \Pi} \sup_{\mathbb{Q} \in \mathcal{M}(\X \times \Y)} \{H(\pi, \mathbb{Q}): \mathcal{W}(\mathbb{Q}, \Pemp) \leq r\} \;\equiv\; \min_{\pi:\hat{\X}\rightarrow \Z} \sup_{\mathbb{Q} \in \mathcal{M}(\hat{\X} \times \Y)} \{H(\pi, \mathbb{Q}): \mathcal{W}(\mathbb{Q}, \Pemp) \leq r\},
\end{align*}
where $\hat{\X}:=\cup_{i=1}^N\{\x_i\}$ and $\mathcal{M}(\hat{\X} \times \Y)$ is the set of all distribution supported on $\hat{\X} \times \Y$.

Prior to the work of \cite{Zhang_Yang_Gao}, many have considered distributionally robust versions of the decision rule optimization problem in the non-contextual setting \citep{yanikoglu2019survey} while \cite{Bertsimas2023dynamic} use LDRs to solve dynamic optimization problems with side information.

\cite{yang_shang2023} point out that the perturbed distributions in the Wasserstein ambiguity set might have a different conditional information structure than the estimated conditional distribution. They introduce a distributionally robust optimization~(DRO) problem with causal transport metric \citep{backhoff2017causal,lassalle2018causal} that places an additional causality constraint on the transport plan compared to the Wasserstein metric. Tractable reformulations of the DRO problem are given under LDRs as well as for one-dimensional convex cost functions. \cite{Rychener_Sutter_2023} present a Bayesian interpretation of decision rule optimization using SGD and show that their algorithm provides an unbiased estimate of the worst-case objective function of a DRO problem as long as a uniqueness condition is satisfied. The authors note that the Wasserstein ambiguity set violates this condition and thus use the Kullback-Leibler~(KL) divergence \citep{kullback1951information} to train the models. 

\section{Sequential learning and optimization}
\label{sec:sequential}
In reviewing contextual optimization approaches that are based on SLO, we distinguish two settings:
\begin{inlinelist}
    \item a more traditional setting where the conditional distribution is learned from data and used directly in the optimization model, and
    \item a setting that attempts to produce decisions that are robust to model misspecification.
\end{inlinelist}
An overview of the methods presented in this section is given in Table~\ref{tab:slo_methods}.
\newcolumntype{E}{>{\centering\arraybackslash}m{1.4cm}}
\begin{table}
\centering
\caption{Overview of contextual optimization papers in the SLO framework.}
\label{tab:slo_methods}
\begin{adjustbox}{max width=\textwidth}
\begin{threeparttable}[ht]
    \begin{tabular}{l*{11}E} 
        \toprule
        & \multicolumn{3}{c}{Method} & \multicolumn{2}{c}{Regularization} & \multicolumn{6}{c}{Learning model}  \\
        \cmidrule(lr){2-4} \cmidrule(lr){5-6} \cmidrule(lr){7-12} & \model{rCSO} & \model{wSAA} & \model{EVB} & {Reg.~CSO} & DRO & General & Linear & Kernel & kNN & DT & RF \\
        \midrule
        \cite{Hannah2010}             & \cross & \tick & \cross & \cross & \cross & \cross & \cross & \tick & \cross & \cross & \cross\\
                \cite{ferreira2016analytics} & \cross & \cross & \tick & \cross & \cross & \cross & \cross & \cross & \cross & \tick & \cross\\
        \cite{Ban2019residual_tree} & \tick & \cross & \cross & \cross & \cross & \cross & \tick & \cross & \cross & \cross & \cross\\
        \cite{Chen_Paschalidis_2019}  & \cross & \tick & \cross & \cross & \tick & \cross & \cross & \cross & \tick & \cross & \cross\\
         \cite{Bertsimas_Kallus_2020}  & \cross & \tick & \cross & \cross & \cross & \cross & \tick & \cross & \tick & \tick & \tick\\
        \cite{Kannan2020residuals}  & \tick & \cross & \cross & \cross & \cross & \tick & \tick & \tick & \tick & \tick & \tick\\
                \cite{Kannan_Bayraksan_Luedtke_2021} & \tick & \cross & \cross & \cross & \tick& \tick & \tick & \tick & \tick & \tick & \tick\\
        \cite{Liu_He_Max_Shen_2021} & \cross & \cross & \tick & \cross & \cross & \cross & \tick & \cross & \cross & \tick & \cross\\
         \cite{Srivastava_Wang_Hanasusanto_Ho_2021} & \cross & \tick & \cross & \tick & \cross & \cross & \cross & \tick & \cross & \cross & \cross\\
                \cite{wang2021distributionally}  & \cross & \tick & \cross & \cross & \tick & \cross & \cross & \tick & \cross & \cross & \cross\\
         \cite{Bertsimas2022bootstrap}  & \cross & \tick & \cross & \cross & \tick & \cross & \cross & \tick & \tick & \cross & \cross\\
        \cite{Deng_Sen_2022}        & \tick & \cross & \cross & \cross & \cross & \tick & \tick & \tick & \tick & \tick & \tick\\
         \cite{Esteban_Perez_Morales_2022}  & \cross & \tick & \cross & \cross & \tick & \cross & \cross & \tick & \tick & \cross & \cross\\
        \cite{Kannan2022} & \tick & \cross & \cross & \cross & \tick& \tick & \tick & \tick & \tick & \tick & \tick\\
        \cite{Lin2022} & \cross & \tick & \cross & \tick & \cross & \cross & \cross & \cross & \tick & \tick & \tick\\
        \cite{Nguyen_Zhang_Blanchet_Delage_Ye_2022}  & \cross & \tick & \cross & \cross & \tick & \cross & \cross & \cross & \tick & \cross & \cross\\
         \cite{Notz2022}               & \cross & \tick & \cross & \cross & \cross & \cross & \cross & \tick & \cross & \tick & \cross\\
         \cite{Zhu_Xie_Sim_2022} & \cross & \cross & \tick & \cross & \tick & \tick & \tick & \tick & \tick & \tick & \tick\\
        \cite{Perakis_Sim_Tang_Xiong_2023} & \tick & \cross & \cross & \cross & \tick& \cross & \tick & \cross & \cross & \cross & \cross\\
        \bottomrule
    \end{tabular}
    \begin{tablenotes}
        \small
        \item Note: we distinguish between regularized CSO models (Reg.~CSO) and DRO-based regularization; an approach is classified as \quoteIt{General} if its learning model is not restricted to specific classes.
    \end{tablenotes}
    \end{threeparttable}
\end{adjustbox}
\end{table}

\subsection{Learning conditional distributions}
\label{sec:non_lin}

Most of the recent literature has employed discrete models for $f_\btheta(\x)$. This is first motivated from a computational perspective by the fact that the CSO Problem \eqref{eq:optim} is easier to solve in this setting. In fact, more often than not, the CSO under a continuous distribution needs to be first replaced by its SAA to be solved \citep{Shapiro2014_lectures}. From a statistical viewpoint, it can also be difficult to assess the probability of outcomes that are not present in the dataset, thus justifying fixing the support of $\y$ to its observed values.

\subsubsection{Residual-based distribution.} A first approach \citep[found in][]{sen2017learning, Kannan2020residuals, Deng_Sen_2022} is to use the errors of a trained regression model (i.e., its residuals) to construct conditional distributions. Let $g_{\hat{\btheta}}$ be a regression model trained to predict the response~$\y$ from the covariate $\x$, thus minimizing an estimation error $\rho$ as in \eqref{eq:regression}. The residual error of sample $i$ is given by $\V{\epsilon}_i = \y_i- g_{\hat{\btheta}}(\x_i) $. The set of residuals measured on the historical data, $\{ \V{\epsilon}_i \}_{i=1}^N$, is then used to form the conditional distribution $f_{\btheta}(\x) = \Prob^{\mymbox{ER}}(\x) := \frac{1}{N}\sum_{i=1}^{N} \Dirac_{\rm{proj}_{\Y}(g_{\hat{\btheta}}(\x)+\V{\epsilon}_i)}$, with $\rm{proj}_{\Y}$ denoting the projection on the support~$\Y$. 
The residual-based CSO~(\model{rCSO}) problem is now given by:
\begin{align}
    \text{(\model{rCSO})} \quad &\min_{\z\in \Z} h(\z, \Prob^{\mymbox{ER}}(\x))
    \label{eq:ER_SAA}
\end{align}

The advantage of residual-based methods is that they can be applied in conjunction with any trained regression model. \cite{Ban2019residual_tree} and \cite{Deng_Sen_2022} build conditional distributions for two-stage and multi-stage CSO problems using the residuals obtained from parametric regression on the historical data.

Notice that, in this approach, the historical data is used twice: to train the regression model $g_\btheta$, and to measure the residuals $\V{\epsilon}_i$. This can lead to an underestimation of the distribution of the residual error. To remove this bias, \cite{Kannan2020residuals} propose a leave-one-out model (also known as \textit{jackknife}). They measure the residuals as $\tilde{\V{\epsilon}}_i = \y_i-g_{\hat{\btheta}_{-i}}(\x_i) $, where $\hat{\btheta}_{-i}$ is trained using all the historical data except the $i$-th sample $(\x_i, \y_i)$. This idea can also be applied to the heteroskedastic case studied in \citet{Kannan_Bayraksan_Luedtke_2021}, where the following conditional distribution is obtained by first estimating the conditional covariance matrix $\hat{Q}(\x)$ (a positive definite matrix for almost every $\x$) and then forming the residuals $\hat{\V{\epsilon}}_i = [\hat{Q}(\x_i)]^{-1}(\y_i-g_{\hat{\btheta}}(\x_i))$: 
\[f_{\btheta}(\x):= \frac{1}{N}\sum_{i=1}^{N} \Dirac_{\rm{proj}_{\Y} ( g_{\hat{\btheta}}(\x) + \hat{Q}(\x) \hat{\V{\epsilon}}_i)}.\]

\subsubsection{Weight-based distribution.} \label{sec:wSAA} A typical approach for formulating the CSO problem is to assign weights to the observations of the uncertain parameters in the historical data and solving the weighted SAA problem (\model{wSAA}) given by \citep{Bertsimas_Kallus_2020}:
\begin{align}
    \label{eq:wsaa}
    \text{(\model{wSAA})} \quad& \min_{z\in \Z} h\left(\z, \sum_{i=1}^N w_i(\x)\Dirac_{\y_i}\right).
\end{align}
In this case, the conditional distribution $f_\btheta(\x) = \sum_{i=1}^N w_i(\x)\Dirac_{\y_i}$ is fully determined by the function used to assign a weight to the historical samples. Different approaches have been proposed to determine the sample weights with ML methods.

\paragraph{Weights based on proximity.} Sample weights can be assigned based on the distance between a \revised{covariate}~$\x$ and each historical sample~${\x_i}$. 
For instance, a $k$-nearest neighbor~(kNN) estimation gives equal weight to the $k$ closest samples in the dataset and zero weight to all the other samples. That is, $w_i^{\mymbox{kNN}}(\x):= (1/k)\Indc[\x_i \in \mathcal{N}_k(\x)]$, where $\mathcal{N}_k(\x)$ denotes the set of $k$ nearest neighbors of $\x$ and $\Indc[\cdot]$ is the indicator function. Even though it may appear simple, this non-parametric approach benefits from asymptotic consistency guarantees on its prescriptive performance.
Another method to determine sample weights is to use kernel density estimators \citep{Hannah2010, Srivastava_Wang_Hanasusanto_Ho_2021, ban_big_2019}. The Nadaraya-Watson (NW) kernel estimator \citep{watson1964smooth, nadaraya1964estimating} employs a weight function:
\begin{align*}
    w_i^{\mymbox{KDE}}(\x) := \frac{K\left((\x-\x_i)/\btheta\right)}{\sum_{j=1}^N K\left((\x-\x_j)/\btheta\right)} 
\end{align*}
where $K$ is a kernel function and $\btheta$ denotes its bandwidth parameter. 
Different kernel functions can be used, e.g., the Gaussian kernel defined as $K(\V{\Delta}) \propto \exp(-{\lVert \V{\Delta}\rVert} ^2)$. \citet{Hannah2010} also use a Bayesian approach that exploits the Dirichlet process mixture to assign sample weights.

\paragraph{Weights based on random forest.}
Weights can also be designed based on random forest \revised{regressors} \citep{Bertsimas_Kallus_2020}. In its simplest setting, the weight function of a decision tree \revised{regressor} is given by:
\begin{equation*}
    w_i^{t}(\x) := \frac{\Indc\left[ \mathcal{R}_t(\x) = \mathcal{R}_t(\x_i)\right] }{\sum_{j=1}^N \Indc\left[{\mathcal{R}_t(\x) = \mathcal{R}_t(\x_j)}\right]}
\end{equation*}
where $\mathcal{R}_t(\x)$ denotes the terminal node of tree $t$ that contains covariate $\x$. Thus, a decision tree assigns equal weights to all the historical samples that end in the same leaf node as $\x$. The random forest weight function generalizes this idea over many random decision trees. Its weight function is defined as:
\begin{equation*}
    w_i^{\mymbox{RF}}(\x) := \frac{1}{T} \sum_{t=1}^T w_i^{t}(\x),
\end{equation*}
where $w_i^{t}$ is the weight function of tree $t$. Random forests are typically trained in order to perform an inference task, e.g. regression, or classification, but can also be used and interpreted as non-parametric conditional density estimators.

\cite{Bertsimas_Kallus_2020} provide conditions for the asymptotic optimality (see Definition~\ref{def:optimal} in Appendix~\ref{sec:appen:definitions}) and consistency (see Definition~\ref{def:const} in Appendix~\ref{sec:appen:definitions}) of prescriptions obtained by solving Problem~\eqref{eq:wsaa} with the weights functions given by kNN, NW kernel estimator, and local linear regression. 

\subsubsection{Expected value-based models.}
As described in Definition~\ref{def:exp_model}, when the cost function is linear, the training pipeline of SLO reduces to conditional mean estimation. For instance, \cite{ferreira2016analytics} train regression trees to forecast daily expected sales for different product categories in an inventory and pricing problem for an online retailer. Alternatively, one may attempt to approximate the conditional density $f_\btheta(\x)$ using a point prediction $g_\btheta(\x)$. For example, \cite{Liu_He_Max_Shen_2021} study a last-mile delivery problem, where customer orders are assigned to drivers, and replace the conditional distribution of the stochastic travel time with a point predictor (e.g. a linear regression or decision tree) that accounts for the number of stops, total distance of the trip, etc. 

\subsection{Regularization and distributionally robust optimization} 

While non-parametric conditional density estimation methods benefit from asymptotic consistency \citep{Bertsimas_Kallus_2020,Notz2022}, they are known to produce overly optimistic policies when the size of the covariate vector is large \citep[see discussions in][]{Bertsimas2022bootstrap}. To circumvent this issue, authors have proposed to either regularize the CSO problem \citep{Srivastava_Wang_Hanasusanto_Ho_2021, Lin2022} or to cast it as a DRO problem. In the latter case, one attempts to minimize the worst-case expected cost over the set of distributions $\mathcal{B}_r(f_{\btheta}(\x))$ that lie at a distance $r$ from the estimated distribution $f_{\btheta}(\x)$:
\begin{align}
    \min_{z\in \Z} \sup_{\mathbb{Q}_y\in \mathcal{B}_r(f_{\btheta}(\x))} h(\z, \mathbb{Q}_y).\label{eq:DRO:gen}
\end{align}

\citet{Bertsimas2022bootstrap} generate bootstrap data from the training set and use it as a proxy for the \quoteIt{out-of-sample disappointment} of a\revised{n action} $\z$ resulting from the out-of-sample cost exceeding the budget given by $\sup_{\mathbb{Q}_y\in \mathcal{B}_r(f_{\btheta}(\x))} h(\z, \mathbb{Q}_y)$. They show that for the NW kernel estimator and KNN estimator, the DRO, under a range of ambiguity sets, can be reformulated as a convex optimization problem. Using KL divergence to measure the distance between the probability distributions, they obtain guarantees (\quoteIt{bootstrap robustness}) with respect to the estimate-then-optimize model taking bootstrap data as a proxy for out-of-sample data. Taking the center of Wasserstein ambiguity set \citep[see][]{kantorovich1958space} to be NW kernel estimator, \cite{wang2021distributionally} show that the distributionally robust newsvendor and conditional value at risk (CVaR) portfolio optimization problems can be reformulated as convex programs. They provide conditions to obtain asymptotic convergence and out-of-sample guarantees on the solutions of the DRO model. 

\cite{Chen_Paschalidis_2019} study a distributionally robust kNN regression problem by combining point estimation of the outcome with a DRO model over a Wasserstein ambiguity set \citep{chen2018robust} and then using kNN to predict the outcome based on the weighted distance metric constructed from the estimates. Extending the methods developed in \cite{NEURIPS2020_adf854f4}, \cite{Nguyen_Zhang_Blanchet_Delage_Ye_2022} study a distributionally robust contextual portfolio allocation problem where worst-case conditional return-risk tradeoff is computed over an optimal transport ambiguity set consisting of perturbations of the joint distribution of covariates and returns. Their approach generalizes the mean-variance and mean-CVaR model, for which the distributionally robust models are shown to be equivalent to semi-definite or second-order cone representable programs. 
\citet{Esteban_Perez_Morales_2022} solve a DRO problem with a novel ambiguity set that is based on trimming the empirical conditional distribution, i.e., reducing the weights over the support points. The authors show the link between trimming a distribution and partial mass transportation problem, and an interested reader can refer to \cite{ESTEBANPEREZ20231047} for an application in the optimal power flow problem. 

A distributionally robust extension of the \model{rCSO} model is presented in \cite{Kannan_Bayraksan_Luedtke_2021} and \cite{Kannan2022}. It hedges against all distributions that lie in the $r$ radius of the (Wasserstein) ambiguity ball centered at the estimated distribution $\Prob^{\mymbox{ER}}(\x)$. \cite{Perakis_Sim_Tang_Xiong_2023} propose a DRO model to solve a two-stage multi-item joint production and pricing problem with a partitioned-moment-based ambiguity set constructed by clustering the residuals estimated from an additive demand model.

\cite{Zhu_Xie_Sim_2022} considers an expected value-based model and suggests an ambiguity set that is informed by the estimation metric used to train $g_{\hat{\btheta}}$. Namely, they consider:
\[\min_{\z\in\Z}\sup_{\btheta\in\mathcal{U}(\hat{\btheta},r)} c(\z,g_\btheta(\x)),\]
with 
\[\mathcal{U}(\hat{\btheta},r):=\{\btheta\in\btheta|\rho(g_\btheta,\Pemp)\leq \rho(g_{\hat{\btheta}},\Pemp) + r\}.\]
They show how finite-dimensional convex reformulations can be obtained when $g_\btheta(\x):=\btheta^T\x$, and promote the use of a \quoteIt{robustness optimization} form.

\section{Integrated learning and optimization}
\label{sec:integrated_learning}
As discussed previously, ILO is an end-to-end framework that includes three components in the training pipeline:
\begin{inlinelist}
    \item a prediction model that maps the covariate to a predicted distribution (or possibly a point prediction),
    \item an optimization model that takes as input a prediction and returns a decision, and
    \item a task-based loss function that captures the downstream optimization problem.
\end{inlinelist}
The parameters of the prediction model are trained to maximize the prescriptive performance of the policy, i.e., it is trained on the task loss incurred by this induced policy rather than the estimation loss.

Next, we discuss several methods for implementing the ILO approach. We start by describing the different models that are used in ILO (Section~\ref{sec:ILOmodels}), and then we present the algorithms used to perform the training. We divide the algorithms into four categories. Namely, training using unrolling (Section~\ref{sec:unrolling}), implicit differentiation (Section~\ref{sec:implicitdiff}), a surrogate differentiable loss function (Section~\ref{sec:surrogateDiffLoss}), and a differentiable optimizer (Section~\ref{sec:surrogateDiffOpt}). An overview of the methods presented in this section is given in Table~\ref{tab:ilo_methods}.

\begin{table}
\centering
\caption{Overview of contextual optimization papers in the ILO framework.}
\label{tab:ilo_methods}
\begin{adjustbox}{max width=\textwidth}
\begin{threeparttable}[ht]
    \begin{tabular}{l*{9}{E}} 
        \toprule
         & \multicolumn{4}{c}{Objective} & \multicolumn{2}{c}{Feasible domain} & \multicolumn{3}{c}{Training} \\
        \cmidrule(lr){2-5} \cmidrule(lr){6-7} \cmidrule(lr){8-10} & LP & QP & Convex & Non convex & Integer & Uncertain & {\,\,Implicit} diff. & Surr. loss & {\,\,Surr.} optim. \\
        \midrule
        \revised{\cite{amos_optnet_2017}} & \revised{\cross} & \revised{\tick} & \revised{\cross}& \revised{\cross}& \revised{\cross}& \revised{\tick}& \revised{\tick}& \revised{\cross}& \revised{\tick} \\
        \cite{donti_task-based_2019}    & \cross & \tick  & \tick& \cross & \cross & \tick & \tick & \cross& \cross \\\revised{\cite{Agrawal_Amos_Barratt_Boyd_Diamond_Kolter_2019}} & \revised{\cross} & \revised{\tick} & \revised{\tick}& \revised{\cross}& \revised{\cross}& \revised{\tick}& \revised{\tick}& \revised{\cross}& \revised{\tick} \\\cite{vlastelica2019diff}       & \tick  & \cross & \cross& \cross & \tick & \cross & \cross & \tick& \cross \\
        \cite{Wilder_Dilkina_Tambe_2019}& \tick  & \cross & \cross& \cross & \revised{\tick} & \cross & \tick & \cross& \cross \\
        \cite{Wilder2019graph}          & \cross & \tick  & \cross& \cross & \tick & \cross & \cross & \tick& \cross \\
        \cite{Berthet2020}              & \tick  & \cross & \cross& \cross & \tick & \cross & \cross &\cross& \tick \\
        \cite{elmachtoub_decision_2020} & \tick  & \cross & \cross& \cross & \cross& \cross & \cross & \tick& \cross \\
        \cite{ferber_mipaal_2020}       & \tick  & \cross & \cross& \cross & \tick & \cross & \tick & \cross& \cross \\
        \cite{Mandi_Guns_2020}          & \tick  & \cross & \cross& \cross & \tick & \cross & \tick & \cross& \cross \\
        \cite{Mandi2020_spo}            & \tick  & \cross & \cross& \cross & \tick & \cross & \cross & \tick& \cross \\
        \cite{Grigas2021integrated}  
        & \cross  & \cross& \tick & \cross & \cross &\cross & \cross &\cross& \tick \\                        \revised{\cite{mulamba2021contrastive}}   & \revised{\tick}  & \revised{\cross} & \revised{\cross}& \revised{\cross} & \revised{\tick}& \revised{\cross} & \revised{\cross} &\revised{\cross}& \revised{\tick} \\
        \cite{chung2022decision}        & \cross  & \cross & \tick & \cross & \tick & \cross & \cross & \tick& \cross \\
        \cite{Cristian2022projectnet}   & \cross  & \cross& \tick & \cross & \cross &\cross & \cross &\cross& \tick \\
        \cite{Dalle2022}                & \tick  & \cross & \cross& \cross & \tick & \cross & \cross &\cross& \tick \\
        \cite{elmachtoub_smart_2022}    & \tick  & \cross & \cross& \cross & \tick & \cross & \cross & \tick& \cross \\
        \cite{Jeong2022exact}           & \tick  & \cross & \cross& \cross & \tick & \cross & \cross & \tick& \cross \\
        \cite{kallus_stochastic_2022}   & \cross & \cross & \tick & \cross & \cross& \cross & \cross & \tick& \cross \\
        \cite{Kong2022energy}           & \cross  & \tick & \tick & \tick  & \tick & \cross & \cross &\cross& \tick \\
        \cite{lawless_note_2022}        & \tick  & \cross & \cross& \cross & \tick & \cross & \cross & \tick& \cross \\
        \cite{loke2022decision}         & \tick  & \cross & \cross& \cross &\cross & \cross & \cross & \tick& \cross \\
        \cite{Mandi2022learning2rank}   & \tick  & \cross & \cross& \cross & \tick & \cross & \cross &\cross& \tick \\
        \cite{Munoz2022}                & \tick  & \cross & \cross& \cross & \cross& \cross & \cross & \tick& \cross \\
        \cite{shah2022decision}         & \revised{\tick}  & \tick & \tick & \tick  & \tick & \cross & \cross &\cross& \tick \\
        \cite{Butler_Kwon_ADMM2023}     & \cross & \tick  & \cross& \cross & \cross & \cross & \tick & \cross& \cross \\
        \cite{Costa_Iyengar_2023}       & \cross  & \tick & \tick & \cross & \cross & \tick & \tick & \tick & \cross     \\
        \cite{Estes2023smart}           & \tick  & \cross & \tick & \cross & \cross& \cross & \cross & \tick& \cross \\
        \cite{kotary2023folded}         & \tick  & \tick  & \tick & \tick  & \cross & \cross & \tick & \cross& \cross \\
        \cite{McKenzie_Fung_Heaton_2023}& \tick  & \cross & \cross& \cross & \cross & \cross & \tick & \cross& \cross \\
        \revised{\cite{Sun_marginlp_2023}}         &\revised{\tick} &\revised{\cross}  & \revised{\cross} & \revised{\cross}  & \revised{\cross}  &\revised{\cross}  & \revised{\cross} &\revised{\tick} & \revised{\cross}  \\
        \cite{Sun_Shi_Wang2023}         & \cross & \tick  & \cross& \tick & \cross & \cross & \tick & \cross& \cross \\
        \bottomrule
    \end{tabular}
    \begin{tablenotes}
    Notes: an approach has a \quoteIt{Convex} objective if it can handle general convex objective functions that are not linear or quadratic such as convex piecewise-linear objective functions; an \quoteIt{Uncertain} feasible domain denotes that some constraints are subject to uncertainty. Implicit diff., surr.~loss and surr.~optim.~denote implicit differentiation, surrogate differentiable loss function and surrogate differentiable optimizer, respectively. \revised{This table covers papers that introduced an ILO framework to solve a class of CSO problems, papers focused on deriving theoretical guarantees and specific applications are therefore intentionally excluded. The tick marks reflect our understanding of the claimed scope of the contributions.}
   \end{tablenotes}
\end{threeparttable}
\end{adjustbox}
\end{table}

\subsection{Models}
\label{sec:ILOmodels}

\cite{Bengio_1997} appears to be the first to train a prediction model using a loss that is influenced by the performance of \revised{an action} prescribed by a conditional expected value-based decision rule. This was done in the context of portfolio management, where an investment decision rule exploits a point prediction of asset returns. Effective wealth accumulation is used to steer the predictor toward predictions that lead to good investments. More recent works attempt to integrate a full optimization model, rather than a rule, into the training pipeline. Next, we summarize how ILO is applied to the two types of contextual optimization models and introduce two additional popular task models that have been considered under ILO, replacing the traditional expected cost task.

\paragraph{Expected value-based model.} To this date, most of the literature has considered performing ILO on an expected value-based optimization model. Namely, following the notation presented in Definition~\ref{def:exp_model} (Section~\ref{sec:learnOptim}), this training pipeline is interested in the loss $\mathcal{L}(\btheta):=H(z^*(\cdot,g_\btheta),\Pemp)=\Expect_{\Pemp}[c(z^*(\x,g_\btheta),\y)]$ with $g_\btheta(\x)$ as a point predictor for $\y$, which we interpret as a prediction of $\Expect[\y|\x]$. This already raises challenges related to the non-convexity of the integrated loss function $\mathcal{L}(\btheta)$ 
and its differentiation with respect to $\btheta$:
\begin{align*}
    \nabla_\btheta \mathcal{L}(\btheta) & =\frac{1}{N}\sum_{i=1}^N\nabla_\btheta c(z^*(\x_i,g_\btheta),\y_i) \\
    & = \frac{1}{N}\sum_{i=1}^N\sum_{j=1}^{\revised{d_{\z}}}\sum_{k=1}^{\revised{d_{\y}}}\frac{\partial c(z^*(\x_i,g_\btheta),\y_i)}{\partial z_j}\left.\frac{\partial z^*_j(\x_i,\hat{\y})}{\partial \hat{y}_k}\right|_{\hat{\y}=g_\theta(\x_i)}\nabla_\btheta [g_\btheta(\x_i)]_k
\end{align*}
with $\frac{\partial z^*_j(\x_i,\hat{\y})}{\partial \hat{y}_k}$ as the most problematic evaluation. For instance, when $z^*(\x_i,g_\btheta)$ is the solution of a linear program (LP), it is well known that its gradient is either null or non-existent as it jumps between extreme points of the feasible polyhedron as the objective is perturbed.

\paragraph{Conditional distribution-based model.}
In the context of learning a conditional distribution model $f_\theta(\x)$, \cite{donti_task-based_2019} appear to be the first to study the ILO problem. They model the distribution of the uncertain parameters using parametric distributions (exponential and normal). For the newsvendor problem, it is shown that the ILO model outperforms decision rule optimization with neural networks and SLO with maximum likelihood estimation~(MLE) when there is model misspecification. Since then, it has become more common to formulate the CSO problem as a weighted SAA model (as discussed in Section~\ref{sec:wSAA}). The prediction model then amounts to identifying a vector of weights to assign to each historical sample. \revised{This approach is taken by \citet{kallus_stochastic_2022}, who train a random forest regressor in an integrated fashion to assign weights, and by \citet{Grigas2021integrated}, who show how to train general differentiable models to predict the probabilities of an uncertain parameter $\y$ with finite support.}

\paragraph{Regret minimization task.} A recent line of work has tackled the ILO problem from the point of view of regret. Indeed, in \citet{elmachtoub_smart_2022}, a contextual point predictor $g_{\btheta}(\x)$ is learned by minimizing the regret associated with implementing the prescribed \revised{action} based on the mean estimator $g_{\btheta}(\x)$ instead of based on the realized parameters $\y$ (a.k.a. the optimal hindsight or wait-and-see decision). Specifically, the value of an expected value-based policy $\pi_\btheta(\x):=z^*(\x,g_{\btheta})$ is measured as the expected regret defined as:
\begin{equation}
    \label{eq:regret}
    H_{\mymbox{Regret}}(\pi_\btheta,\Prob):=\Expect_{\Prob}[c(\pi_\btheta(\x),\y)- c(z^*(\x,\y),\y)].
\end{equation}
Minimizing the expected regret returns the same optimal parameter vector $\btheta$ as the ILO problem \eqref{eq:integ_opt:objt}. This is due to the fact that:
\[ H_{\mymbox{Regret}}(\pi,\Pemp)=\Expect_{\Pemp}[c(\pi(\x),\y)- c(z^*(\x,\y),\y)]=H(\pi,\Pemp)- \Expect_{\Pemp}[c(z^*(\x,\y),\y)].\]
Hence, both $H_{\mymbox{Regret}}(\pi,\Pemp)$ and $H(\pi,\Pemp)$ have the same set of minimizers.
    
\paragraph{Optimal action imitation task.}
ILO has some connections to inverse optimization, i.e., the problem of learning the parameters of an optimization model given data about its optimal solution \citep[see][where both problems are addressed using the same method]{Sun_marginlp_2023}. Indeed, one can replace the original objective of ILO with an objective that seeks to produce a $z^*(\x,f_\btheta)$ that is as close as possible to the optimal hindsight action and, therefore, closer to the regret objective. Specifically, to learn a policy that \quoteIt{imitates} the optimal hindsight action, one can first augment the data set with $\z_i^*:=z^*(\x_i,\y_i)$ to get $\{(\x_i,\y_i,\z_i^*)\}_{i=1}^N$. 
Thereafter, a prediction model $f_\btheta(\x)$ is learned in a way that the \revised{action} $z^*(\x_i,f_\btheta)$ is as close as possible to $\z_i^*$ for all samples in the training set \citep{Kong2022energy}: 
\begin{equation}
    H_{\mymbox{Imitation}}(\pi,\Pemp'):=\Expect_{\Pemp'}[d(\pi(\x),\z^*)]=\Expect_{\Pemp}[d(\pi(\x),z^*(\x,\y))]    
\end{equation}
with $\Pemp'$ as the empirical distribution on the lifted tuple $(\x,\Vxi,z^*(\x,\Vxi))$ based on the augmented data set and a distance function $d(\z,\z^*)$. We note that there is no reason to believe that the best imitator under a general distance function, e.g., $\|\z-\z^*\|_2$, performs well under our original metric $H(\pi,\Pemp)$. One exception is for $d(\z,\z^*):=c(\z,\Vxi)-c(\z^*,\Vxi)$, where we allow the distance to also depend on $\Vxi$, for which we recover the regret minimization approach, and therefore the same solution as with $H(\pi,\Pemp)$. Readers that have an interest in general inverse optimization methods should consult \cite{Chan2021inverse} for an extensive recent review of the field.

\subsection{Training by unrolling} \label{sec:unrolling}

An approach to obtain the Jacobian matrix $\frac{\partial z^*(\x,\hat{\y})}{\partial \hat{\y}}$ is unrolling \citep{domke2012generic}, which involves approximating the optimization problem with an iterative solver (e.g., first-order gradient-based method). Each operation is stored on the computational graph, which then allows, in principle, for computing gradients through classical back-propagation methods. Unfortunately, this approach requires extensive amounts of memory. 
Besides this, the large size of the computational graph exacerbates the vanishing and exploding gradient problems typically associated with training neural networks \citep{monga2021algorithm}.

\subsection{Training using implicit differentiation}\label{sec:implicitdiff}
Implicit differentiation allows for a memory-efficient backpropagation as opposed to unrolling \citep[we refer to][for discussion on training constant memory implicit models using a fixed-point -- FP -- equation and feedforward networks of infinite depths]{Bai_Kolter_Koltun_2019}. \cite{amos_optnet_2017} appear to be the first to have employed implicit differentiation methods to train an ILO model, which they refer to as \model{OptNet}. They consider expected value-based optimization models that take the form of constrained quadratic programs~(QP) with equality and inequality constraints. They show how the implicit function theorem \revised{\citep[IFT --][]{Halkin1974}} can be used to differentiate $z^*(\x,g_\btheta)$ with respect to $\btheta$ using the Karush–Kuhn–Tucker (KKT) conditions that are satisfied at optimality. Further, they provide a custom solver based on a primal-dual interior method to simultaneously solve multiple QPs on GPUs in batch form, permitting 100-times speedups compared to Gurobi and CPLEX. This approach is extended to conditional stochastic and strongly convex optimization models in \cite{donti_task-based_2019}. They use sequential quadratic programming (\model{SQP}) to obtain quadratic approximations of the objective functions of the convex program at each iteration until convergence to the solution and then differentiate the last iteration of \model{SQP} to obtain the Jacobian. For a broader view of implicit differentiation, we refer to the surveys by \cite{duvenaud2020deep} and \cite{Blondel_Berthet_Cuturi2022}.

To solve large-scale QPs with linear equality and box inequality constraints, \cite{Butler_Kwon_ADMM2023} use the ADMM algorithm to decouple the differentiation procedure for primal and dual variables, thereby decomposing the large problem into smaller subproblems. Their procedure relies on implicit differentiation of the FP equations of the alternating direction method of multipliers~(ADMM) algorithm (\model{ADMM-FP}). They show that unrolling the iterations of the ADMM algorithm on the computational graph \citep{sun2016deep,xie2019differentiable} results in higher computation time than \revised{\model{ADMM-FP}}. Their empirical results on a portfolio optimization problem with $254$ assets suggest that computational time can be reduced by a factor of almost five by using \model{ADMM-FP} compared to \model{OptNet}, mostly due to the use of the ADMM algorithm in the forward pass. Note that the experiments in \cite{Butler_Kwon_ADMM2023} were conducted on a CPU.

To extend \model{OptNet} to a broader class of problems, \cite{Agrawal_Amos_Barratt_Boyd_Diamond_Kolter_2019} introduce \model{cvxpylayers} that relies on converting disciplined convex programs in the domain-specific language used by CVXPY into conic programs. They implicitly differentiate the residual map of the homogeneous self-dual embedding associated with the conic program.

\cite{McKenzie_Fung_Heaton_2023} note that using KKT conditions for constrained optimization problems with DNN-based policies is computationally costly as \quoteIt{\model{cvxpylayers} struggles with solving problems containing more than 100 variables} \citep[see also][]{Butler_Kwon_ADMM2023}. An alternative is to use projected gradient descent (\model{PGD}) where DNN-based policies are updated using an iterative solver and projected onto the constraint set $\Z$ at each iteration and the associated FP system \citep{Donti_Roderick_Fazlyab_Kolter_2021, Chen_Donti_Baker_Kolter_Berges_2021, Blondel_Berthet_Cuturi2022} is used to obtain the Jacobian.

Since a closed-form solution for the projection onto $\Z$ is unavailable in many cases, the projection step may be costly, and in some cases, \model{PGD} may not even converge to a feasible point \citep{Rychener_Sutter_2023}. To avoid computing the projection in the forward pass, \cite{McKenzie_Fung_Heaton_2023} solve the expected value-based CSO problem using Davis-Yin operator splitting \citep{davis2017threeoperator} while the backward pass uses the Jacobian-free backpropagation \citep{Fung_Heaton_Li_Mckenzie_Osher_Yin_2022} in which the Jacobian matrix is replaced with an identity matrix \revised{\cite[see also][ where a similar approach is used for expected value-based models]{Sahoo2023backpropagation}}.

To mitigate the issues with unrolling, \cite{kotary2023folded} propose FP folding (\model{fold-opt}) that allows analytically differentiating the FP system of general iterative solvers, e.g., \model{ADMM}, \model{SQP}, and \model{PGD}. By unfolding (i.e., partial unrolling), some of the steps of unrolling are grouped in analytically differentiable update function $\mathcal{T}: \Re^{\revised{d_{\y}}}\rightarrow\Re^{\revised{d_{\y}}}$:
\[ z_{k+1}(\x, \hat{\y}) = \mathcal{T}(z_k(\x, \hat{\y}), \hat{\y}).\]
Realizing that $z^*(\x,\hat{\y})$ is the FP of the above system, 
they use the IFT to obtain a linear system (a differential FP condition) that can be solved to obtain the Jacobian. This effectively decouples the forward and backward pass enabling the use of black box solvers like Gurobi for the forward pass while \model{cvxpylayers} is restricted to operator splitting solvers like ADMM. An added benefit of using \model{fold-opt} is that it can solve non-convex problems. In the case of portfolio optimization, the authors note that the superior performance of their model with respect to \model{cvxpylayers} can be explained by the precise calculations made in the forward pass by Gurobi.

While speedups can be obtained for sparse problems, \cite{Sun_Shi_Wang2023} remark that the complexity associated with differentiating the KKT conditions is cubic in the total number of decision variables and constraints in general. They propose an alternating differentiation framework (called \model{Alt-Diff}) to solve parameterized convex optimization problems with polyhedral constraints using ADMM that decouples the objective and constraints. This procedure results in a smaller Jacobian matrix when there are many constraints since the gradient computations for primal, dual, and slack variables are done alternatingly. The gradients are shown to converge to those obtained by differentiating the KKT conditions. The authors employ truncation of iterations to compensate for the slow convergence of ADMM when compared to interior-point methods and provide theoretical upper bounds on the error in the resulting gradients. \model{Alt-Diff} is shown to achieve the same accuracy with truncation and lower computational time when compared to \model{cvxpylayers} for an energy generation scheduling problem.

Motivated by \model{OptNet}, several extensions have been proposed to solve linear and combinatorial problems. \citet{Wilder_Dilkina_Tambe_2019} solve LP-representable combinatorial optimization problems and LP relaxations of combinatorial problems during the training phase. Their model, referred to as \model{QPTL} (Quadratic Programming Task Loss), adds a quadratic penalty term to the objective function of the linear problem. This has two advantages: it recovers a differentiable linear-quadratic program, and the added term acts as a regularizer, which might avoid overfitting. To solve a general mixed-integer LP (MILP), \citet{ferber_mipaal_2020} develop a cutting plane method \model{MIPaal}, which adds a given number of cutting planes in the form of constraints $S\z\leq \V{s}$ to the LP relaxation of the MILP. Instead of adding a quadratic term, \cite{Mandi_Guns_2020} propose \model{IntOpt} based on the interior point method to solve LPs that adds a log barrier term to the objective function and differentiates the homogeneous self-dual formulation of the LP. Their experimental analyses show that this approach performs better on energy cost-aware scheduling problems than \model{QPTL} using the data from \cite{ifrim2012properties}.

\revised{\cite{Costa_Iyengar_2023}} introduce an ILO framework with the \revised{weighted average of} Sharpe ratio \revised{and MSE loss} as a task loss and replace the optimization problem with a surrogate DRO problem. By using convex duality, they reformulate the minimax problem as a minimization problem and learn the parameters (e.g., size of ambiguity set) using implicit differentiation instead of cross-validation~(CV). More specifically, the DRO model uses a deviation risk measure (e.g., variance) to control variability in the portfolio returns associated with the prediction errors $\V{\epsilon}_i = \y_i - g_\theta(\x_i)$:
\[ \argmin_{\z} \max_{\mathbb{Q} \in \mathcal{B}_{r}^\phi(\Pemp)} \revised{\Expect_{\mathbb{Q}}\left[\left(\V{\epsilon}^\top \z - \Expect_{\mathbb{Q}}[\V{\epsilon}^\top \z]  \right)^2\right]},\]
where the distribution of errors lies in $\phi$-divergence (e.g., Hellinger distance) based ambiguity set %
$\mathcal{B}^\phi_r(\Pemp) = \{ \mathbb{Q}:  \Expect_{\hat{\Prob}}[ \phi( \mathbb{Q}/ \hat{\Prob})] \leq r \}$ \revised{centered at  $\Pemp = \frac{1}{N}\sum_{i=1}^N \delta_{\V{\epsilon}_i}$} .

\revised{For convex problems, the optimality conditions are given by KKT conditions, which can be represented as $F(\btheta, \z)=0$ where $F:\mathbb{R}^{d_{\x}} \times \mathbb{R}^{d_{\z}} \rightarrow \mathbb{R}^{m}$, where $m$ is proportional to the number of constraints that define $\Z$. From the classical IFT \citep{dontchev2009implicit}, we know that if $F$ is continuously differentiable and the Jacobian matrix  with respect to $\z$, denoted by $\nabla_{\z} F(\btheta, \z(\btheta))$, is non-singular at the point $(\bar{\btheta}, \bar{\z})$, then there exists a neighborhood around $\bar{\btheta}$ for which the gradient of the optimal solution with respect to the parameters is given by:
\[\frac{\partial\z^*(\btheta)}{\partial \btheta}= -(\nabla_{\z} F(\btheta, \z(\btheta)))^{-1} \nabla_{\btheta} F(\btheta, \z(\btheta)).\]   
When the Jacobian matrix $\nabla_{\z} F(\btheta, \z(\btheta))$ is singular, classical IFT cannot be applied. This occurs in linear programs and can also arise in smooth QPs as shown in \cite{bolte2021nonsmooth}. 
\cite{bolte2021nonsmooth} obtain a generalization of IFT to non-smooth functions using conservative Jacobians that generalize Clarke Jacobians \citep{clarke1990optimization} for locally Lipschitz function $F$. They also derive conservative Jacobians for conic optimization layers  \citep{Agrawal_Amos_Barratt_Boyd_Diamond_Kolter_2019}. 

Further, \cite{bolte2021nonsmooth} illustrate using \model{cvxpylayers} that in a bilevel program which is a composition of a quadratic function with the solution map of a linear program, gradient descent does not converge but gets stuck in a ``limit cycle of non-critical points'' even though invertibility condition does not hold only on a set of measure $0$ (defined by a line) where the solution map moves from extreme point to another. As this example illustrates, the convergence of gradient methods based on IFT can be impacted by the non-invertibility of the Jacobian matrix and non-smoothness which is difficult to verify a priori. As a result, research efforts have been directed toward designing surrogate loss functions and perturbation-based models for CSO problems that could circumvent the need to use the IFT.}

\subsection{Training using a surrogate differentiable loss function}
\label{sec:surrogateDiffLoss}
As discussed in Section~\ref{sec:ILOmodels}, minimizing directly the task loss in~\eqref{eq:integ_opt:objt} or the regret in~\eqref{eq:regret} is computationally difficult in most cases. For instance, the loss may be piecewise-constant as a function of the parameters of a prediction model and, thus, may have no informative gradient. To address this issue, several surrogate loss functions with good properties, e.g., differentiability and convexity, have been proposed to train ILO models.

\subsubsection{\model{SPO+}.}
In \revised{CSO problems, \cite{elmachtoub_smart_2022} first tackle the potential non-uniqueness of $z^*(\x,g_{\btheta})$ by introducing a Smart \quoteIt{Predict, then Optimize} (\model{SPO}) model where the decision-maker chooses} to minimize the empirical average of the regret under the worst-case optimal solution as defined below:
\begin{equation}
    \begin{split}
    \text{(\model{SPO})}\quad \min_{\btheta} \max_{\pi} & \, H_{\mymbox{Regret}}(\pi,\Pemp), \\
     \text{s.t. } & \, \pi(\x)\in\argmin_{\z\in\Z}c(\z,g_{\btheta}(\x)), \, \forall \x .
    \end{split}
\end{equation}
In the expected value-based model, they show that the \model{SPO} objective reduces to training the prediction model according to the ERM problem:
\begin{align*}
   \btheta^\star \in\argmin_{\btheta} \rho_{{\mymbox{\model{SPO}}}}(g_{\btheta}, \Pemp):=  \Expect_{\Pemp} \big[\ell_{\mymbox{\model{SPO}}}(g_{\btheta}(\x), \y)\big],
\end{align*}
with:
\begin{equation*}
    \ell_{\mymbox{\model{SPO}}}(\hat{\y}, \y):=\sup_{\bar{\z}\in\argmin_{\z\in\Z} c(\z,\hat{\y})}c(\bar{\z},\y)- c(z^*(\x,\y),\y).
\end{equation*}

Since the \model{SPO} loss function is nonconvex and discontinuous in $\hat{\y}$ \citep[Lemma 1]{ho-nguyen_risk_2022}, \cite{elmachtoub_smart_2022} focus on the linear objective $c(\z,\y):=\y^T\z$ and replace the \model{SPO} loss with a convex \revised{upper bound} called \model{SPO+}:
\begin{equation*}
    \ell_{\mymbox{\model{SPO}$+$}}(\hat{\y}, \y):=\sup_{\z\in\Z} (\y-2\hat{\y})^T\z+2 \hat{\y}^T z^*(\x,\y)-\y^T z^*(\x,\y),
\end{equation*}
which has a closed-form expression for its subgradient
\begin{equation}
   2\left(z^*(\x,\Vxi) - z^*(\x,2\hat{\Vxi}-\Vxi)\right)\in\nabla_{\hat{\Vxi}}\ell_{\mymbox{\model{SPO}$+$}}(\hat{\Vxi},\Vxi). \label{eq:grad:spo}
\end{equation}

\cite{loke2022decision} propose a decision-driven regularization model~(\model{DDR}) that combines prediction accuracy and decision quality in a single optimization problem with loss function as follows:
\begin{align*}
    \ell_{\mymbox{\model{DDR}}}(\hat{\y},\y) =   d(\hat{\y},\y)-\lambda \min_{\z \in \Z} \{ \mu \y^\top \z +(1-\mu) \hat{\y}^\top \z \},
\end{align*}
and \model{SPO+} being a special case with $\mu=-1,\; \lambda=1,$ and $d(\hat{\y},\y) = 2\hat{\y}^\top z^*(\x, \y)-\y^T z^*(\x,\y)$.

\paragraph{\model{SPO+} for combinatorial problems.}
Evaluating the gradient of \model{SPO+} loss in \eqref{eq:grad:spo} requires solving the optimization problem \eqref{eq:dp} to obtain $z^*(\x,2\hat{\y}-\y)$ for each data point. This can be computationally demanding when the optimization model in \eqref{eq:dp} is an NP-hard problem.
\cite{Mandi2020_spo} propose a \model{SPO-relax} approach that computes the gradient of \model{SPO+} loss by solving instead a continuous relaxation when \eqref{eq:dp} is a MILP. They also suggest speeding up the resolution using a warm-start for learning with a pre-trained model that uses MSE as the loss function. Another way proposed to speed up the computation is warm-starting the solver. For example, $z^*(\x,\y)$ can be used as a starting point for MILP solvers or to cut away a large part of the feasible space. \cite{Mandi2020_spo} show that for weighted and unweighted knapsack problems as well as energy-cost aware scheduling problems \citep[CSPLib, Problem 059,][]{csplib059}, \model{SPO-relax} results in faster convergence and similar performance compared to \model{SPO+} loss. Also, \model{SPO-relax} provides low regret solutions and faster convergence compared to \model{QPTL} in the aforementioned three problems, except in the weighted knapsack problem with low capacity. 

With a focus on exact solution approaches, \cite{Jeong2022exact} study the problem of minimizing the regret in \eqref{eq:regret} assuming a linear prediction model $g_{\btheta}(\x) = \btheta \x$ with $\btheta\in\Re^{\revised{d_{\z}}\times \revised{d_{\x}}}$. Under the assumption that $z^*(\x,g_{\btheta})$ is unique for all $\btheta$ and $\x$, the authors reformulate the bilevel \model{SPO} problem as a single-level MILP using symbolic variable elimination. They show that their model can achieve up to two orders of magnitude improvement in expected regret compared to \model{SPO+} \revised{on the training set}.
\cite{Munoz2022} applies a similar idea of representing the set of optimal solutions with a MILP. They rely on the KKT conditions of the problem defining $z^*(\x,g_{\btheta})$ to transform the bilevel integrated problem into a single-level MILP. 
Finally, \citet{Estes2023smart} use the \model{SPO} loss function to solve a two-stage LP with right-hand side uncertainty. They propose a lexicographical ordering rule to select the minimal solution when there are multiple optima and approximate the resulting piecewise-linear loss function, \model{lex-SPO}, by a convex surrogate to find the point predictor.

\paragraph{\model{SPO} Trees.} 
\cite{elmachtoub_decision_2020} propose a model (\model{SPOT}) to construct decision trees that segment the \revised{covariates} based on the \model{SPO} loss function while retaining the interpretability in the end-to-end learning framework. Their model outperforms classification and regression trees (CART) in the numerical experiments on a news recommendation problem using \revised{a real-world} dataset and on the shortest path problem with synthetic data \revised{(also} used in \cite{elmachtoub_smart_2022}\revised{)}.

\paragraph{Guarantees.}
\cite{elmachtoub_smart_2022} show that under certain conditions, the minimizers of the \model{SPO} loss, \model{SPO+} loss and MSE loss are almost always equal to $\Expect_{\Prob(\y \rvert \x)}[\y]$ given that $\Expect_{\Prob(\y \rvert \x)}[\y]\in \mathcal{H}$. Thus, \model{SPO+} is Fisher consistent (see Definition \ref{def:fisher} in Appendix~\ref{sec:appen:definitions}) with respect to \revised{the} \model{SPO} loss. This means that minimizing the surrogate loss also minimizes the true loss function. \cite{ho-nguyen_risk_2022} show that for \revised{some examples of} a multiclass classification problem, \model{SPO+} is Fisher inconsistent, while MSE loss is consistent. However, complete knowledge of the distribution is a limitation in practice where the decision-maker has access to only the samples from the distribution. As a result, \cite{ho-nguyen_risk_2022} and \cite{liu2021risk} provide calibration bounds that hold for a class of distributions $\mathcal{D}$ on $\X \times \Y$ and ensure that a lower excess risk of predictor for MSE and \model{SPO+}, respectively,  translates to lower excess \model{SPO} risk (see Definition \ref{def:uniform_calibration} in Appendix~\ref{sec:appen:definitions}).

In many ML applications, one seeks to derive finite-sample guarantees, which are given in the form of a generalization bound, i.e., an upper bound on the difference between the true risk of a loss function and its empirical risk estimate for a given sample size $N$. A generalization bound for the \model{SPO} loss function is given in \cite{El_Balghiti_Elmachtoub_Grigas_Tewari_2022} \citep[extension of][]{El_Balghiti_Elmachtoub_Grigas_Tewari_2019} based on Rademacher complexity (see Definition \ref{def:rademacher} in Appendix \ref{sec:appen:definitions}) of the \model{SPO} loss composed with the prediction functions $g_\btheta \in \mathcal{H}$. More specifically, the bound achieved in \cite{El_Balghiti_Elmachtoub_Grigas_Tewari_2019} is $\text{O}\left(\sqrt{\frac{\log(N)}{N}}\right)$, and tighter bounds with respect to \revised{action} and feature dimension are obtained using \model{SPO} function's structure and if $\Z$ satisfies a \quoteIt{strength} property. \cite{Hu_Kallus_Mao_2022} show that for linear CSO problems, the generalization bound for MSE loss and \model{SPO} loss is $\text{O}\left(\sqrt{\frac{1}{N}}\right)$ while faster convergence rates for the SLO model compared to ILO model are obtained under certain low-noise assumptions. \cite{Elmachtoub_Lam_Zhang_Zhao_2023} show that for non-linear optimization problems, SLO models stochastically dominate ILO in terms of their asymptotic optimality gaps when the hypothesis class covers the true distribution. When the model is misspecified, they show that ILO outperforms SLO \revised{asymptotically in a general nonlinear setting}. 

\subsubsection{Surrogate loss for a stochastic forest.}

\cite{kallus_stochastic_2022} propose an algorithm called \model{StochOptForest}, which generalizes the random-forest based local parameter estimation procedure in \cite{Athey_Tibshirani_Wager_2019}. A second-order perturbation analysis of stochastic optimization problems allows them to scale to larger CSO problems since they can avoid solving an optimization problem at each candidate split. The policies obtained using their model are shown to be asymptotically consistent, and the benefit of \revised{ILO} is illustrated by comparing their approach to the random forests of \citet{Bertsimas_Kallus_2020} on a set of problems with synthetic and real-world data.

\subsubsection{Other surrogates.}
\cite{Wilder2019graph} introduce \model{ClusterNet} to solve hard combinatorial graph optimization problems by learning incomplete graphs. The model combines graph convolution networks to embed the graphs in a continuous space and uses a soft version of k-means clustering to obtain a differential proxy for the combinatorial problems, e.g., community detection and facility location. Numerical experiments on a synthetic data set show that \model{ClusterNet} outperforms the two-stage SLO approach of first learning the graph and then optimizing, as well as other baselines used in community detection and facility location.

Focusing on combinatorial problems, \cite{vlastelica2019diff} propose a differentiable black-box (\model{DBB}) approach to
tackle the issue that the Jacobian of $z^*(\x,g_\btheta)$ is zero almost everywhere by approximating the true loss function using an interpolation controlled in a way that balances between \quoteIt{informativeness of the gradient} and \quoteIt{faithfulness to the original function}.
Algorithmically, this  is done by perturbing the prediction $g_\btheta(\x)$ in the direction $\nabla_{\z} c(z^*(\x,g_\btheta),\y)$ and obtaining a gradient of the surrogate loss based on the effect of this perturbation on the resulting perturbed action.

\revised{\cite{chung2022decision} introduce a computationally tractable ILO model to solve non-linear CSO problems. Using the first-order Taylor expansion of the task loss around the prediction, they introduce a reweighted MSE loss function where weights are determined by taking the gradient of task loss with respect to the prediction. To solve a large-scale multi-facility inventory allocation problem with few samples for each facility, they use a single random forest that can predict the demand across facilities and products. Assuming that each tree in the random forest provides an independent and identically distributed realization of the uncertain parameter, they obtain the conditional distribution of uncertain parameter, $f_{\btheta_0} = \frac{1}{T}\sum_{t=1}^{T} \Dirac_{\hat{\y}^t}$, where $\hat{\y}^t$ is the prediction of tree $t$. For each feature $\x_i$ and conditional distribution $f_{\btheta_0}$, they obtain an  optimal allocation, $z_i^j = z^*(\x,f_{\btheta_0})$ for facility $j$ that minimizes the average unmet demand.  In the last step, they retrain the random forest  to minimize the reweighted MSE loss function:
\begin{align}
    \label{eq:chung_loss2}
   \argmin_{\btheta} \sum_{i=1}^N\sum_{j=1}^M \sum_{t=1}^T \Indc\left[\hat{y}^{t,j}_i \geq s_i^j + z_i^j \right]  \; \lvert f_{\btheta}(\x_i)-\hat{y}^{t, j}_i\rvert,
\end{align}
where $M$ is the number of facilities, $\hat{y}^{t,j}_i$ and  $s^j_i$ denote the demand and inventory levels, respectively, at facility $j$. The above model \eqref{eq:chung_loss2} solves the optimization problem once during training, and is shown to be scalable for a medical allocation problem in Sierra Leone when compared to \citet{kallus_stochastic_2022} where splitting of the feature space is done based on the task loss.
}

\revised{\cite{lawless_note_2022} introduce a loss function similar to \citet{chung2022decision} that weighs} the prediction error with a regret term as follows:
\begin{equation}
\label{eq:weightedregret}
    d(g_\btheta(\x), \y) =  [c(z^*(\x,g_\btheta),\y)- c(z^*(\x,\y),\y)] (\y-g_\btheta(\x))^2.
\end{equation}
Learning optimal $\btheta$ from the above formulation involves an $\argmin$ differentiation. So, the authors provide a two-step polynomial time algorithm to approximately solve the above problem. It first computes a pilot estimator \revised{$g_{{\btheta}_0}$} by solving \eqref{eq:regression} with $d(g_{\btheta}(\x), \y) = (g_{\btheta}(\x) - \y)^2 $ and then solving \eqref{eq:regression} with the distance function in \eqref{eq:weightedregret} where \revised{$c(z^*(\x,g_{\btheta}),\y)$} is substituted with \revised{$c(z^*(\x,g_{{\btheta}_0}),\y)$}. The authors show that their simple algorithm performs comparably to {\model{SPO+}}.

We conclude this subsection on surrogate loss functions by mentioning the efforts in \cite{Sun_marginlp_2023} to learn a cost point estimator (in an expected value-based model) to imitate the hindsight optimal solution. This is done by designing a surrogate loss function that penalizes how much the optimal basis optimality conditions are violated. They derive generalization error bounds for this new loss function and employ them to provide a bound on the sub-optimality of the minimal $\btheta$.

\subsection{Training using a surrogate differentiable optimizer}
\label{sec:surrogateDiffOpt}
\subsubsection{Differentiable perturbed optimizer.} 
One way of obtaining a differentiable optimizer is 
to apply a stochastic perturbation to the parameters predicted by the ML model. Taking the case of expected value-based models as an example, the key idea is that although the gradient of the solution of the \revised{contextual} problem with respect to the predicted parameters $\hat{\y}:=g_{\btheta}(\x)$ is zero almost everywhere, if we perturb the predictor using a noise with differentiable density, then the expectation of the solution of the perturbed contextual problem, 
\[ \bar{z}^\varepsilon(\x,g_{\btheta}) = \Expect_{\Psi}[\tilde{z}^\varepsilon(\x, g_{\btheta},\Psi)] \mbox{ with }\tilde{z}^\varepsilon(\x, g_{\btheta},\Psi) := \argmin_{\z \in \Z} c(\z,g_{\btheta}(\x) + \varepsilon \Psi)  ,\]
where $\varepsilon>0$ controls the amount of perturbation, and more generally of the expected cost of the associated random policy $\Expect_\Psi[H(\tilde{z}^\varepsilon(\cdot,g_\btheta,\Psi),\Pemp)]$ can be shown to be smooth and differentiable. This idea is proposed and exploited in \citet{Berthet2020}, which focus on a bi-linear cost $c(\z,\Vxi):=\Vxi^T\z$ thus simplifying $\Expect_\Psi[H(\tilde{z}^\varepsilon(\cdot,g_\btheta,\Psi),\Pemp)]=H(\bar{z}^\varepsilon(\cdot,g_\btheta),\Pemp)$. Further, they show that when an imitation ILO model is used with a special form of Bregman divergence to capture the difference between $z^*(\x,\y)$ and $\tilde{z}^\varepsilon(\x,\hat{\y},\Psi)$, the gradient of $H_{\mymbox{Imitation}}(\tilde{z}^\varepsilon(\cdot,g_\btheta,\Psi),\Pemp')$ can be computed directly without needing to determine the Jacobian of $\bar{z}^\varepsilon(\x, g_{\btheta})$ \citep{Blondel2020learning}:
\[    H_{\mymbox{Imitation}}(\tilde{z}^\varepsilon(\cdot,g_\btheta,\Psi),\Pemp'):= \Expect_{\Pemp}[\ell_{\model{PFYL}}(g_\btheta(\x), \y)]\]
where $\ell_{\model{PFYL}}$ is a \revised{perturbed} Fenchel-Young loss \revised{(\model{PFYL})} given by:
\[\ell_{\model{PFYL}}(\hat{\y}, \y) := \hat{\y}^Tz^*(\x,\y) - \Expect_\Psi[ (\hat{\y}+\varepsilon\Psi)^T\tilde{\z}^\varepsilon(\x, \hat{\y},\Psi)]  + \varepsilon \Omega_{\model{PFYL}}(z^*(\x,\y)),\]
and $\Omega_{\model{PFYL}}(\z)$ is the Fenchel dual of $F(\y):=- \Expect_\Psi[ (\y+\Psi)^T\tilde{z}^\varepsilon(\x, \y,\Psi)] $. 
The gradient of the Fenchel-Young loss with respect to the model prediction is given by:
\[    \nabla_{\hat{\y}} \ell_{\model{PFYL}}(\hat{\y}, \y) = z^*(\x,\y) - \bar{\z}^\varepsilon(\x, \hat{\y}).\]
Evaluating this gradient can be done through Monte Carlo evaluations by sampling perturbations and solving the corresponding perturbed problems.

\citet{Dalle2022} introduce a multiplicative perturbation with the advantage that it preserves the sign of $g_{\btheta}(\x)$ without adding any bias: 
\[    \tilde{\z}^\varepsilon(\x, g_{\btheta},\Psi) := \argmin_{\z \in \Z} c(\z,g_{\btheta}(\x) \odot \exp(\varepsilon \Psi - \varepsilon^2/2)), \]
where $\odot$ is the Hadamard dot-product and the exponential is taken elementwise. \citet{Dalle2022} and \citet{Sun2023unified} also show that there is a one-to-one equivalence between the perturbed optimizer approach and using a regularized randomized version of the CSO problem for combinatorial problems with linear objective functions. 
Finally, \cite{Dalle2022} show an intimate connection between the perturbed minimizer approach proposed by \cite{Berthet2020} and surrogate loss functions approaches such as \model{SPO+} by casting them as special cases of a more general surrogate loss formulation.

\cite{mulamba2021contrastive} and \cite{Kong2022energy} consider an \quoteIt{energy-based} perturbed optimizer defined by its density of the form:
\begin{align}
    \tilde{\z}^\varepsilon(\x,f_\btheta)\sim \frac{\exp(- h(\z,f_\btheta(\x))/\varepsilon)}{\int \exp(- h(\z',f_\btheta(\x))/\varepsilon) d\z'},
    \label{eq:energyPerturb}
\end{align}
with $\varepsilon=1$, in the context of an imitation ILO problem. This general form of perturbed optimizer captures a varying amount of perturbation through $\varepsilon$, with $\tilde{\z}^\varepsilon(\x,f_\btheta)$ converging in distribution to $z^*(\x,f_\btheta)$ as $\varepsilon$ goes to zero. They employ the negative log-likelihood to measure the divergence between $\tilde{\z}^\varepsilon(\x,f_\btheta)$ and the hindsight optimal solution $z^*(\x,\y)$. Given the difficulties associated with calculating the partition function in the denominator of \eqref{eq:energyPerturb}, \cite{mulamba2021contrastive} devise a surrogate loss function based on noise-contrastive estimation, which replaces likelihood with relative likelihood when compared to a set of sampled suboptimal solutions. This scheme is shown to improve the performance over \model{SPO+} and \model{DBB} in terms of expected regret performance for linear combinatorial CSO.

Based on the noise contrastive estimation approach of \cite{mulamba2021contrastive}, \cite{Mandi2022learning2rank} note that ILO for combinatorial problems can be viewed as a learning-to-rank problem. They propose surrogate loss functions, with closed-form expressions for gradients, that are used to train to rank feasible points in terms of performance on the downstream optimization problem. Unlike \cite{mulamba2021contrastive}, \cite{Kong2022energy} tackles the partition function challenge by employing a self-normalized importance sampler that provides a discrete approximation. To avoid overfitting, the authors also introduce a regularization that penalizes the KL divergence between the perturbed optimizer distribution and a subjective posterior distribution over perturbed optimal hindsight actions $\Prob(\tilde{\z}^\varepsilon(\x,\y)\rvert \Vxi)$:
\begin{align*}
    &H_{\mymbox{Imitation}}(\tilde{\z}^\varepsilon(\cdot,f_\btheta),\Pemp'):=\\
    &\qquad\qquad\qquad-\Expect_{\Pemp}[\log(\Prob(\tilde{\z}^\varepsilon(\x,f_\btheta)=z^*(\x,\y))|\x,\y)]+\lambda \Expect_{\Pemp} [\text{KL}(\Prob(\tilde{\z}^\varepsilon(\x,\y)\rvert \Vxi) \rVert \tilde{\z}^\varepsilon(\x,f_\btheta)|\x,\y)].
\end{align*}
The authors show that their model outperforms ILO trained using \model{SQP} and \model{cvxpylayers} in terms of computational time and gives lower task loss than sequential models trained using MLE and policy learning with neural networks.

\subsubsection{Supervised learning.}
\cite{Grigas2021integrated} solve a CSO problem with a convex and non-negative decision regularizer $\Omega(\z)$ assuming that the uncertain parameter $\y$ has discrete support. Their model, called \model{ICEO-$\lambda$}, is thus trained by solving:
\begin{subequations}
    \begin{align}
        \model{(ICEO-$\lambda$)}\quad\quad \min_{\btheta}\quad&H(z^*_\lambda(\cdot, f_{\btheta}),\Pemp)   +\lambda \Expect_{\Pemp} [\Omega(z_\lambda^*(\x, f_{\btheta}))]\\
        \subto \quad & z_\lambda^*(\x, f_{\btheta}) = \argmin_{\z} c(\z,f_{\btheta}(\x))+\lambda \Omega(\z), \forall \x.
   \end{align}
\end{subequations}
The regularization ensures uniqueness and Lipschitz property of $z^*_\lambda(\x,f_{\btheta})$ with respect to $f_{\btheta}$ and leads to finite-sample guarantees. To circumvent the challenge associated with non-differentiability of $z^*_\lambda(\x, f_\btheta)$ with respect to $\btheta$, they replace $z^*_\lambda(\x, f_{\btheta})$ with a smooth approximation $\tilde{z}_\lambda(\x, f_{\btheta})$ that is learned using a random data set $(\boldsymbol{p}_i,\z_i)$ generated by sampling $\boldsymbol{p}_i$ from the probability simplex over the discrete support and then finding the optimal solution $\z_i$. They show asymptotic optimality and consistency of their solutions when the hypothesis class is well-specified. They compare their approach to other ILO pipelines and to the SLO approach that estimates the conditional distribution using cross-entropy.

\cite{Cristian2022projectnet} introduce the \model{ProjectNet} model to solve uncertain constrained linear programs in an end-to-end framework by training an optimal policy network, which employs a differentiable approximation of the step of projection to feasibility. 

Another approach, related to \cite{Berthet2020}, that generalizes beyond LPs is given in \cite{shah2022decision} that constructs locally optimized decision losses (\model{LODL}) with supervised learning to directly evaluate the performance of the predictors on the downstream optimization task. To learn a convex \model{LODL} for each data point, this approach first generates labels in the neighborhood of label $\y_i$ in the training set, e.g., by adding Gaussian noise, and then chooses the parameter that minimizes the MSE between \model{LODL} and the downstream \revised{task} loss. The \model{LODL} is used in place of the task-specific surrogate optimization layers and outperforms SLO on three resource allocation problems (linear top-1 item selection problem, web advertising, and portfolio optimization). The numerical experiments indicate that handcrafted surrogate functions only perform better for the web advertising problem.

\subsection{Applications}
In this subsection, we discuss the applications of the ILO framework to a wide range of real-world problems.

\cite{tian2023smart} and \cite{Tian_Yan_Wang_Laporte_2023} use \model{SPOT} and noise-contrastive estimation method \citep{mulamba2021contrastive}, respectively, to solve the maritime transportation problem. A comprehensive tutorial on prescriptive analytics methods for logistics is given in \cite{tian:tutorial}.   \model{SPO} has been used in solving last-mile delivery \citep{Chu_Zhang_Bai_Chen_2021} and ship inspection problems \citep{Yan_Wang_Fagerholt_2020,yan2021,yan2023}. \cite{Demirovic_Stuckey_Bailey2019} and \cite{Demirovic_Stuckey_Guns_Bailey_Leckie_Ramamohanarao_Chan_2020} minimize the same expected regret as \model{SPO} for specific applications related to ranking optimization and dynamic programming problems, respectively.

\cite{perrault_end--end_2020} solve a Stackelberg security game with the ILO framework by learning the attack probability distribution over a discrete set of targets to maximize a surrogate for the defender's expected utility. They show that their model results in higher expected utility for the defender on synthetic and human subjects data than the sequential models that learn the attack probability by minimizing the cross entropy loss. \cite{wang_automatically_2020} replace the large-scale optimization problem with a low dimensional surrogate by reparameterizing the feasible space of decisions. They observe significant performance improvements for non-convex problems compared to the strongly convex case.

\cite{Stratigakos_Camal_Michiorri_Kariniotakis_2022} \revised{solve} an integrated forecasting and optimization model for trading in renewable energy that trains an ensemble of prescriptive trees by randomly splitting the feature space $\X$ based on the task-specific cost function. \cite{Sang_Xu_Long_Hu_Sun_2022} introduce an \revised{ILO framework for} electricity price prediction for energy storage system arbitrage. They present a hybrid loss function to measure prediction and decision errors and a hybrid stochastic gradient descent learning method. \cite{sang_2023} solve a voltage regulation problem using a similar hybrid loss function, and backpropagation is done by implicitly differentiating the optimality conditions of a second-order cone program.

\cite{Liu_Yin_Bai_Grimm_2023} use \revised{a} DNN to model the routing behavior of users in a transportation network and learn the parameters by minimizing the mismatch between the flow prescribed by the variational inequality and the observed flow. The backward pass is obtained by applying the IFT to the variational inequality. \cite{Wahdany_Schmitt_Cremer_2023} propose an integrated model for wind-power forecasting that learns the parameters of a neural network to optimize the energy system costs under the system constraints. \citet{Vohra2023end} apply similar ideas to develop end-to-end renewable energy generation forecasts, using multiple contextual sources such as satellite images and meteorological time series.

\cite{Butler_Kwon_portfolio_2023} solves the contextual mean-variance portfolio (MVP) optimization problem by learning the parameters of the linear prediction model using the ILO framework. The covariance matrix is estimated using the exponentially weighted moving average model. They provide analytical solutions to unconstrained and equality-constrained MVP optimization problems and show that they outperform SLO models based on ordinary least squares regression. These analytical solutions lead to lower variance when compared with the exact solutions of the corresponding inequality-constrained MVP optimization problem.

\section{Active research directions}
\label{sec:directions}
We now summarize \revised{active and future} research directions for further work in contextual optimization.

\paragraph{Uncertainty in constraints.} Most studies on contextual optimization assume that there is no uncertainty in the constraints. If constraints are also uncertain, the SAA solutions that ignore the covariates information might not be feasible \citep{Rahimian2022}. \cite{Bertsimas_Kallus_2020} have highlighted the challenges in using ERM in a constrained CSO problem. \cite{Rahimian2022} solve a conditional chance-constrained program that ensures with a high probability that the solution remains feasible \revised{under the conditional distribution given the realized \revised{covariates}}. Although they do not focus on contextual optimization, interesting links can be found with the literature on constraint learning \citep{Fajemisin2023optimization} and inverse optimization \citep{Chan2021inverse}.

\paragraph{Risk aversion.}
There has been a growing interest in studying contextual optimization in the risk-averse setting. Specifically, one can consider replacing the risk-neutral expectation from \eqref{eq:CSO} with a risk measure such as value-at-risk. By doing so, one would expect, with a high probability, that a decision-maker's loss is lower than a particular threshold. One can easily represent such a risk measure using an uncertainty set which represents the set of all possible outcomes that may occur in the future. The resulting uncertainty set should be carefully chosen. It should capture the most relevant scenarios to balance the trade-off between avoiding risks and obtaining returns. The recently proposed Conditional Robust Optimization (CRO) paradigm by \cite{Chenreddy_Bandi_Delage_2022} \citep[see also][]{ohmori2021predictive, Sun_Shi_Wang2023, Pervsak2023contextual} consists in learning a conditional set $\mathcal{U}(\x)$ to solve the following problem:
\begin{align}
    \text{(CRO)} \quad & \min_{\z\in \Z} \max_{\Vxi \in \mathcal{U}(\x)} c(\z,\Vxi),
    \label{CRO}
\end{align}
where $\mathcal{U}(\x)$ is an uncertainty set designed to contain with high probability the realization of $\Vxi$
conditionally on observing $\x$. Their approach solves the CRO problem sequentially where $\mathcal{U}(\x)$ is learned first and is subsequently used to solve the downstream RO problem. A challenging problem is to learn the uncertainty set to minimize the downstream cost function.

\paragraph{Toolboxes and benchmarking.} Several toolboxes and packages have been proposed recently to train decision pipelines. \cite{Agrawal_Amos_Barratt_Boyd_Diamond_Kolter_2019} provide the \model{cvxpylayers} library, which includes a subclass of convex optimization problems as differentiable layers in auto-differentiation libraries in PyTorch, TensorFlow, and JAX. Other libraries for differentiating non-linear optimization problems for end-to-end learning include higher \citep{grefenstette2019generalized}, JAXopt \citep{Blondel_Berthet_Cuturi2022}, TorchOpt \citep{ren2022torchopt}, and Theseus \citep{pineda2022theseus}. \cite{Tang_Khalil_2022} introduce \revised{an open-source software package called PyEPO (Pytorch-based End-to-End Predict-then-Optimize) implemented} in Python for ILO of problems that are linear in uncertain parameters. They implement various existing methods, such as \model{SPO+}, \model{DBB},  and \model{PFYL}. They also include new benchmarks and comprehensive experiments highlighting the advantages of integrated learning. \citet{Dalle2022} provide similar tools for combinatorial problems in Julia. 

Comparisons of existing approaches in fixed simulation settings are scarce, especially with real-world data. \citet{Buttler2022meta} provide a meta-analysis of selected methods on an unconstrained newsvendor problem on four data sets from the retail and food sectors. They highlight that there is no single method that clearly outperforms all the others on the four data sets. \revised{\cite{mandi2023decision} carried out a comprehensive benchmarking of ILO frameworks tailored for expected value-based models on seven distinct problems using public datasets.}

\paragraph{Endogenous uncertainty.} 
While there has been some progress in studying problems where the decision affects the uncertain parameters \citep{Basciftciejor2021, liu2022coupled}, the literature on decision-dependent uncertainty with covariates is sparse \citep{Bertsimas_Kallus_2020, Bertsimas_Koduri_2022}. An example could be a facility location problem where demand changes once a facility is located in a region or a price-setting newsvendor problem whose demand depends on the price \citep{Liu2023solving}. In these problems, the causal relationship between demand and prices is unknown. These examples offer interesting parallels with the research on heterogeneous treatment effects such as \citet{Wager2018estimation}, which introduce causal forests for estimating treatment effects and provide asymptotic consistency results. \citet{Alley2023pricing} study a price-setting problem and provide a new loss function to isolate the causal effects of price on demand from the conditional effects due to other \revised{covariates}.

\paragraph{Data privacy.} Another issue is that the data might come from multiple sources and contain sensitive private information, so it cannot be directly provided in its original form to the system operator. Differential privacy techniques \citep[see, e.g.,][]{Abadi2016} can be used to obfuscate data but may impact predictive and prescriptive performance. \cite{Mieth_Morales_Poor_2023} determine the data quality after obfuscation in an optimal power flow problem with a Wasserstein ambiguity set and use a DRO model to determine the data value for decision-making.

\paragraph{Interpretability \& explainability.} Decision pipelines must be trusted to be implemented. This is evident from the European Union legislation \quoteIt{General Data Protection Regulation} that requires entities using automated systems to provide \quoteIt{meaningful information about the logic involved} in making decisions, known popularly as the \quoteIt{right to explanation} \citep{Doshi-Velez_Kim_2017,kaminski2019right}. For instance, a citizen has the right to ask a bank for an explanation in the case of loan denial. While interpretability has received much attention in predictive ML applications \citep{Rudin2019}, it remains largely unexplored in a contextual optimization, i.e., prescriptive context. Interpretability requires transparent decision pipelines that are intelligible to users, e.g., built over simple models such as decision trees or rule lists. In contrast, explainability may be achieved with an additional algorithm on top of a black box or complex model. Feature importance has been analyzed in a prescriptive context by \cite{serrano_bilevel_2022}. They introduce an integrated approach that solves a bilevel program with an integer master problem optimizing (cross-)validation accuracy. To achieve explainability, \citet{Forel2023explainable} adapt the concept of counterfactual explanations to explain a given data-driven decision through differences of context that make this decision optimal, or better suited than a given expert decision. Having identified these differences, it becomes possible to correct or complete the contextual information, if necessary, or otherwise to give explanative elements supporting different decisions. \revised{Another research direction could be to train tree-based models (such as optimal classification trees) to approximate the policy of a complex learning-and-optimization pipeline. This has interesting connections with model distillation, i.e., the idea in the ML community of approximating a large model by a smaller one, and the work of \citet{Bertsimas_Stellato_2021}, which learns the mapping from the problem parameters to optimal decisions through interpretable models.}

\paragraph{Fairness.} Applying decisions based on contextual information can raise fairness issues when the context is made of protected attributes. This has been studied especially in pricing problems, to \revised{ensure} that different customers or groups of customers are proposed prices that \revised{don't} differ \revised{significantly}\citep{Cohen2021dynamic, Cohen2022price}.

\paragraph{Finite sample guarantees for ILO.} \revised{In \cite{Grigas2021integrated}, the authors derive finite-sample guarantees for ILO under the assumption of discrete support for the uncertain parameter.} An open problem is to derive generalization bounds on the performance of ILO models for non-linear problems \revised{where the uncertain parameters have continuous support.}

\revised{\paragraph{Correcting for in-sample bias of data-driven optimization.}
When devising an optimal policy based on a finite number of samples, it is desired that low in-sample risk translates to low out-of-sample risk. However, decision rule optimization in \eqref{eq:policy_opt} or learning and optimization model in \eqref{eq:optim} are known to produce optimistically biased estimates of the true expected cost of the prescribed policy \citep{ban_big_2019, Costa_Iyengar_2023, gupta2022debiasing}. While one can replace this estimation with an unbiased one if data was reserved for this purpose, this is usually considered a wasteful use of data given that it could instead have been used to obtain a better-performing policy. Recent research has identified ways of circumventing this issue by estimating and correcting for the in-sample bias in contextual \citep{gupta2022debiasing} and non-contextual \citep{Ito_Yabe_Fujimaki_2018, Gupta_Rusmevichientong_2021, Iyengar_Lam_Wang_2023} stochastic optimization problems under the assumption that errors in the estimation of uncertain parameters are normally distributed. In addition to correcting the bias using Stein's lemma, \cite{Gupta_Rusmevichientong_2021} show that certain ``Bayes-inspired'' and regularized policies achieve the same performance as optimal in-sample policy in small-data large-scale regimes. 

A promising future research direction could be to build a general framework to learn the in-sample policies that directly minimize the debiased objective functions. In this regard, one might find inspiration from the work of \cite{Gupta_Rusmevichientong_2021} addressing a similar issue in the non-contextual setting.}%

\paragraph{Multi-agent decision-making.} A multi-agent perspective becomes necessary in transportation and operations management problems, where different agents have access to different sources of information (i.e. covariates). In this regard, some recent work by \cite{Heaton_McKenzie_Li_Fung_Osher_Yin_2022} identifies the Nash equilibrium of contextual games using implicit differentiation of variational inequalities and {jacobian-free backpropagation}.

\paragraph{Costly label acquisition.}
In many applications, it is costly to \revised{gather observations of uncertain vectors and covariate pairs}. For instance, in personalized pricing, surveys can be sent to customers to obtain information on the sensitivity of purchasing an item with respect to its price. However, creating, sending, and collecting the surveys may have a cost. \cite{Liu_Grigas_Liu_Shen_2023} develop an active learning approach to obtain labels to solve the \model{SPO} problem, while the more general case of developing active learning methods for non-linear contextual optimization is an interesting future direction. \citet{Besbes_Ma_Mouchtaki_2023} provide theoretical results on the trade-off between the quality and quantity of data in a newsvendor problem, thus guiding decision-makers on how to invest in data acquisition strategies.

\paragraph{Multi-stage contextual optimization.} Most works on contextual optimization focus on single and two-stage problems. \citet{Ban2019residual_tree} and \citet{Rios2015multi} use the residuals of the regression model to build multi-stage scenario trees and solve multi-stage CSO problems. \citet{Bertsimas2023dynamic} generalize the weighted SAA model for multi-stage problems. \cite{Qi_MS2023} propose an end-to-end learning framework to solve a real-world multistage inventory replenishment problem.

An active area of research is sequential decision-making with uncertainty. \revised{Inverse reinforcement learning \citep{ng2000algorithms} focuses on learning rewards that are consistent with observed trajectories. In the econometrics literature on dynamic discrete choice modeling the focus lies more broadly on estimating structural parameters of} Markov decision processes (MDPs) \revised{\citep{rust1998,AguiMira2010} including rewards, transition functions and discount factors.} \revised{On both topics, e}stimates are \revised{typically} obtained through MLE employing a soft version of the Bellman operator \citep[e.g.,][]{rust1987optimal, ziebart2008maximum}. 

\revised{In the context of} model-based reinforcement learning, \revised{so-called decision awareness \citep[i.e. explicitly training components of a reinforcement learning system to help the agent improve the total amount of collected reward,][]{NikishinEtAl22_workshop} is receiving increasing attention \citep[e.g.,][]{JosephEtAl13, FarahEtAl17, Farah18, GrimmEtAl20}. For example,} \cite{nikishin2022control} introduce an approach that combines learning and planning to optimize expected returns for both tabular and non-tabular MDPs. They employ the soft version of the Bellman operator for efficient parameter learning using the IFT and show that their state-action value function has a lower approximation error than that of MLE in tabular MDPs.

Another interesting research direction is to challenge the assumption that the joint distribution of the \revised{covariate} and uncertain parameters is stationary. \citet{Neghab2022integrated} study a newsvendor model with a hidden Markov model underlying the distribution of the covariates and demand.

Finally, an area that requires attention is the deployment of models for real-world applications by tackling computational hurdles associated with decision-aware learning in MDPs, such as large state-action pairs and high-dimensional policy spaces \citep{wang_scalable_2023}. An example is a service call scheduling problem that is formulated as a restless multi-armed bandit~(RMAB) problem in \cite{Mate_Madaan_Taneja_Madhiwalla_Verma_Singh_Hegde_Varakantham_Tambe_2022} to improve maternal and child health in a non-profit organization. They model each beneficiary as an arm, apply a clustering method to learn the dynamics, and then use the Whittle Index policy to solve the RMAB. \cite{wang_scalable_2023} use decision-focused learning to solve RMAB, where the computational difficulty in differentiating the Whittle index policy of selecting the top-k arms, is mitigated by making a soft-top-k selection of arms which is an optimal transport problem \citep{Xie_Dai_Chen_Dai_Zhao_Zha_Wei_Pfister_2020}.

\revised{
\section{Conclusions}
\label{sec:conclusion}
In this survey, we summarize advancements in contextual optimization methods developed to solve decision-making problems under uncertainty. A salient feature of the problems studied is that a covariate that is correlated with the uncertain parameter is revealed to the decision-maker before a decision is made. Therefore, historical data on both covariates and corresponding uncertain parameter values can be used to prescribe decisions based on the covariate information. We showed that contextual optimization literature can be categorized into three data-driven frameworks for learning policies, namely, decision rule optimization, sequential learning and optimization, and integrated learning and optimization. While decision rule optimization explicitly parameterizes the policy as a function of the covariates, learning and optimization require estimating the conditional distribution (or a sufficient statistic in the case of expected value-based models) of the uncertain parameter given the covariates in a sequential or integrated manner.
By providing a parametric description of these three data-driven frameworks, we have introduced a uniform notation and terminology for analyzing different methods. In particular, we have elaborately described integrated learning and optimization models using both their training pipeline and modeling choice, further emphasizing that these two aspects are intertwined. Furthermore, we have provided a list of active fields of research in CSO problems.
}

\section*{Acknowledgments}
The first author’s research is supported by the Group for Research in Decision Analysis (GERAD) postdoctoral fellowship and Fonds de Recherche du Québec–Nature et Technologies (FRQNT) postdoctoral research scholarship [Grant 301065]. Erick Delage was partially supported by the Canadian Natural Sciences and Engineering Research Council [Grant RGPIN-2022-05261] and by the Canada Research Chair program [950-230057]. Emma Frejinger was partially supported by the Canada Research Chair program [950-232244]. Finally, the authors also acknowledge the support of IVADO and the Canada First Research Excellence Fund (Apogée/CFREF).

\bibliographystyle{informs2014}
\bibliography{biblio}

\begin{thebibliography}{194}
\providecommand{\natexlab}[1]{#1}
\providecommand{\url}[1]{\texttt{#1}}
\providecommand{\urlprefix}{URL }

\bibitem[{Abadi et~al.(2016)Abadi, Chu, Goodfellow, McMahan, Mironov, Talwar, \protect\BIBand{} Zhang}]{Abadi2016}
Abadi M, Chu A, Goodfellow I, McMahan HB, Mironov I, Talwar K, Zhang L (2016) Deep learning with differential privacy. \emph{Proceedings of the 2016 ACM SIGSAC Conference on Computer and Communications Security}, 308–318, CCS '16 (NY, USA: Association for Computing Machinery).

\bibitem[{Agrawal et~al.(2019)Agrawal, Amos, Barratt, Boyd, Diamond, \protect\BIBand{} Kolter}]{Agrawal_Amos_Barratt_Boyd_Diamond_Kolter_2019}
Agrawal A, Amos B, Barratt S, Boyd S, Diamond S, Kolter JZ (2019) Differentiable convex optimization layers. \emph{Advances in Neural Information Processing Systems}, volume~32 (Curran Associates, Inc.).

\bibitem[{Aguirregabiria \protect\BIBand{} Mira(2010)}]{AguiMira2010}
Aguirregabiria V, Mira P (2010) Dynamic discrete choice structural models: A survey. \emph{Journal of Econometrics} 156(1):38--67.

\bibitem[{Alley et~al.(2023)Alley, Biggs, Hariss, Herrmann, Li, \protect\BIBand{} Perakis}]{Alley2023pricing}
Alley M, Biggs M, Hariss R, Herrmann C, Li ML, Perakis G (2023) Pricing for heterogeneous products: Analytics for ticket reselling. \emph{Manufacturing \& Service Operations Management} 25(2):409--426.

\bibitem[{Amos \protect\BIBand{} Kolter(2017)}]{amos_optnet_2017}
Amos B, Kolter JZ (2017) {OptNet}: Differentiable optimization as a layer in neural networks. \emph{International Conference on Machine Learning}, volume~70, 136--145 (PMLR).

\bibitem[{Aronszajn(1950)}]{aronszajn1950theory}
Aronszajn N (1950) Theory of reproducing kernels. \emph{Transactions of the American mathematical society} 68(3):337--404.

\bibitem[{Athey et~al.(2019)Athey, Tibshirani, \protect\BIBand{} Wager}]{Athey_Tibshirani_Wager_2019}
Athey S, Tibshirani J, Wager S (2019) Generalized random forests. \emph{The Annals of Statistics} 47(2):1148 -- 1178.

\bibitem[{Backhoff et~al.(2017)Backhoff, Beiglbock, Lin, \protect\BIBand{} Zalashko}]{backhoff2017causal}
Backhoff J, Beiglbock M, Lin Y, Zalashko A (2017) Causal transport in discrete time and applications. \emph{SIAM Journal on Optimization} 27(4):2528--2562.

\bibitem[{Bai et~al.(2019)Bai, Kolter, \protect\BIBand{} Koltun}]{Bai_Kolter_Koltun_2019}
Bai S, Kolter JZ, Koltun V (2019) Deep equilibrium models. \emph{Advances in Neural Information Processing Systems}, volume~32 (Curran Associates, Inc.).

\bibitem[{Ban et~al.(2019)Ban, Gallien, \protect\BIBand{} Mersereau}]{Ban2019residual_tree}
Ban GY, Gallien J, Mersereau AJ (2019) Dynamic procurement of new products with covariate information: The residual tree method. \emph{Manufacturing \& Service Operations Management} 21(4):798--815.

\bibitem[{Ban \protect\BIBand{} Rudin(2019)}]{ban_big_2019}
Ban GY, Rudin C (2019) The {Big} {Data} {Newsvendor}: {Practical} {Insights} from {Machine} {Learning}. \emph{Operations Research} 67(1):90--108.

\bibitem[{Basciftci et~al.(2021)Basciftci, Ahmed, \protect\BIBand{} Shen}]{Basciftciejor2021}
Basciftci B, Ahmed S, Shen S (2021) Distributionally robust facility location problem under decision-dependent stochastic demand. \emph{European Journal of Operational Research} 292(2):548--561.

\bibitem[{Bazier-Matte \protect\BIBand{} Delage(2020)}]{Bazier-Matte_Delage_2020}
Bazier-Matte T, Delage E (2020) Generalization bounds for regularized portfolio selection with market side information. \emph{INFOR: Information Systems and Operational Research} 58(2):374–401.

\bibitem[{Bengio(1997)}]{Bengio_1997}
Bengio Y (1997) Using a financial training criterion rather than a prediction criterion. \emph{International Journal of Neural Systems} 8(4):433–443.

\bibitem[{Bengio et~al.(2021)Bengio, Lodi, \protect\BIBand{} Prouvost}]{Bengio_2021}
Bengio Y, Lodi A, Prouvost A (2021) Machine learning for combinatorial optimization: A methodological tour d’horizon. \emph{European Journal of Operational Research} 290(2):405--421.

\bibitem[{Berthet et~al.(2020)Berthet, Blondel, Teboul, Cuturi, Vert, \protect\BIBand{} Bach}]{Berthet2020}
Berthet Q, Blondel M, Teboul O, Cuturi M, Vert JP, Bach F (2020) Learning with differentiable perturbed optimizers. \emph{Advances in Neural Information Processing Systems}, volume~33, 9508--9519 (Curran Associates, Inc.).

\bibitem[{Bertsimas et~al.(2019)Bertsimas, Dunn, \protect\BIBand{} Mundru}]{doi:10.1287/ijoo.2018.0005}
Bertsimas D, Dunn J, Mundru N (2019) Optimal prescriptive trees. \emph{INFORMS Journal on Optimization} 1(2):164--183.

\bibitem[{Bertsimas \protect\BIBand{} Kallus(2020)}]{Bertsimas_Kallus_2020}
Bertsimas D, Kallus N (2020) From predictive to prescriptive analytics. \emph{Management Science} 66(3):1025–1044.

\bibitem[{Bertsimas \protect\BIBand{} Koduri(2022)}]{Bertsimas_Koduri_2022}
Bertsimas D, Koduri N (2022) Data-driven optimization: A reproducing kernel hilbert space approach. \emph{Operations Research} 70(1):454–471.

\bibitem[{Bertsimas et~al.(2023)Bertsimas, McCord, \protect\BIBand{} Sturt}]{Bertsimas2023dynamic}
Bertsimas D, McCord C, Sturt B (2023) Dynamic optimization with side information. \emph{European Journal of Operational Research} 304(2):634--651.

\bibitem[{Bertsimas \protect\BIBand{} Stellato(2021)}]{Bertsimas_Stellato_2021}
Bertsimas D, Stellato B (2021) The voice of optimization. \emph{Machine Learning} 110(2):249–277.

\bibitem[{Bertsimas \protect\BIBand{} Van~Parys(2022)}]{Bertsimas2022bootstrap}
Bertsimas D, Van~Parys B (2022) Bootstrap robust prescriptive analytics. \emph{Mathematical Programming} 195(1-2):39--78.

\bibitem[{Besbes et~al.(2023)Besbes, Ma, \protect\BIBand{} Mouchtaki}]{Besbes_Ma_Mouchtaki_2023}
Besbes O, Ma W, Mouchtaki O (2023) Quality vs. quantity of data in contextual decision-making: Exact analysis under newsvendor loss. \emph{arXiv preprint arXiv:2302.08424} .

\bibitem[{Birge \protect\BIBand{} Louveaux(2011)}]{birge2011introduction}
Birge JR, Louveaux F (2011) \emph{Introduction to Stochastic Programming} (New York, NY: Springer New York).

\bibitem[{Blondel et~al.(2022)Blondel, Berthet, Cuturi, Frostig, Hoyer, Llinares-L{\'o}pez, Pedregosa, \protect\BIBand{} Vert}]{Blondel_Berthet_Cuturi2022}
Blondel M, Berthet Q, Cuturi M, Frostig R, Hoyer S, Llinares-L{\'o}pez F, Pedregosa F, Vert JP (2022) Efficient and modular implicit differentiation. \emph{Advances in Neural Information Processing Systems}, volume~35, 5230--5242 (Curran Associates, Inc.).

\bibitem[{Blondel et~al.(2020)Blondel, Martins, \protect\BIBand{} Niculae}]{Blondel2020learning}
Blondel M, Martins AF, Niculae V (2020) Learning with {Fenchel-Young} losses. \emph{Journal of Machine Learning Research} 21(1):1314--1382.

\bibitem[{Bolte et~al.(2021)Bolte, Le, Pauwels, \protect\BIBand{} Silveti-Falls}]{bolte2021nonsmooth}
Bolte J, Le T, Pauwels E, Silveti-Falls T (2021) Nonsmooth implicit differentiation for machine-learning and optimization. \emph{Advances in Neural Information Processing Systems} 34:13537--13549.

\bibitem[{Butler \protect\BIBand{} Kwon(2023{\natexlab{a}})}]{Butler_Kwon_ADMM2023}
Butler A, Kwon RH (2023{\natexlab{a}}) Efficient differentiable quadratic programming layers: an {ADMM} approach. \emph{Computational Optimization and Applications} 84(2):449–476.

\bibitem[{Butler \protect\BIBand{} Kwon(2023{\natexlab{b}})}]{Butler_Kwon_portfolio_2023}
Butler A, Kwon RH (2023{\natexlab{b}}) Integrating prediction in mean-variance portfolio optimization. \emph{Quantitative Finance} 23(3):429--452.

\bibitem[{Buttler et~al.(2022)Buttler, Philippi, Stein, \protect\BIBand{} Pibernik}]{Buttler2022meta}
Buttler S, Philippi A, Stein N, Pibernik R (2022) A meta analysis of data-driven newsvendor approaches. \emph{ICLR 2022 Workshop on Setting up ML Evaluation Standards to Accelerate Progress}.

\bibitem[{Chan et~al.(2021)Chan, Mahmood, \protect\BIBand{} Zhu}]{Chan2021inverse}
Chan TC, Mahmood R, Zhu IY (2021) Inverse optimization: Theory and applications. \emph{arXiv preprint arXiv:2109.03920} .

\bibitem[{Chen et~al.(2021)Chen, Donti, Baker, Kolter, \protect\BIBand{} Bergés}]{Chen_Donti_Baker_Kolter_Berges_2021}
Chen B, Donti PL, Baker K, Kolter JZ, Bergés M (2021) Enforcing policy feasibility constraints through differentiable projection for energy optimization. \emph{Proceedings of the Twelfth ACM International Conference on Future Energy Systems}, 199–210 (Virtual Event Italy: ACM).

\bibitem[{Chen \protect\BIBand{} Paschalidis(2019)}]{Chen_Paschalidis_2019}
Chen R, Paschalidis I (2019) Selecting optimal decisions via distributionally robust nearest-neighbor regression. \emph{Advances in Neural Information Processing Systems}, volume~32 (Curran Associates, Inc.).

\bibitem[{Chen \protect\BIBand{} Paschalidis(2018)}]{chen2018robust}
Chen R, Paschalidis IC (2018) A robust learning approach for regression models based on distributionally robust optimization. \emph{Journal of Machine Learning Research} 19(13):1--48.

\bibitem[{Chen et~al.(2023)Chen, Tanneau, \protect\BIBand{} Van~Hentenryck}]{Chen_Tanneau_VanHentenryck2023}
Chen W, Tanneau M, Van~Hentenryck P (2023) End-to-end feasible optimization proxies for large-scale economic dispatch. \emph{arXiv preprint arXiv:2304.11726} .

\bibitem[{Chenreddy et~al.(2022)Chenreddy, Bandi, \protect\BIBand{} Delage}]{Chenreddy_Bandi_Delage_2022}
Chenreddy AR, Bandi N, Delage E (2022) Data-driven conditional robust optimization. \emph{Advances in Neural Information Processing Systems}, volume~35, 9525–9537 (Curran Associates, Inc.).

\bibitem[{Chu et~al.(2023)Chu, Zhang, Bai, \protect\BIBand{} Chen}]{Chu_Zhang_Bai_Chen_2021}
Chu H, Zhang W, Bai P, Chen Y (2023) Data-driven optimization for last-mile delivery. \emph{Complex \& Intelligent Systems} 9(9):2271–2284.

\bibitem[{Chung et~al.(2022)Chung, Rostami, Bastani, \protect\BIBand{} Bastani}]{chung2022decision}
Chung TH, Rostami V, Bastani H, Bastani O (2022) Decision-aware learning for optimizing health supply chains. \emph{arXiv preprint arXiv:2211.08507} .

\bibitem[{Ciocan \protect\BIBand{} Mi\v{s}i\'{c}(2022)}]{doi:10.1287/mnsc.2020.3592}
Ciocan DF, Mi\v{s}i\'{c} VV (2022) Interpretable optimal stopping. \emph{Management Science} 68(3):1616--1638.

\bibitem[{Clarke(1990)}]{clarke1990optimization}
Clarke FH (1990) \emph{Optimization and nonsmooth analysis} (SIAM).

\bibitem[{Cohen et~al.(2022)Cohen, Elmachtoub, \protect\BIBand{} Lei}]{Cohen2022price}
Cohen MC, Elmachtoub AN, Lei X (2022) Price discrimination with fairness constraints. \emph{Management Science} 68(12):8536--8552.

\bibitem[{Cohen et~al.(2021)Cohen, Miao, \protect\BIBand{} Wang}]{Cohen2021dynamic}
Cohen MC, Miao S, Wang Y (2021) Dynamic pricing with fairness constraints. \emph{Available at \url{https://dx.doi.org/10.2139/ssrn.3930622}} .

\bibitem[{Costa \protect\BIBand{} Iyengar(2023)}]{Costa_Iyengar_2023}
Costa G, Iyengar GN (2023) Distributionally robust end-to-end portfolio construction. \emph{Quantitative Finance} 23(10):1465–1482.

\bibitem[{Cristian et~al.(2022)Cristian, Harsha, Perakis, Quanz, \protect\BIBand{} Spantidakis}]{Cristian2022projectnet}
Cristian R, Harsha P, Perakis G, Quanz BL, Spantidakis I (2022) End-to-end learning via constraint-enforcing approximators for linear programs with applications to supply chains. \emph{AI for Decision Optimization Workshop of the AAAI Conference on Artificial Intelligence} .

\bibitem[{Cybenko(1989)}]{cybenko1989approximation}
Cybenko G (1989) Approximation by superpositions of a sigmoidal function. \emph{Mathematics of Control, Signals and Systems} 2(4):303--314.

\bibitem[{Dalle et~al.(2022)Dalle, Baty, Bouvier, \protect\BIBand{} Parmentier}]{Dalle2022}
Dalle G, Baty L, Bouvier L, Parmentier A (2022) Learning with combinatorial optimization layers: a probabilistic approach. \emph{arXiv preprint arXiv:2207.13513} .

\bibitem[{Davis \protect\BIBand{} Yin(2017)}]{davis2017threeoperator}
Davis D, Yin W (2017) A three-operator splitting scheme and its optimization applications. \emph{Set-Valued and Variational Analysis} 25(4):829--858.

\bibitem[{Demirovi{\'c} et~al.(2019)Demirovi{\'c}, J~Stuckey, Bailey, Chan, Leckie, Ramamohanarao, Guns, \protect\BIBand{} Kraus}]{Demirovic_Stuckey_Bailey2019}
Demirovi{\'c} E, J~Stuckey P, Bailey J, Chan J, Leckie C, Ramamohanarao K, Guns T, Kraus S (2019) Predict+ optimise with ranking objectives: Exhaustively learning linear functions. \emph{Proceedings of the Thirtieth International Joint Conference on Artificial Intelligence}, 1078--1085 (International Joint Conferences on Artificial Intelligence Organization).

\bibitem[{Demirović et~al.(2020)Demirović, Stuckey, Guns, Bailey, Leckie, Ramamohanarao, \protect\BIBand{} Chan}]{Demirovic_Stuckey_Guns_Bailey_Leckie_Ramamohanarao_Chan_2020}
Demirović E, Stuckey PJ, Guns T, Bailey J, Leckie C, Ramamohanarao K, Chan J (2020) Dynamic programming for predict+optimise. \emph{Proceedings of the AAAI Conference on Artificial Intelligence} 34(0202):1444–1451.

\bibitem[{Deng \protect\BIBand{} Sen(2022)}]{Deng_Sen_2022}
Deng Y, Sen S (2022) Predictive stochastic programming. \emph{Computational Management Science} 19(1):65–98.

\bibitem[{Domke(2012)}]{domke2012generic}
Domke J (2012) Generic methods for optimization-based modeling. \emph{International Conference on Artificial Intelligence and Statistics}, 318--326 (PMLR).

\bibitem[{Dontchev et~al.(2009)Dontchev, Rockafellar, \protect\BIBand{} Rockafellar}]{dontchev2009implicit}
Dontchev AL, Rockafellar RT, Rockafellar RT (2009) \emph{Implicit functions and solution mappings: A view from variational analysis}, volume 616 (Springer).

\bibitem[{Donti et~al.(2017)Donti, Amos, \protect\BIBand{} Kolter}]{donti_task-based_2019}
Donti P, Amos B, Kolter JZ (2017) Task-based end-to-end model learning in stochastic optimization. \emph{Advances in Neural Information Processing Systems}, volume~30 (Curran Associates, Inc.).

\bibitem[{Donti et~al.(2021)Donti, Roderick, Fazlyab, \protect\BIBand{} Kolter}]{Donti_Roderick_Fazlyab_Kolter_2021}
Donti PL, Roderick M, Fazlyab M, Kolter JZ (2021) Enforcing robust control guarantees within neural network policies. \emph{International Conference on Learning Representations}.

\bibitem[{Doshi-Velez \protect\BIBand{} Kim(2017)}]{Doshi-Velez_Kim_2017}
Doshi-Velez F, Kim B (2017) Towards a rigorous science of interpretable machine learning. \emph{arXiv preprint arXiv:1702.08608} .

\bibitem[{Duvenaud et~al.(2020)Duvenaud, Kolter, \protect\BIBand{} Johnson}]{duvenaud2020deep}
Duvenaud D, Kolter JZ, Johnson M (2020) Deep implicit layers tutorial - neural {ODEs}, deep equilibrium models, and beyond. \emph{Neural Information Processing Systems Tutorial} .

\bibitem[{El~Balghiti et~al.(2019)El~Balghiti, Elmachtoub, Grigas, \protect\BIBand{} Tewari}]{El_Balghiti_Elmachtoub_Grigas_Tewari_2019}
El~Balghiti O, Elmachtoub AN, Grigas P, Tewari A (2019) Generalization bounds in the predict-then-optimize framework. \emph{Advances in Neural Information Processing Systems}, volume~32 (Curran Associates, Inc.).

\bibitem[{El~Balghiti et~al.(2022)El~Balghiti, Elmachtoub, Grigas, \protect\BIBand{} Tewari}]{El_Balghiti_Elmachtoub_Grigas_Tewari_2022}
El~Balghiti O, Elmachtoub AN, Grigas P, Tewari A (2022) Generalization bounds in the predict-then-optimize framework. \emph{Mathematics of Operations Research} Forthcoming.

\bibitem[{Elmachtoub \protect\BIBand{} Grigas(2022)}]{elmachtoub_smart_2022}
Elmachtoub AN, Grigas P (2022) Smart “{Predict}, then {Optimize}”. \emph{Management Science} 68(1):9--26.

\bibitem[{Elmachtoub et~al.(2023)Elmachtoub, Lam, Zhang, \protect\BIBand{} Zhao}]{Elmachtoub_Lam_Zhang_Zhao_2023}
Elmachtoub AN, Lam H, Zhang H, Zhao Y (2023) Estimate-then-optimize versus integrated-estimation-optimization: A stochastic dominance perspective. \emph{arXiv preprint arXiv:2304.06833} .

\bibitem[{Elmachtoub et~al.(2020)Elmachtoub, Liang, \protect\BIBand{} McNellis}]{elmachtoub_decision_2020}
Elmachtoub AN, Liang JCN, McNellis R (2020) Decision trees for decision-making under the predict-then-optimize framework. \emph{International Conference on Machine Learning}, 2858--2867 (PMLR).

\bibitem[{Esteban-Pérez \protect\BIBand{} Morales(2022)}]{Esteban_Perez_Morales_2022}
Esteban-Pérez A, Morales JM (2022) Distributionally robust stochastic programs with side information based on trimmings. \emph{Mathematical Programming} 195(1):1069–1105.

\bibitem[{Esteban-Pérez \protect\BIBand{} Morales(2023)}]{ESTEBANPEREZ20231047}
Esteban-Pérez A, Morales JM (2023) Distributionally robust optimal power flow with contextual information. \emph{European Journal of Operational Research} 306(3):1047--1058.

\bibitem[{Estes \protect\BIBand{} Richard(2023)}]{Estes2023smart}
Estes AS, Richard JPP (2023) Smart predict-then-optimize for two-stage linear programs with side information. \emph{INFORMS Journal on Optimization} Forthcoming.

\bibitem[{Fajemisin et~al.(2023)Fajemisin, Maragno, \protect\BIBand{} den Hertog}]{Fajemisin2023optimization}
Fajemisin AO, Maragno D, den Hertog D (2023) Optimization with constraint learning: a framework and survey. \emph{European Journal of Operational Research} .

\bibitem[{Farahmand(2018)}]{Farah18}
Farahmand AM (2018) Iterative value-aware model learning. Bengio S, Wallach H, Larochelle H, Grauman K, Cesa-Bianchi N, Garnett R, eds., \emph{Advances in Neural Information Processing Systems}, volume~31.

\bibitem[{Farahmand et~al.(2017)Farahmand, Barreto, \protect\BIBand{} Nikovski}]{FarahEtAl17}
Farahmand AM, Barreto A, Nikovski D (2017) {Value-Aware Loss Function for Model-based Reinforcement Learning}. \emph{Proceedings of the 20th International Conference on Artificial Intelligence and Statistics}, volume~54, 1486--1494.

\bibitem[{Ferber et~al.(2020)Ferber, Wilder, Dilkina, \protect\BIBand{} Tambe}]{ferber_mipaal_2020}
Ferber A, Wilder B, Dilkina B, Tambe M (2020) {MIPaaL}: {Mixed} {Integer} {Program} as a {Layer}. \emph{Proceedings of the AAAI Conference on Artificial Intelligence} 34(02):1504--1511.

\bibitem[{Ferreira et~al.(2016)Ferreira, Lee, \protect\BIBand{} Simchi-Levi}]{ferreira2016analytics}
Ferreira KJ, Lee BH, Simchi-Levi D (2016) Analytics for an online retailer: Demand forecasting and price optimization. \emph{Manufacturing \& Service Operations Management} 18(1):69--88.

\bibitem[{Forel et~al.(2023)Forel, Parmentier, \protect\BIBand{} Vidal}]{Forel2023explainable}
Forel A, Parmentier A, Vidal T (2023) Explainable data-driven optimization: From context to decision and back again. \emph{arXiv preprint arXiv:2301.10074} .

\bibitem[{Fung et~al.(2022)Fung, Heaton, Li, Mckenzie, Osher, \protect\BIBand{} Yin}]{Fung_Heaton_Li_Mckenzie_Osher_Yin_2022}
Fung SW, Heaton H, Li Q, Mckenzie D, Osher S, Yin W (2022) {JFB}: Jacobian-free backpropagation for implicit networks. \emph{Proceedings of the AAAI Conference on Artificial Intelligence} 36(66):6648–6656.

\bibitem[{Grefenstette et~al.(2019)Grefenstette, Amos, Yarats, Htut, Molchanov, Meier, Kiela, Cho, \protect\BIBand{} Chintala}]{grefenstette2019generalized}
Grefenstette E, Amos B, Yarats D, Htut PM, Molchanov A, Meier F, Kiela D, Cho K, Chintala S (2019) Generalized inner loop meta-learning. \emph{arXiv preprint arXiv:1910.01727} .

\bibitem[{Grigas et~al.(2021)Grigas, Qi, \protect\BIBand{} Shen}]{Grigas2021integrated}
Grigas P, Qi M, Shen M (2021) Integrated conditional estimation-optimization. \emph{arXiv preprint arXiv:2110.12351} .

\bibitem[{Grimm et~al.(2020)Grimm, Barreto, Singh, \protect\BIBand{} Silver}]{GrimmEtAl20}
Grimm C, Barreto A, Singh S, Silver D (2020) The value equivalence principle for model-based reinforcement learning. \emph{Advances in Neural Information Processing Systems}, volume~33.

\bibitem[{Gupta et~al.(2022)Gupta, Huang, \protect\BIBand{} Rusmevichientong}]{gupta2022debiasing}
Gupta V, Huang M, Rusmevichientong P (2022) Debiasing in-sample policy performance for small-data, large-scale optimization. \emph{Operations Research} .

\bibitem[{Gupta \protect\BIBand{} Rusmevichientong(2021)}]{Gupta_Rusmevichientong_2021}
Gupta V, Rusmevichientong P (2021) Small-data, large-scale linear optimization with uncertain objectives. \emph{Management Science} 67(1):220–241.

\bibitem[{Halkin(1974)}]{Halkin1974}
Halkin H (1974) Implicit functions and optimization problems without continuous differentiability of the data. \emph{SIAM Journal on Control} 12(2):229--236.

\bibitem[{Hannah et~al.(2010)Hannah, Powell, \protect\BIBand{} Blei}]{Hannah2010}
Hannah L, Powell W, Blei D (2010) Nonparametric density estimation for stochastic optimization with an observable state variable. \emph{Advances in Neural Information Processing Systems}, volume~23 (Curran Associates, Inc.).

\bibitem[{Hastie et~al.(2009)Hastie, Tibshirani, \protect\BIBand{} Friedman}]{Hastie2009}
Hastie T, Tibshirani R, Friedman JH (2009) \emph{The Elements of Statistical Learning: Data Mining, Inference, and Prediction}, volume~2 (Springer).

\bibitem[{Heaton et~al.(2022)Heaton, McKenzie, Li, Fung, Osher, \protect\BIBand{} Yin}]{Heaton_McKenzie_Li_Fung_Osher_Yin_2022}
Heaton H, McKenzie D, Li Q, Fung SW, Osher S, Yin W (2022) Learn to predict equilibria via fixed point networks. \emph{arXiv:2106.00906} .

\bibitem[{Ho-Nguyen \protect\BIBand{} K{\i}l{\i}n{\c{c}}-Karzan(2022)}]{ho-nguyen_risk_2022}
Ho-Nguyen N, K{\i}l{\i}n{\c{c}}-Karzan F (2022) Risk guarantees for end-to-end prediction and optimization processes. \emph{Management Science} 68(12):8680--8698.

\bibitem[{Hofmann et~al.(2008)Hofmann, Schölkopf, \protect\BIBand{} Smola}]{hofmann2008kernel}
Hofmann T, Schölkopf B, Smola A (2008) Kernel methods in machine learning. \emph{Annals of Statistics} 36(3):1171--1220.

\bibitem[{Hu et~al.(2022)Hu, Kallus, \protect\BIBand{} Mao}]{Hu_Kallus_Mao_2022}
Hu Y, Kallus N, Mao X (2022) Fast rates for contextual linear optimization. \emph{Management Science} 68(6):4236–4245.

\bibitem[{Huber et~al.(2019)Huber, M{\"u}ller, Fleischmann, \protect\BIBand{} Stuckenschmidt}]{Huber2019}
Huber J, M{\"u}ller S, Fleischmann M, Stuckenschmidt H (2019) A data-driven newsvendor problem: From data to decision. \emph{European Journal of Operational Research} 278(3):904--915.

\bibitem[{Ifrim et~al.(2012)Ifrim, O'Sullivan, \protect\BIBand{} Simonis}]{ifrim2012properties}
Ifrim G, O'Sullivan B, Simonis H (2012) Properties of energy-price forecasts for scheduling. \emph{International Conference on Principles and Practice of Constraint Programming} (Springer).

\bibitem[{Ito et~al.(2018)Ito, Yabe, \protect\BIBand{} Fujimaki}]{Ito_Yabe_Fujimaki_2018}
Ito S, Yabe A, Fujimaki R (2018) Unbiased objective estimation in predictive optimization. \emph{Proceedings of the 35th International Conference on Machine Learning}, 2176–2185 (PMLR).

\bibitem[{Iyengar et~al.(2023)Iyengar, Lam, \protect\BIBand{} Wang}]{Iyengar_Lam_Wang_2023}
Iyengar G, Lam H, Wang T (2023) Optimizer’s information criterion: Dissecting and correcting bias in data-driven optimization. \emph{arXiv preprint arXiv:2306.10081} .

\bibitem[{Jeong et~al.(2022)Jeong, Jaggi, Butler, \protect\BIBand{} Sanner}]{Jeong2022exact}
Jeong J, Jaggi P, Butler A, Sanner S (2022) An exact symbolic reduction of linear smart {P}redict+{O}ptimize to mixed integer linear programming. \emph{International Conference on Machine Learning}, volume 162, 10053--10067 (PMLR).

\bibitem[{Joseph et~al.(2013)Joseph, Geramifard, Roberts, How, \protect\BIBand{} Roy}]{JosephEtAl13}
Joseph J, Geramifard A, Roberts JW, How JP, Roy N (2013) Reinforcement learning with misspecified model classes. \emph{IEEE International Conference on Robotics and Automation}, 939--946.

\bibitem[{Kallus \protect\BIBand{} Mao(2022)}]{kallus_stochastic_2022}
Kallus N, Mao X (2022) Stochastic optimization forests. \emph{Management Science} 69(4):1975--1994.

\bibitem[{Kaminski(2019)}]{kaminski2019right}
Kaminski ME (2019) The right to explanation, explained. \emph{Berkeley Technology Law Journal} 34(1):189--218.

\bibitem[{Kannan et~al.(2021)Kannan, Bayraksan, \protect\BIBand{} Luedtke}]{Kannan_Bayraksan_Luedtke_2021}
Kannan R, Bayraksan G, Luedtke J (2021) Heteroscedasticity-aware residuals-based contextual stochastic optimization. \emph{arXiv preprint arXiv:2101.03139} .

\bibitem[{Kannan et~al.(2020)Kannan, Bayraksan, \protect\BIBand{} Luedtke}]{Kannan2020residuals}
Kannan R, Bayraksan G, Luedtke JR (2020) Residuals-based distributionally robust optimization with covariate information. \emph{arXiv preprint arXiv:2012.01088} .

\bibitem[{Kannan et~al.(2022)Kannan, Bayraksan, \protect\BIBand{} Luedtke}]{Kannan2022}
Kannan R, Bayraksan G, Luedtke JR (2022) Data-driven sample average approximation with covariate information. \emph{arXiv preprint arXiv:2207.13554} .

\bibitem[{Kantorovich \protect\BIBand{} Rubinshtein(1958)}]{kantorovich1958space}
Kantorovich LV, Rubinshtein GS (1958) On a space of totally additive functions. \emph{Vestnik Leningradskogo Universiteta} 13(7):52--59.

\bibitem[{Keshavarz(2022)}]{parisa}
Keshavarz P (2022) \emph{Interpretable Contextual Newsvendor Models: A Tree-Based Method to Solving Data-Driven Newsvendor Problems}. Master's thesis, University of Ottawa.

\bibitem[{Kong et~al.(2022)Kong, Cui, Zhuang, Feng, Prakash, \protect\BIBand{} Zhang}]{Kong2022energy}
Kong L, Cui J, Zhuang Y, Feng R, Prakash BA, Zhang C (2022) End-to-end stochastic optimization with energy-based model. \emph{Advances in Neural Information Processing Systems}, volume~35, 11341--11354 (Curran Associates, Inc.).

\bibitem[{Kotary et~al.(2023)Kotary, Dinh, \protect\BIBand{} Fioretto}]{kotary2023folded}
Kotary J, Dinh MH, Fioretto F (2023) Folded optimization for end-to-end model-based learning. \emph{International Joint Conference on Artificial Intelligence}.

\bibitem[{Kotary et~al.(2021)Kotary, Fioretto, Van~Hentenryck, \protect\BIBand{} Wilder}]{Kotary_Fioretto_Van_Hentenryck_Wilder_2021}
Kotary J, Fioretto F, Van~Hentenryck P, Wilder B (2021) End-to-end constrained optimization learning: A survey. \emph{Proceedings of the Thirtieth International Joint Conference on Artificial Intelligence}, 4475–4482 (International Joint Conferences on Artificial Intelligence Organization).

\bibitem[{Kullback \protect\BIBand{} Leibler(1951)}]{kullback1951information}
Kullback S, Leibler RA (1951) On information and sufficiency. \emph{The Annals of Mathematical Statistics} 22(1):79--86.

\bibitem[{Lassalle(2018)}]{lassalle2018causal}
Lassalle R (2018) Causal transport plans and their {M}onge--{K}antorovich problems. \emph{Stochastic Analysis and Applications} 36(3):452--484.

\bibitem[{Lawless \protect\BIBand{} Zhou(2022)}]{lawless_note_2022}
Lawless C, Zhou A (2022) A note on task-aware loss via reweighing prediction loss by decision-regret. \emph{arXiv preprint arXiv:2211.05116} .

\bibitem[{Lin et~al.(2022)Lin, Chen, Li, \protect\BIBand{} Shen}]{Lin2022}
Lin S, Chen Y, Li Y, Shen ZJM (2022) Data-driven newsvendor problems regularized by a profit risk constraint. \emph{Production and Operations Management} 31(4):1630--1644.

\bibitem[{Liu \protect\BIBand{} Grigas(2021)}]{liu2021risk}
Liu H, Grigas P (2021) Risk bounds and calibration for a smart predict-then-optimize method. \emph{Advances in Neural Information Processing Systems}, volume~34, 22083--22094 (Curran Associates, Inc.).

\bibitem[{Liu et~al.(2022)Liu, Li, \protect\BIBand{} Sen}]{liu2022coupled}
Liu J, Li G, Sen S (2022) Coupled learning enabled stochastic programming with endogenous uncertainty. \emph{Mathematics of Operations Research} 47(2):1681--1705.

\bibitem[{Liu et~al.(2023{\natexlab{a}})Liu, Grigas, Liu, \protect\BIBand{} Shen}]{Liu_Grigas_Liu_Shen_2023}
Liu M, Grigas P, Liu H, Shen ZJM (2023{\natexlab{a}}) Active learning in the predict-then-optimize framework: A margin-based approach. \emph{arXiv preprint arXiv:2305.06584} .

\bibitem[{Liu et~al.(2021)Liu, He, \protect\BIBand{} Max~Shen}]{Liu_He_Max_Shen_2021}
Liu S, He L, Max~Shen ZJ (2021) On-time last-mile delivery: Order assignment with travel-time predictors. \emph{Management Science} 67(7):4095–4119.

\bibitem[{Liu \protect\BIBand{} Zhang(2023)}]{Liu2023solving}
Liu W, Zhang Z (2023) Solving data-driven newsvendor pricing problems with decision-dependent effect. \emph{arXiv preprint arXiv:2304.13924} .

\bibitem[{Liu et~al.(2023{\natexlab{b}})Liu, Yin, Bai, \protect\BIBand{} Grimm}]{Liu_Yin_Bai_Grimm_2023}
Liu Z, Yin Y, Bai F, Grimm DK (2023{\natexlab{b}}) End-to-end learning of user equilibrium with implicit neural networks. \emph{Transportation Research Part C: Emerging Technologies} 150:104085.

\bibitem[{Liyanage \protect\BIBand{} Shanthikumar(2005)}]{Liyanage_Shanthikumar_2005}
Liyanage LH, Shanthikumar JG (2005) A practical inventory control policy using operational statistics. \emph{Operations Research Letters} 33(4):341–348.

\bibitem[{Loke et~al.(2022)Loke, Tang, \protect\BIBand{} Xiao}]{loke2022decision}
Loke GG, Tang Q, Xiao Y (2022) Decision-driven regularization: A blended model for predict-then-optimize. \emph{Available at \url{https://dx.doi.org/10.2139/ssrn.3623006}} .

\bibitem[{Lu et~al.(2017)Lu, Pu, Wang, Hu, \protect\BIBand{} Wang}]{lu2017expressive}
Lu Z, Pu H, Wang F, Hu Z, Wang L (2017) The expressive power of neural networks: A view from the width. \emph{Advances in Neural Information Processing Systems} (Curran Associates, Inc.).

\bibitem[{Mandi et~al.(2022)Mandi, Bucarey, Tchomba, \protect\BIBand{} Guns}]{Mandi2022learning2rank}
Mandi J, Bucarey V, Tchomba MMK, Guns T (2022) Decision-focused learning: Through the lens of learning to rank. \emph{International Conference on Machine Learning}, 14935--14947 (PMLR).

\bibitem[{Mandi \protect\BIBand{} Guns(2020)}]{Mandi_Guns_2020}
Mandi J, Guns T (2020) Interior point solving for {LP}-based prediction+optimisation. \emph{Advances in Neural Information Processing Systems}, volume~33, 7272–7282 (Curran Associates, Inc.).

\bibitem[{Mandi et~al.(2023)Mandi, Kotary, Berden, Mulamba, Bucarey, Guns, \protect\BIBand{} Fioretto}]{mandi2023decision}
Mandi J, Kotary J, Berden S, Mulamba M, Bucarey V, Guns T, Fioretto F (2023) Decision-focused learning: Foundations, state of the art, benchmark and future opportunities. \emph{arXiv preprint arXiv:2307.13565} .

\bibitem[{Mandi et~al.(2020)Mandi, Stuckey, Guns et~al.}]{Mandi2020_spo}
Mandi J, Stuckey PJ, Guns T, et~al. (2020) Smart predict-and-optimize for hard combinatorial optimization problems. \emph{Proceedings of the AAAI Conference on Artificial Intelligence} 34(02):1603--1610.

\bibitem[{Mart{\'\i}nez-de Albeniz \protect\BIBand{} Belkaid(2021)}]{martinez2021comes}
Mart{\'\i}nez-de Albeniz V, Belkaid A (2021) Here comes the sun: Fashion goods retailing under weather fluctuations. \emph{European Journal of Operational Research} 294(3):820--830.

\bibitem[{Mate et~al.(2022)Mate, Madaan, Taneja, Madhiwalla, Verma, Singh, Hegde, Varakantham, \protect\BIBand{} Tambe}]{Mate_Madaan_Taneja_Madhiwalla_Verma_Singh_Hegde_Varakantham_Tambe_2022}
Mate A, Madaan L, Taneja A, Madhiwalla N, Verma S, Singh G, Hegde A, Varakantham P, Tambe M (2022) Field study in deploying restless multi-armed bandits: Assisting non-profits in improving maternal and child health. \emph{Proceedings of the AAAI Conference on Artificial Intelligence} 36(1111):12017–12025.

\bibitem[{McKenzie et~al.(2023)McKenzie, Fung, \protect\BIBand{} Heaton}]{McKenzie_Fung_Heaton_2023}
McKenzie D, Fung SW, Heaton H (2023) Faster predict-and-optimize with three-operator splitting. \emph{arXiv preprint arXiv:2301.13395} .

\bibitem[{Mieth et~al.(2023)Mieth, Morales, \protect\BIBand{} Poor}]{Mieth_Morales_Poor_2023}
Mieth R, Morales JM, Poor HV (2023) Data valuation from data-driven optimization. \emph{arXiv preprint arXiv:2305.01775} .

\bibitem[{Mi{\v{s}}i{\'c} \protect\BIBand{} Perakis(2020)}]{mivsic2020data}
Mi{\v{s}}i{\'c} VV, Perakis G (2020) Data analytics in operations management: A review. \emph{Manufacturing \& Service Operations Management} 22(1):158--169.

\bibitem[{Monga et~al.(2021)Monga, Li, \protect\BIBand{} Eldar}]{monga2021algorithm}
Monga V, Li Y, Eldar YC (2021) Algorithm unrolling: Interpretable, efficient deep learning for signal and image processing. \emph{IEEE Signal Processing Magazine} 38(2):18--44.

\bibitem[{Mulamba et~al.(2021)Mulamba, Mandi, Diligenti, Lombardi, Lopez, \protect\BIBand{} Guns}]{mulamba2021contrastive}
Mulamba M, Mandi J, Diligenti M, Lombardi M, Lopez VB, Guns T (2021) Contrastive losses and solution caching for predict-and-optimize. \emph{30th International Joint Conference on Artificial Intelligence (IJCAI-21): IJCAI-21}, 2833--2840 (International Joint Conferences on Artificial Intelligence).

\bibitem[{Mu{\~n}oz et~al.(2022)Mu{\~n}oz, Pineda, \protect\BIBand{} Morales}]{Munoz2022}
Mu{\~n}oz MA, Pineda S, Morales JM (2022) A bilevel framework for decision-making under uncertainty with contextual information. \emph{Omega} 108:102575.

\bibitem[{Nadaraya(1964)}]{nadaraya1964estimating}
Nadaraya E (1964) On estimating regression. \emph{Theory of Probability \& its Applications} 9(1):141--142.

\bibitem[{Neghab et~al.(2022)Neghab, Khayyati, \protect\BIBand{} Karaesmen}]{Neghab2022integrated}
Neghab DP, Khayyati S, Karaesmen F (2022) An integrated data-driven method using deep learning for a newsvendor problem with unobservable features. \emph{European Journal of Operational Research} 302(2):482--496.

\bibitem[{Ng et~al.(2000)Ng, Russell et~al.}]{ng2000algorithms}
Ng AY, Russell S, et~al. (2000) Algorithms for inverse reinforcement learning. \emph{Icml}, volume~1, 2.

\bibitem[{Nguyen et~al.(2020)Nguyen, Zhang, Blanchet, Delage, \protect\BIBand{} Ye}]{NEURIPS2020_adf854f4}
Nguyen VA, Zhang F, Blanchet J, Delage E, Ye Y (2020) Distributionally robust local non-parametric conditional estimation. \emph{Advances in Neural Information Processing Systems}, volume~33, 15232--15242 (Curran Associates, Inc.).

\bibitem[{Nguyen et~al.(2021)Nguyen, Zhang, Blanchet, Delage, \protect\BIBand{} Ye}]{Nguyen_Zhang_Blanchet_Delage_Ye_2022}
Nguyen VA, Zhang F, Blanchet J, Delage E, Ye Y (2021) Robustifying conditional portfolio decisions via optimal transport. \emph{arXiv preprint arXiv:2103.16451} .

\bibitem[{Nikishin et~al.(2022{\natexlab{a}})Nikishin, Abachi, Agarwal, \protect\BIBand{} Bacon}]{nikishin2022control}
Nikishin E, Abachi R, Agarwal R, Bacon PL (2022{\natexlab{a}}) Control-oriented model-based reinforcement learning with implicit differentiation. \emph{Proceedings of the AAAI Conference on Artificial Intelligence} 36(7):7886--7894.

\bibitem[{Nikishin et~al.(2022{\natexlab{b}})Nikishin, D'Oro, Precup, Barreto, massoud Farahmand, Bacon, \protect\BIBand{} Hall}]{NikishinEtAl22_workshop}
Nikishin E, D'Oro P, Precup D, Barreto A, massoud Farahmand A, Bacon PL, Hall G (2022{\natexlab{b}}) Decision awareness in reinforcement learning. Workshop abstract. International Conference on Machine Learning.

\bibitem[{Notz \protect\BIBand{} Pibernik(2022)}]{Notz2022}
Notz PM, Pibernik R (2022) Prescriptive analytics for flexible capacity management. \emph{Management Science} 68(3):1756--1775.

\bibitem[{Ohmori(2021)}]{ohmori2021predictive}
Ohmori S (2021) A predictive prescription using minimum volume k-nearest neighbor enclosing ellipsoid and robust optimization. \emph{Mathematics} 9(2):119.

\bibitem[{Oroojlooyjadid et~al.(2020)Oroojlooyjadid, Snyder, \protect\BIBand{} Tak{\'a}{\v{c}}}]{Oroojlooyjadid2020}
Oroojlooyjadid A, Snyder LV, Tak{\'a}{\v{c}} M (2020) Applying deep learning to the newsvendor problem. \emph{IISE Transactions} 52(4):444--463.

\bibitem[{Perakis et~al.(2023)Perakis, Sim, Tang, \protect\BIBand{} Xiong}]{Perakis_Sim_Tang_Xiong_2023}
Perakis G, Sim M, Tang Q, Xiong P (2023) Robust pricing and production with information partitioning and adaptation. \emph{Management Science} 69(3):1398–1419.

\bibitem[{Perrault et~al.(2020)Perrault, Wilder, Ewing, Mate, Dilkina, \protect\BIBand{} Tambe}]{perrault_end--end_2020}
Perrault A, Wilder B, Ewing E, Mate A, Dilkina B, Tambe M (2020) End-to-{End} {Game}-{Focused} {Learning} of {Adversary} {Behavior} in {Security} {Games}. \emph{Proceedings of the AAAI Conference on Artificial Intelligence} 34(02):1378--1386.

\bibitem[{Per{\v{s}}ak \protect\BIBand{} Anjos(2023)}]{Pervsak2023contextual}
Per{\v{s}}ak E, Anjos MF (2023) Contextual robust optimisation with uncertainty quantification. \emph{Proceedings of the 20th International Conference on the Integrationg of Constraint Programming, Artificial Intelligence, and Operations Research}, 124--132 (Springer).

\bibitem[{Pineda et~al.(2022)Pineda, Fan, Monge, Venkataraman, Sodhi, Chen, Ortiz, DeTone, Wang, Anderson et~al.}]{pineda2022theseus}
Pineda L, Fan T, Monge M, Venkataraman S, Sodhi P, Chen RT, Ortiz J, DeTone D, Wang A, Anderson S, et~al. (2022) Theseus: A library for differentiable nonlinear optimization. \emph{Advances in Neural Information Processing Systems}, volume~35, 3801--3818 (Curran Associates, Inc.).

\bibitem[{Qi \protect\BIBand{} Shen(2022)}]{qi2022integrating}
Qi M, Shen ZJ (2022) Integrating prediction/estimation and optimization with applications in operations management. \emph{Tutorials in Operations Research: Emerging and Impactful Topics in Operations}, 36--58 (INFORMS).

\bibitem[{Qi et~al.(2023)Qi, Shi, Qi, Ma, Yuan, Wu, \protect\BIBand{} Shen}]{Qi_MS2023}
Qi M, Shi Y, Qi Y, Ma C, Yuan R, Wu D, Shen ZJM (2023) A practical end-to-end inventory management model with deep learning. \emph{Management Science} 69(2):759--773.

\bibitem[{Rahimian \protect\BIBand{} Pagnoncelli(2022)}]{Rahimian2022}
Rahimian H, Pagnoncelli B (2022) Data-driven approximation of contextual chance-constrained stochastic programs. \emph{Available at \url{https://optimization-online.org/?p=20569}} .

\bibitem[{Ren et~al.(2022)Ren, Feng, Liu, Pan, Fu, Mai, \protect\BIBand{} Yang}]{ren2022torchopt}
Ren J, Feng X, Liu B, Pan X, Fu Y, Mai L, Yang Y (2022) Torchopt: An efficient library for differentiable optimization. \emph{arXiv preprint arXiv:2211.06934} .

\bibitem[{Rios et~al.(2015)Rios, Wets, \protect\BIBand{} Woodruff}]{Rios2015multi}
Rios I, Wets RJ, Woodruff DL (2015) Multi-period forecasting and scenario generation with limited data. \emph{Computational Management Science} 12(2):267--295.

\bibitem[{Rockafellar \protect\BIBand{} Wets(2009)}]{rockafellar1998variational}
Rockafellar RT, Wets RJB (2009) \emph{Variational Analysis} (Berlin: Springer).

\bibitem[{Rudin(2019)}]{Rudin2019}
Rudin C (2019) {Stop explaining black box machine learning models for high stakes decisions and use interpretable models instead}. \emph{Nature Machine Intelligence} 1(5):206--215.

\bibitem[{Rust(1987)}]{rust1987optimal}
Rust J (1987) Optimal replacement of gmc bus engines: An empirical model of harold zurcher. \emph{Econometrica: Journal of the Econometric Society} 999--1033.

\bibitem[{Rust(1988)}]{rust1998}
Rust J (1988) Maximum likelihood estimation of discrete control processes. \emph{SIAM Journal on Control and Optimization} 26(5):1006--1024.

\bibitem[{Rychener et~al.(2023)Rychener, Kuhn, \protect\BIBand{} Sutter}]{Rychener_Sutter_2023}
Rychener Y, Kuhn D, Sutter T (2023) End-to-end learning for stochastic optimization: A bayesian perspective. \emph{International Conference on Machine Learning}.

\bibitem[{Sahoo et~al.(2023)Sahoo, Paulus, Vlastelica, Musil, Kuleshov, \protect\BIBand{} Martius}]{Sahoo2023backpropagation}
Sahoo SS, Paulus A, Vlastelica M, Musil V, Kuleshov V, Martius G (2023) Backpropagation through combinatorial algorithms: Identity with projection works. \emph{International Conference on Learning Representations}.

\bibitem[{Sang et~al.(2022)Sang, Xu, Long, Hu, \protect\BIBand{} Sun}]{Sang_Xu_Long_Hu_Sun_2022}
Sang L, Xu Y, Long H, Hu Q, Sun H (2022) Electricity price prediction for energy storage system arbitrage: A decision-focused approach. \emph{IEEE Transactions on Smart Grid} 13(4):2822–2832.

\bibitem[{Sang et~al.(2023)Sang, Xu, Long, \protect\BIBand{} Wu}]{sang_2023}
Sang L, Xu Y, Long H, Wu W (2023) Safety-aware semi-end-to-end coordinated decision model for voltage regulation in active distribution network. \emph{IEEE Transactions on Smart Grid} 14(3):1814--1826.

\bibitem[{Sen \protect\BIBand{} Deng(2017)}]{sen2017learning}
Sen S, Deng Y (2017) Learning enabled optimization: Towards a fusion of statistical learning and stochastic programming. \emph{Available at \url{https://optimization-online.org/?p=14456}} .

\bibitem[{Serrano et~al.(2022)Serrano, Minner, Schiffer, \protect\BIBand{} Vidal}]{serrano_bilevel_2022}
Serrano B, Minner S, Schiffer M, Vidal T (2022) Bilevel optimization for feature selection in the data-driven newsvendor problem. \emph{arXiv preprint arXiv:2209.05093} .

\bibitem[{Shah et~al.(2022)Shah, Wang, Wilder, Perrault, \protect\BIBand{} Tambe}]{shah2022decision}
Shah S, Wang K, Wilder B, Perrault A, Tambe M (2022) Decision-focused learning without decision-making: Learning locally optimized decision losses. \emph{Advances in Neural Information Processing Systems}, volume~35, 1320--1332 (Curran Associates, Inc.).

\bibitem[{Shapiro et~al.(2014)Shapiro, Dentcheva, \protect\BIBand{} Ruszczyński}]{Shapiro2014_lectures}
Shapiro A, Dentcheva D, Ruszczyński A (2014) \emph{Lectures on Stochastic Programming: Modeling and Theory} (Philadelphia: SIAM), second edition.

\bibitem[{Shlezinger et~al.(2022)Shlezinger, Eldar, \protect\BIBand{} Boyd}]{Shlezinger_Eldar_Boyd_2022}
Shlezinger N, Eldar YC, Boyd SP (2022) Model-based deep learning: On the intersection of deep learning and optimization. \emph{IEEE Access} 10:115384–115398.

\bibitem[{Simonis et~al.(2014)Simonis, O’Sullivan, Mehta, Hurley, \protect\BIBand{} Cauwer}]{csplib059}
Simonis H, O’Sullivan B, Mehta D, Hurley B, Cauwer MD (2014) {CSPLib} problem 059: Energy-cost aware scheduling. \url{http://www.csplib.org/Problems/prob059}, accessed on June 14, 2023.

\bibitem[{Smith \protect\BIBand{} Winkler(2006)}]{smith2006optimizers}
Smith JE, Winkler RL (2006) The optimizer's curse: Skepticism and postdecision surprise in decision analysis. \emph{Management Science} 52(3):311--322.

\bibitem[{Srivastava et~al.(2021)Srivastava, Wang, Hanasusanto, \protect\BIBand{} Ho}]{Srivastava_Wang_Hanasusanto_Ho_2021}
Srivastava PR, Wang Y, Hanasusanto GA, Ho CP (2021) On data-driven prescriptive analytics with side information: A regularized {N}adaraya-{W}atson approach. \emph{arXiv preprint arXiv:2110.04855} .

\bibitem[{Stratigakos et~al.(2022)Stratigakos, Camal, Michiorri, \protect\BIBand{} Kariniotakis}]{Stratigakos_Camal_Michiorri_Kariniotakis_2022}
Stratigakos A, Camal S, Michiorri A, Kariniotakis G (2022) Prescriptive trees for integrated forecasting and optimization applied in trading of renewable energy. \emph{IEEE Transactions on Power Systems} 37(6):4696–4708.

\bibitem[{Sun et~al.(2023{\natexlab{a}})Sun, Liu, \protect\BIBand{} Li}]{Sun_marginlp_2023}
Sun C, Liu S, Li X (2023{\natexlab{a}}) Maximum optimality margin: A unified approach for contextual linear programming and inverse linear programming. \emph{arXiv preprint arXiv:2301.11260} .

\bibitem[{Sun et~al.(2023{\natexlab{b}})Sun, Shi, Wang, Tuan, Poor, \protect\BIBand{} Tao}]{Sun_Shi_Wang2023}
Sun H, Shi Y, Wang J, Tuan HD, Poor HV, Tao D (2023{\natexlab{b}}) Alternating differentiation for optimization layers. \emph{International Conference on Learning Representations}.

\bibitem[{Sun et~al.(2016)Sun, Li, Xu et~al.}]{sun2016deep}
Sun J, Li H, Xu Z, et~al. (2016) Deep {ADMM-Net} for compressive sensing {MRI}. \emph{Advances in Neural Information Processing Systems}, volume~29 (Curran Associates, Inc.).

\bibitem[{Sun et~al.(2023{\natexlab{c}})Sun, Leung, Li, \protect\BIBand{} Wu}]{Sun2023unified}
Sun X, Leung CH, Li Y, Wu Q (2023{\natexlab{c}}) A unified perspective on regularization and perturbation in differentiable subset selection. \emph{International Conference on Artificial Intelligence and Statistics}, 4629--4642 (PMLR).

\bibitem[{Sutton et~al.(1999)Sutton, McAllester, Singh, \protect\BIBand{} Mansour}]{Sutton1999}
Sutton RS, McAllester D, Singh S, Mansour Y (1999) Policy gradient methods for reinforcement learning with function approximation. \emph{Advances in Neural Information Processing Systems}, volume~12 (MIT Press).

\bibitem[{Tang \protect\BIBand{} Khalil(2022)}]{Tang_Khalil_2022}
Tang B, Khalil EB (2022) {PyEPO}: A {PyTorch}-based end-to-end predict-then-optimize library for linear and integer programming. \emph{arXiv preprint arXiv:2206.14234} .

\bibitem[{Tian et~al.(2023{\natexlab{a}})Tian, Yan, Liu, \protect\BIBand{} Wang}]{tian2023smart}
Tian X, Yan R, Liu Y, Wang S (2023{\natexlab{a}}) A smart predict-then-optimize method for targeted and cost-effective maritime transportation. \emph{Transportation Research Part B: Methodological} 172:32--52.

\bibitem[{Tian et~al.(2023{\natexlab{b}})Tian, Yan, Wang, \protect\BIBand{} Laporte}]{Tian_Yan_Wang_Laporte_2023}
Tian X, Yan R, Wang S, Laporte G (2023{\natexlab{b}}) Prescriptive analytics for a maritime routing problem. \emph{Ocean \& Coastal Management} 242:106695.

\bibitem[{Tian et~al.(2023{\natexlab{c}})Tian, Yan, Wang, Liu, \protect\BIBand{} Zhen}]{tian:tutorial}
Tian X, Yan R, Wang S, Liu Y, Zhen L (2023{\natexlab{c}}) Tutorial on prescriptive analytics for logistics: What to predict and how to predict. \emph{Electronic Research Archive} 31(4):2265--2285.

\bibitem[{Van~Parys et~al.(2021)Van~Parys, Esfahani, \protect\BIBand{} Kuhn}]{vanparys2021data}
Van~Parys BP, Esfahani PM, Kuhn D (2021) From data to decisions: Distributionally robust optimization is optimal. \emph{Management Science} 67(6):3387--3402.

\bibitem[{Vlastelica et~al.(2019)Vlastelica, Paulus, Musil, Martius, \protect\BIBand{} Rolinek}]{vlastelica2019diff}
Vlastelica M, Paulus A, Musil V, Martius G, Rolinek M (2019) Differentiation of blackbox combinatorial solvers. \emph{International Conference on Learning Representations}.

\bibitem[{Vohra et~al.(2023)Vohra, Rajaei, \protect\BIBand{} Cremer}]{Vohra2023end}
Vohra R, Rajaei A, Cremer JL (2023) End-to-end learning with multiple modalities for system-optimised renewables nowcasting. \emph{arXiv preprint arXiv:2304.07151} .

\bibitem[{Wager \protect\BIBand{} Athey(2018)}]{Wager2018estimation}
Wager S, Athey S (2018) Estimation and inference of heterogeneous treatment effects using random forests. \emph{Journal of the American Statistical Association} 113(523):1228--1242.

\bibitem[{Wahdany et~al.(2023)Wahdany, Schmitt, \protect\BIBand{} Cremer}]{Wahdany_Schmitt_Cremer_2023}
Wahdany D, Schmitt C, Cremer JL (2023) More than accuracy: end-to-end wind power forecasting that optimises the energy system. \emph{Electric Power Systems Research} 221:109384.

\bibitem[{Wang et~al.(2023)Wang, Verma, Mate, Shah, Taneja, Madhiwalla, Hegde, \protect\BIBand{} Tambe}]{wang_scalable_2023}
Wang K, Verma S, Mate A, Shah S, Taneja A, Madhiwalla N, Hegde A, Tambe M (2023) Scalable decision-focused learning in restless multi-armed bandits with application to maternal and child health. \emph{arXiv preprint arXiv:2202.00916} .

\bibitem[{Wang et~al.(2020)Wang, Wilder, Perrault, \protect\BIBand{} Tambe}]{wang_automatically_2020}
Wang K, Wilder B, Perrault A, Tambe M (2020) Automatically learning compact quality-aware surrogates for optimization problems. \emph{Advances in Neural Information Processing Systems}, volume~33, 9586--9596 (Curran Associates, Inc.).

\bibitem[{Wang et~al.(2021)Wang, Chen, \protect\BIBand{} Wang}]{wang2021distributionally}
Wang T, Chen N, Wang C (2021) Distributionally robust prescriptive analytics with {Wasserstein} distance. \emph{arXiv preprint arXiv:2106.05724} .

\bibitem[{Watson(1964)}]{watson1964smooth}
Watson G (1964) Smooth regression analysis. \emph{Sankhy\={a}: The Indian Journal of Statistics, Series A} 26(4):359--372.

\bibitem[{Wilder et~al.(2019{\natexlab{a}})Wilder, Dilkina, \protect\BIBand{} Tambe}]{Wilder_Dilkina_Tambe_2019}
Wilder B, Dilkina B, Tambe M (2019{\natexlab{a}}) Melding the data-decisions pipeline: Decision-focused learning for combinatorial optimization. \emph{Proceedings of the AAAI Conference on Artificial Intelligence} 33(01):1658–1665.

\bibitem[{Wilder et~al.(2019{\natexlab{b}})Wilder, Ewing, Dilkina, \protect\BIBand{} Tambe}]{Wilder2019graph}
Wilder B, Ewing E, Dilkina B, Tambe M (2019{\natexlab{b}}) End to end learning and optimization on graphs. \emph{Advances in Neural Information Processing Systems}, volume~32 (Curran Associates, Inc.).

\bibitem[{Xie et~al.(2019)Xie, Wu, Liu, Zhong, \protect\BIBand{} Lin}]{xie2019differentiable}
Xie X, Wu J, Liu G, Zhong Z, Lin Z (2019) Differentiable linearized {ADMM}. \emph{International Conference on Machine Learning}, 6902--6911 (PMLR).

\bibitem[{Xie et~al.(2020)Xie, Dai, Chen, Dai, Zhao, Zha, Wei, \protect\BIBand{} Pfister}]{Xie_Dai_Chen_Dai_Zhao_Zha_Wei_Pfister_2020}
Xie Y, Dai H, Chen M, Dai B, Zhao T, Zha H, Wei W, Pfister T (2020) Differentiable top-k with optimal transport. \emph{Advances in Neural Information Processing Systems}, volume~33, 20520–20531 (Curran Associates, Inc.).

\bibitem[{Xu \protect\BIBand{} Cohen(2018)}]{xu-cohen-2018-stock}
Xu Y, Cohen SB (2018) Stock movement prediction from tweets and historical prices. \emph{Proceedings of the 56th Annual Meeting of the Association for Computational Linguistics (Volume 1: Long Papers)}, 1970--1979 (Association for Computational Linguistics).

\bibitem[{Yan et~al.(2021)Yan, Wang, Cao, \protect\BIBand{} Sun}]{yan2021}
Yan R, Wang S, Cao J, Sun D (2021) Shipping domain knowledge informed prediction and optimization in port state control. \emph{Transportation Research Part B: Methodological} 149:52--78.

\bibitem[{Yan et~al.(2020)Yan, Wang, \protect\BIBand{} Fagerholt}]{Yan_Wang_Fagerholt_2020}
Yan R, Wang S, Fagerholt K (2020) A semi-“smart predict then optimize” (semi-{SPO}) method for efficient ship inspection. \emph{Transportation Research Part B: Methodological} 142:100–125.

\bibitem[{Yan et~al.(2023)Yan, Wang, \protect\BIBand{} Zhen}]{yan2023}
Yan R, Wang S, Zhen L (2023) An extended smart ``predict, and optimize'' (spo) framework based on similar sets for ship inspection planning. \emph{Transportation Research Part E: Logistics and Transportation Review} 173:103109.

\bibitem[{Yang et~al.(2023)Yang, Zhang, Chen, Gao, \protect\BIBand{} Hu}]{yang_shang2023}
Yang J, Zhang L, Chen N, Gao R, Hu M (2023) Decision-making with side information: A causal transport robust approach. \emph{Available at \url{https://optimization-online.org/?p=20639}} .

\bibitem[{Yan{\i}ko{\u{g}}lu et~al.(2019)Yan{\i}ko{\u{g}}lu, Gorissen, \protect\BIBand{} den Hertog}]{yanikoglu2019survey}
Yan{\i}ko{\u{g}}lu {\.I}, Gorissen BL, den Hertog D (2019) A survey of adjustable robust optimization. \emph{European Journal of Operational Research} 277(3):799--813.

\bibitem[{Zhang et~al.(2021)Zhang, Zhang, Cucuringu, \protect\BIBand{} Zohren}]{zhang2021universal}
Zhang C, Zhang Z, Cucuringu M, Zohren S (2021) A universal end-to-end approach to portfolio optimization via deep learning. \emph{arXiv preprint arXiv:2111.09170} .

\bibitem[{Zhang et~al.(2023{\natexlab{a}})Zhang, Yang, \protect\BIBand{} Gao}]{Zhang_Yang_Gao}
Zhang L, Yang J, Gao R (2023{\natexlab{a}}) Optimal robust policy for feature-based newsvendor. \emph{Management Science} Forthcoming.

\bibitem[{Zhang \protect\BIBand{} Gao(2017)}]{ZhangNNNewsvendor}
Zhang Y, Gao J (2017) Assessing the performance of deep learning algorithms for newsvendor problem. Liu D, Xie S, Li Y, Zhao D, El-Alfy ESM, eds., \emph{Neural Information Processing}, 912--921 (Cham: Springer International Publishing).

\bibitem[{Zhang et~al.(2023{\natexlab{b}})Zhang, Liu, \protect\BIBand{} Zhao}]{zhang2023_pwaffine}
Zhang Y, Liu J, Zhao X (2023{\natexlab{b}}) Data-driven piecewise affine decision rules for stochastic programming with covariate information. \emph{arXiv preprint arXiv:2304.13646} .

\bibitem[{Zhu et~al.(2022)Zhu, Xie, \protect\BIBand{} Sim}]{Zhu_Xie_Sim_2022}
Zhu T, Xie J, Sim M (2022) Joint estimation and robustness optimization. \emph{Management Science} 68(3):1659–1677.

\bibitem[{Ziebart et~al.(2008)Ziebart, Maas, Bagnell, \protect\BIBand{} Dey}]{ziebart2008maximum}
Ziebart BD, Maas AL, Bagnell JA, Dey AK (2008) Maximum entropy inverse reinforcement learning. \emph{Proceedings of the AAAI Conference on Artificial Intelligence} 8:1433--1438.

\end{thebibliography}

\clearpage

\begin{APPENDICES}

\section{Theoretical guarantees}
\label{sec:appen:definitions}
In this section, we provide some definitions and theoretical results for completeness.
\begin{definition}\citep{Bertsimas_Kallus_2020}
\label{def:optimal}
    We say that a policy $\pi^N(\x)$ obtained using $N$ samples is asymptotically optimal if, with probability 1, we have that for $\Prob(\x)$-almost-everywhere $\x\in \X$:
  \[\lim_{N\rightarrow \infty} h(\pi^N(\x),\Prob(\y \rvert \x)) = h(\pi^*(\x),\Prob(\y \rvert \x)).\]
\end{definition}
\begin{definition}\citep{Bertsimas_Kallus_2020}
\label{def:const}
       We say that a policy $\pi^N(\x)$ obtained using $N$ samples is consistent if, with probability 1, we have that for $\Prob(\x)$-almost-everywhere $\x\in \X$:
\[\lVert \pi^N(\x)-\Z^*(\x)\rVert =0 \;\;\;\text{where}\; \lVert \pi^N(\x)-\Z^*(\x)\rVert = \inf_{\z\in \Z^*(\x)} \lVert \pi^N(\x)-\z\rVert,\]
and
\[\Z^*(\x) =\{\z \rvert \z \in \argmin_{\z' \in \Z}\Expect_{\Prob} \big[h(\z', \Prob(\y \rvert \x))\big].\]
\end{definition}

The above conditions imply that, as the number of samples tends to infinity, the performance of the \revised{policy} under almost all covariates matches the optimal conditional cost.
\begin{definition} \label{def:fisher}
    A surrogate loss function $\ell$ is Fisher consistent with respect to the SPO  loss if the following condition holds for all $\x\in \X$:
    \[ \argmin_{\btheta} 
 \Expect_{\Prob(\y \rvert \x)}[\ell (g_\btheta(\x),\y)]\subseteq \argmin_{\btheta}   \Expect_{\Prob(\y \rvert \x)} [\ell_{\model{SPO}} (g_\btheta(\x),\y)].
\]
\end{definition}
The Fisher consistency condition defined above requires complete knowledge of the joint distribution $\Prob$. Instead, an interesting issue is to determine if the surrogate loss $\ell$ is calibrated with respect to \model{SPO}, that is,  whether low surrogate excess risk translates to small excess true risk.
\begin{definition}\citep{ho-nguyen_risk_2022} \label{def:uniform_calibration}
    A loss function $\ell$ is uniformly calibrated with respect to \model{SPO} for a class of distributions $\mathcal{D}$ on $\X \times \mathcal{Y}$ if for $\epsilon>0$,  there exists a function $\Delta_\ell (\cdot): \Re^+ \rightarrow \Re^+$ such that for all $x \in \X$:
    \begin{align*}
        &\Expect_{\Prob(\y \rvert \x)}[\ell (g_\btheta(\x),\y)] - \inf_{\btheta'}\Expect_{\Prob(\y \rvert \x)}[\ell (g_{\btheta'}(\x),\y)] < \Delta_{\ell}(\epsilon) \\&\quad \Rightarrow \Expect_{\Prob(\y \rvert \x)}[\ell_{\model{SPO}} (g_\btheta(\x),\y)] - \inf_{\btheta'}\Expect_{\Prob(\y \rvert \x)}[\ell_{\model{SPO}} (g_{\btheta'}(\x),\y)] < \epsilon.
    \end{align*}
\end{definition}

\cite{ho-nguyen_risk_2022} introduce a \quoteIt{calibration function} with which it is simpler to verify the uniform calibration of MSE loss and show that the calibration function is $\text{O}(\epsilon^2)$. \cite{liu2021risk} assumed that the conditional distribution of $\y$ given $\x$ is bounded from below by the density of a normal distribution and obtain a $\text{O}(\epsilon^2)$ calibration function for polyhedral $\Z$ and $\text{O}(\epsilon)$ when $\mathcal{Z}$ is a level set of a strongly convex and smooth function. 

\begin{definition}\label{def:rademacher}
The empirical Rademacher complexity of a hypothesis class $\mathcal{H}$ under the loss function $\ell$ is given by:
\begin{equation*} 
    \mathbb{E}_{\sigma}\Bigg[\sup_{g \in \mathcal{H}}\frac{1}{N}\Bigg|\sum_{i=1}^{N}\sigma_i \ell(g(\x_i),\y_i)\Bigg|\Bigg],
\end{equation*}
where $\sigma_1, \sigma_2, \cdots, \sigma_N$ are independent and identically distributed Rademacher random variables, i.e., $\Prob(\sigma_i=1) = \Prob(\sigma_i=-1)=\frac{1}{2}$, $\forall i \in \{1,2,\cdots, N\}$.
\end{definition}
 
\section{List of abbreviations}\label{appen:abbrev}
In Table~\ref{appen:abbrev}, we provide the list of abbreviations used in this survey.
\begin{table}[ht]
\centering
\caption{List of abbreviations}
\label{tab:abbrevs}
\begin{tabular}{ll} 
    \toprule
    Abbreviation & Description \\
    \midrule
    ADMM & alternating direction method of multipliers \\
    CSO & \revised{contextual} stochastic optimization \\
    CVaR & conditional value at risk \\
    DRO & distributionally robust optimization \\
    DNN & deep neural network \\
    ERM & empirical risk minimization \\
    FP & fixed point \\
    ILO & integrated learning and optimization \\
    IFT & implicit function theorem \\
    kNN & $k$-nearest neighbor\\
    KKT & Karush–Kuhn–Tucker \\
    KL & Kullback-Leibler \\
    LDR & linear decision rule \\
    LP & linear program\\
    \revised{MDP} & \revised{Markov decision process}\\
    ML & machine learning \\
    MLE & maximum likelihood estimation \\
    MILP & mixed-integer linear program\\
    NW & Nadaraya-Watson \\
    RKHS & reproducing kernel Hilbert space \\
    SAA & sample average approximation \\
    SLO & sequential learning and optimization \\
    QP & quadratic program \\
    \bottomrule
    \end{tabular}
\end{table}

\end{APPENDICES}

\end{document}